\documentclass[twoside]{amsart}
\usepackage{amsmath,amssymb,amsfonts}
\usepackage[bookmarksnumbered,plainpages,hypertex]{hyperref}

\newtheorem{theorem}{\sc Theorem}[section]
\newtheorem{proposition}[theorem]{\sc Proposition}

\newtheorem{lemma}[theorem]{\sc Lemma}
\newtheorem{corollary}[theorem]{\sc Corollary}
\theoremstyle{definition}
\newtheorem{definition}[theorem]{\sc Definition}
\newtheorem{definitions}[theorem]{\sc Definitions}
\newtheorem{example}[theorem]{\sc Example}

\newtheorem{problem}[theorem]{\sc Problem}
\theoremstyle{remark}
\newtheorem{remark}[theorem]{\sc Remark}

\newtheorem{claim}[theorem]{}
\newtheorem{noname}[theorem]{}

\begin{document}
\title{A First Sight Towards Primitively Generated Connected Braided
Bialgebras}
\author{Alessandro Ardizzoni}
\address{University of Ferrara, Department of Mathematics, Via Machiavelli
35, Ferrara, I-44100, Italy}
\email{alessandro.ardizzoni@unife.it}
\urladdr{http://www.unife.it/utenti/alessandro.ardizzoni}
\thanks{This paper was written while the author was member of GNSAGA with
partial financial support from MIUR}

\begin{abstract}
The main aim of this paper is to investigate the structure of primitively
generated connected braided bialgebras $A$ with respect to the braided
vector space $P$ consisting of their primitive elements. When the Nichols
algebra of $P$ is obtained dividing out the tensor algebra $T( P) $ by the
two-sided ideal generated by its primitive elements of degree at least two,
we show that $A$ can be recovered as a sort of universal enveloping algebra
of $P$. One of the main applications of our construction is the description,
in terms of universal enveloping algebras, of connected braided bialgebras
whose associated graded coalgebra is a quadratic algebra.
\end{abstract}

\maketitle
\tableofcontents

\pagestyle{headings}

\section{Introduction}

Let $K$ be a fixed field, let $H$ be a pointed Hopf algebra over $K$ (this
means that all its simple subcoalgebras are one dimensional). Denote by $G$
the set of grouplike elements in $H$. It is well known that the graded
coalgebra $\mathrm{gr}\left( H\right) ,$ associated to the coradical
filtration of $H,$ is a Hopf algebra itself and can be described as a
Radford-Majid bosonization by $KG$ of a suitable graded connected braided
bialgebra $R,$ called diagram of $H,$ in the braided monoidal category of
Yetter-Drinfeld modules over $KG$. This is the starting point of the so
called \emph{lifting\ method} for the classification of finite dimensional
pointed Hopf algebras, introduced by N. Andruskiewitsch and H.J. Schneider
(see e.g. \cite{AS- Lifting}). Accordingly to this method, first one has to
describe $R$ by generators and relations, then to lift the informations
obtained to $H.$ It is worth noticing that in the finite dimensional case
there is a conjecture asserting that, in characteristic zero, $R$ is always
primitively generated \cite[Conjecture 2.7]{AS} (see also \cite[Conjecture
1.4]{AS- FiniteQuantCartan}).

Therefore, in many cases, for proving certain properties of Hopf algebras,
it is enough to do it in the connected case. The price that one has to pay
is to work with Hopf algebras in a braided category, and not with ordinary
Hopf algebras. Actually nowadays people recognize that it is more convenient
to work with braided bialgebras, that were introduced in \cite{Ta}.

Motivated by these observations, in this paper we will investigate the
structure of primitively generated connected braided bialgebras. Let $\left(
A,c_{A}\right) $ be such a bialgebra. Note that $\left( A,c_{A}\right) $ is
in particular a braided vector space i.e. $c_{A}:A\otimes A\rightarrow
A\otimes A$ is a map, called braiding, fulfilling the quantum Yang-Baxter
equation (\ref{ec: braided equation}). Furthermore $c_{A}$ induces a
braiding $c_{P}:P\otimes P\rightarrow P\otimes P$ on the space $P$ of
primitive elements in $A$. The tensor algebra $T\left( P\right) $ carries a
braided bialgebra structure depending on $c_{P},$ that will be denoted by $%
T\left( P,c_{P}\right) ,$ and there is a unique surjective braided bialgebra
map $\varphi :T\left( P,c_{P}\right) \rightarrow A$ that lifts the inclusion 
$P\subseteq A.$ Let $E\left( P,c_{P}\right) $ be the space of primitive
elements in $T\left( P,c_{P}\right) $ of degree at least two. Since $\varphi 
$ preserves primitive elements, we can define the map $b_{P}:E\left(
P,c_{P}\right) \rightarrow P,b\left( z\right) =\varphi \left( z\right) \ $%
and consider the quotient 
\begin{equation*}
U\left( P,c_{P},b_{P}\right) :=\frac{T\left( P,c_{P}\right) }{\left( \left( 
\mathrm{Id}-b_{P}\right) \left[ E\left( P,c_{P}\right) \right] \right) }.
\end{equation*}%
Now $\varphi $ quotients to a surjective braided bialgebra map $\overline{%
\varphi }:U\left( P,c_{P},b_{P}\right) \rightarrow A.$ It is remarkable here
that the canonical $K$-linear map $i_{U}:P\rightarrow U\left(
P,c_{P},b_{P}\right) $ is injective as $\overline{\varphi }\circ i_{U}$ is
the inclusion $P\subseteq A.$ It is also notable that $E\left(
P,c_{P}\right) $ may contain homogeneous elements of arbitrary degree.

By a famous result due to Heyneman and Radford (see \cite[Theorem 5.3.1]%
{Montgomery}), $\overline{\varphi }$ is an isomorphism if the space of
primitive elements in $U\left( P,c_{P},b_{P}\right) $ identifies with $P$
via $i_{U}$ (the converse is trivial).

We prove that this property holds whenever $P$ belongs to a large class $%
\mathcal{S}$ of braided vector spaces. It stems from our construction that $%
\mathcal{S}$ can be taken to be the class of braided vector spaces $(V,c)$
such that the Nichols algebra $\emph{B}\left( V,c\right) $, i.e. the image
of the canonical graded braided bialgebra homomorphism from the tensor
algebra $T\left( V,c\right) $ into the quantum shuffle algebra $T^{c}\left(
V,c\right) $, is obtained dividing out $T\left( V,c\right) $ by the
two-sided ideal generated by $E\left( V,c\right) $. We point out that the
class $\mathcal{S}$ is so large to enclose both the class of all braided
vector spaces of diagonal type whose Nichols algebra is a domain of finite
Gelfand-Kirillov dimension and also the class of all two dimensional braided
vector spaces of abelian group type whose symmetric algebra has dimension at
most $31.$ We note that in either one of these classes there are braided
vector spaces $(V,c)$ whose Nichols algebra is not quadratic and whose
braiding has minimal polynomial of degree greater than two.

Meaningful examples of elements in $\mathcal{S}$ are given. For instance,
any braided vector space with braiding of Hecke type with regular mark (e.g.
the braiding is a symmetry and the characteristic of $K$ is zero) is in $%
\mathcal{S}$.

An application of our construction is that $A\simeq U\left(
P,c_{P},b_{P}\right) $ whenever $A$ is a connected braided bialgebra such
that $\mathrm{gr}\left( A\right) $ is a quadratic algebra with respect to
its natural braided bialgebra structure.

We point out that the structure and properties of $U\left(
P,c_{P},b_{P}\right) $ are encoded in the datum $\left( P,c_{P},b_{P}\right)
.$ Indeed this leads to the introduction of what will be called a braided
Lie algebra $\left( V,c,b\right) $ for any braided vector space $\left(
V,c\right) $ and of the related universal enveloping algebra $U\left(
V,c,b\right) $. When $c$ is a symmetry, i.e. $c^{2}=\mathrm{Id}_{V\otimes
V}, $ and the characteristic of $K$ is zero, our enveloping algebra reduces
to the one introduced in \cite{Gu- Gen Trans Lie}. In this case, the crucial property, mentioned above, that the space of primitive elements in $U\left(V,c,b\right) $ identifies with $V$ via $i_{U}$ was proved in \cite[Lemma 6.2]{Kharchenko-Connected} and it was applied to obtain a Cartier-Kostant-Milnor-Moore type theorem for connected braided Hopf algebras with involutive braidings \cite[Theorem 6.1]{Kharchenko-Connected}. Other notions of Lie
algebra and enveloping algebra, extending the ones in \cite{Gu- Gen Trans Lie} to the non
symmetric case, appeared in the literature. Let us mention some of them.

\begin{itemize}
\item Lie algebras for braided vector spaces $\left( V,c\right) $ where $c$
is a braiding of Hecke type i.e. $\left( c+\mathrm{Id}_{V\otimes V}\right)
\left( c-q\mathrm{Id}_{V\otimes V}\right) =0$ for some $q\in K\backslash
\left\{ 0\right\} $ which is called the mark of $c$ \cite[Definition 7.1]{Wa}. See \cite[Theorem 5.5]{AMS-MM2} for a strong version of Cartier-Kostant-Milnor-Moore theorem in this setting.

\item Lie algebras for braided vector spaces $\left( V,c\right) $ where $c$
is not of Hecke-type but constructed by means of braidings of Hecke type 
\cite[Definition 1]{Gurevich: Hecke sym}. There a suitable quadratic algebra
is requiring to be Koszul. 

\item Lie algebras for objects in the braided monoidal category of
Yetter-Drinfeld modules over a Hopf algebra with bijective antipode \cite[%
Definition 4.1]{Pareigis-On Lie} (here called Pareigis-Lie algebras - see
Definition \ref{def: Pareigis-Lie alg}). In the present paper we also stress
the relation between our notion of Lie algebra and the one by Pareigis,
which involves partially defined $n$-ary bracket multiplications as well.

\item Lie algebras defined by considering quantum operations (see \cite[%
Definition 2.2]{Kharchenko-AnAlgSkew}) as primitive polynomials in the
tensor algebra. When the underline braided vector space is an object in the
category of Yetter-Drinfeld modules over some group algebra, to any Lie
algebra of this kind, the author associates a universal enveloping algebra
which is not connected (see \cite{Kharchenko-SkewPrim}). Note that the
universal enveloping algebra we introduce in the present paper is connected
as the scope here is to investigate the structure of primitively generated
connected braided bialgebras over $K$.
\end{itemize}

Finally, many authors have tried to characterize quantized enveloping
algebras $U_{q}(\mathfrak{g})$ of Drinfeld and Jimbo or other quantum groups
as bialgebras generated by a finite-dimensional Lie algebra like object via
some kind of enveloping algebra construction (see the introduction of \cite%
{Majid-Quantum}). One idea is to attempt to build this on $\mathfrak{g}$
itself but with some kind of deformed bracket obeying suitable new axioms.
Alternatively one can consider some subspace of $U_{q}(\mathfrak{g})$
endowed with some kind of `quantum Lie bracket' based on the quantum adjoint
action. Contributions to these problems can be found e.g. in \cite%
{Woronowicz, Lyubashenki-Sudbery, Delius-Gardner-Gould, Gomez-Majid,
Gurevich-Saponov} and in the references therein.

Now, since $U_{q}(\mathfrak{g})$ is not connected but just pointed, we have
no hope to describe it directly as an enveloping algebra of our kind (see 
\cite[Section 2]{Kharchenko-ACombinatorial} for a different approach).
Furthermore, it is remarkable that, by \cite[Theorem 7.8]{Masuoka-Abelian}, $%
U_{q}(\mathfrak{g}) $ can be viewed as a cocycle deformation of $\mathrm{gr}%
\left( U_{q}(\mathfrak{g})\right) $. By the foregoing, $\mathrm{gr}\left(
U_{q}(\mathfrak{g})\right) $ is the Radford-Majid bosonization of the
diagram $Q$ of $U_{q}(\mathfrak{g})$ by $KG,$ where $G$ denotes the set of
grouplike elements in $U_{q}(\mathfrak{g}).$ One could believe that $Q$ is
an enveloping algebra of our kind as it is a primitively generated connected
braided bialgebra. The point is that $Q$ is much more. It is a Nichols
algebra whence defined by quotienting the tensor algebra $T\left( P\left(
Q\right) \right) $ by homogeneous relations. Thus no enveloping algebra is
required to describe $Q$.\medskip \newline
\textbf{Methodology.} We proceed as follows. Let $(V,c)$ be a braided vector
space and denote by $E\left( V,c\right) $ the space of primitive elements in 
$T\left( V,c\right) $ of degree at least two. A bracket on $\left(
V,c\right) $ is a $K$-linear map $b:E\left( V,c\right) \rightarrow V$ that
commutes with the braiding of $T(V,c)$ in the sense that (\ref{form:
c-bracket}) is satisfied. If $b$ is a bracket on $\left( V,c\right) ,$ then
the \emph{universal enveloping algebra} of $\left( V,c\right) $ is defined
to be 
\begin{equation*}
U\left( V,c,b\right) :=\frac{T\left( V,c\right) }{\left( \left( \mathrm{Id}%
-b\right) \left[ E\left( V,c\right) \right] \right) }.
\end{equation*}%
Thus $U\left( V,c,b\right) $ carries a unique braided bialgebra structure
such that the canonical projection $\pi _{U}:T\left( V,c\right) \rightarrow
U\left( V,c,b\right) $ is a braided bialgebra homomorphism (Theorem \ref%
{teo: U bialgebra}). We say that $\left( V,c,b\right) $ is a \emph{braided
Lie algebra} whenever $\left( V,c\right) $ is a braided vector space, $%
b:E\left( V,c\right) \rightarrow V$ is a bracket on $\left( V,c\right) $ and
the canonical $K$-linear map $i_{U}:V\rightarrow U\left( V,c,b\right) $ is
injective (Definition \ref{def: c-Lie alg}). Natural examples of braided Lie
algebras arise as follows:

\begin{enumerate}
\item[1)] If $\left( V,c\right) $ is a braided vector space, then $\left(
V,c,0\right) $ is a braided Lie algebra (Proposition \ref{pro: (V,c,0) Lie}%
). The corresponding universal enveloping algebra is the \emph{symmetric
algebra} $S\left( V,c\right) .$

\item[2)] If $\left( A,c_{A}\right) $ is a braided algebra, then $\left(
A,c_{A},b_{A}\right) $ is a braided Lie algebra, where $b_{A}\ $acts on the $%
t$-th graded component of $E\left( A,c_{A}\right) $ as the restriction of
the iterated multiplication of $A\ $(Proposition \ref{pro: alg is Lie}).

\item[3)] If $\left( A,c_{A}\right) $ is a connected braided bialgebra and $%
P=P\left( A\right) $ is the space of primitive elements of $A,$ then $P$
forms a braided Lie algebra $\left( P,c_{P},b_{P}\right) $ with structures
induced by $\left( A,c_{A},b_{A}\right) $ (Theorem \ref{teo: univ U}). $%
\left( P,c_{P},b_{P}\right) $ will be called the \emph{infinitesimal braided
Lie algebra of }$A$.
\end{enumerate}

Theorem \ref{teo: magnum} contains the main characterization of the above
notion of braided Lie algebra. Let $\left( V,c\right) \in \mathcal{S}$ and
let $b:E\left( V,c\right) \rightarrow V$ be a bracket. Then, the following
assertions are equivalent.

\begin{itemize}
\item $\left( V,c,b\right) $ is a braided Lie algebra.

\item $U\left( V,c,b\right) $ is of PBW type in the sense of Definition \ref%
{def: PBW}.

\item $i_{U}:V\rightarrow U\left( V,c,b\right) $ induces an isomorphism
between $V$ and $P\left( U\left( V,c,b\right) \right) .$
\end{itemize}


Therefore, if $\left( V,c,b\right) $ is a braided Lie algebra whose
underlying braided vector space lies in $\mathcal{S}$, then $P\left( U\left(
V,c,b\right) \right) $ identifies with $V$ via $i_{U}$. The main consequence
of this property, which really justifies our construction, is Theorem \ref%
{teo: generated} concerning primitively generated connected braided
bialgebras $\left( A,c_{A}\right) .$ If $\left( P,c_{P}\right) \in \mathcal{S%
},$ then $A$ is isomorphic to $U\left( P,c_{P},b_{P}\right) $ as a braided
bialgebra, where $\left( P,c_{P},b_{P}\right) $ is the infinitesimal braided
Lie algebra of $A$. The proof of this fact uses the universal property of
the universal enveloping algebra (Theorem \ref{teo: univ U}) to obtain a
bialgebra projection of $U\left( P,c_{P},b_{P}\right) $ onto $A$; this
projection comes out to be injective as $P\left( U\left( V,c,b\right)
\right) $ identifies with $V$ via $i_{U}$.

In Section \ref{section: some braided}, the first meaningful examples of
elements in $\mathcal{S}$ are given. Between them, braided vector spaces of
diagonal type whose Nichols algebra is a domain of finite Gelfand-Kirillov
dimension (Theorem \ref{teo: GK dim}) and two dimensional braided vector
spaces of abelian group type whose symmetric algebra has dimension at most $%
31$ (Proposition \ref{pro: abelian grp}). In Examples \ref{es: AG}, \ref{ex:
twodim sdeg 2} and \ref{ex: Cartan}, we provide instances of braided vector
spaces which are not in $\mathcal{S}$ (the last two examples are due to V.
K. Kharchenko).

The main applications and examples we are interested in are given in the
last Sections \ref{section: Hecke}, \ref{section: Quadratic} and \ref%
{section: Pareigis}. Explicitly, in Section \ref{section: Hecke} we deal
with the class of braided vector spaces $\left( V,c\right) $ such that $c$
is a braiding of Hecke type i.e. it satisfies the equation $\left( c+\mathrm{%
Id}_{V\otimes V}\right) \left( c-q\mathrm{Id}_{V\otimes V}\right) =0$ for
some regular $q\in K\backslash \left\{ 0\right\} $ (see Definition \ref{def:
regular element}) called mark of $c$. In Theorem \ref{teo: sdeg Hecke} we
prove that this class is contained in $\mathcal{S}$. Let $\left(
V,c,b\right) $ be a braided Lie algebra such that $c$ is of Hecke type with
regular mark $q$. Then, in Theorem \ref{teo: Hecke Lie}, using the results
in \cite{AMS-MM2, AMS-MM}, we show that the universal enveloping
algebra of $\left( V,c,b\right) $ can be rewritten as follows 
\begin{equation*}
U\left( V,c,b\right) =\frac{T\left( V,c\right) }{\left( c\left( z\right) -qz-%
\left[ z\right] _{b}\mid z\in V\otimes V\right) }
\end{equation*}%
where $\left[ -\right] _{b}:V\otimes V\rightarrow V$ is given by $\left[ z%
\right] _{b}=b\left( c\left( z\right) -qz\right) $. Moreover if $q\neq 1$
and $\mathrm{char}\left( K\right) \neq 2,$ then $\left[ -\right] _{b}=0.$

Section \ref{section: Quadratic}, is devoted to study the class of braided
vector spaces $\left( V,c\right) $ such that the Nichols algebra $\emph{B}%
\left( V,c\right) $ is a quadratic algebra (see Definition \ref{def:
quadratic}). This also comes out to be a subclass of $\mathcal{S}$ (Theorem %
\ref{teo: quadratic => sdeg 1}). In the literature there is plenty of
examples of braided vector spaces of this kind (see Remark \ref{rem: B(V)
quadratic}). The main result of the section is Theorem \ref{teo: quadratic}
concerning connected braided bialgebras $\left( A,c_{A}\right) $ such that
the graded coalgebra $\mathrm{gr}\left( A\right) $ associated to the
coradical filtration of $A$ is a quadratic algebra with respect to its
natural braided bialgebra structure. Such an $A$ is proved to be isomorphic
as a braided bialgebra to the enveloping algebra of its infinitesimal
braided Lie algebra. The proof relies on the fact that the Nichols algebra
associated to this braided Lie algebra results to be quadratic. Braided Lie
algebras with this property will be further investigated in \cite{AS:
Ardi-Stumbo}.

In Section \ref{section: Pareigis}, we stress the relation between our
notions of Lie algebra and universal enveloping algebra, and those
introduced by Pareigis in \cite{Pareigis-On Lie}. The main result of the
section is Theorem \ref{teo: BrLie => PLie} establishing that a Pareigis-Lie
algebra $\left( V,c,\left[ -\right] \right) $ is associated to any braided
Lie algebra $\left( V,c,b\right) .$ Moreover there is a canonical braided
bialgebra projection $p:U_{P}\left( V,c,\left[ -\right] \right) \rightarrow
U\left( V,c,b\right) $ where $U_{P}\left( V,c,\left[ -\right] \right) $
denotes the universal enveloping algebra of $\left( V,c,\left[ -\right]
\right) $ (see Definition \ref{def: Pareigis-Lie alg}). We prove this
projection to be an isomorphism if condition (\ref{form: Pi su}) holds.

Note that, if $\left( V,c\right) $ lies in $\mathcal{S}$, then $p$ is an
isomorphism if and only if $P\left( U_{P}\left( V,c,\left[ -\right] \right)
\right) $ identifies with $V$ through the canonical map $i_{U_{P}}:V%
\rightarrow U_{P}\left( V,c,\left[ -\right] \right) $ (see Remark \ref{rem:
Pareigig prim}). \medskip \newline
\textbf{Strongness degree. }The class $\mathcal{S}$ can be described as the
class of braided vector spaces of strongness degree at most one. This
terminology requires some explanation. First, we point out that many of the
constructions and results we mentioned above concerning the tensor algebra $%
T\left( V,c\right) $ are indeed stated in the paper for a generic graded
braided bialgebra $B.$ In particular, in Definition \ref{def: symmetric}, we
associate a symmetric algebra $S\left( B\right) $ to any graded braided
bialgebra $B$. If $B$ is also endowed with a bracket $b$, then we can define
the universal enveloping algebra $U\left( B,b\right) $ (see Definition \ref%
{def: bracket}). Recall that a graded coalgebra $(C=\oplus _{n\in \mathbb{N}%
}C^{n},\Delta _{C},\varepsilon _{C})$ is \emph{strongly }$\mathbb{N}$\emph{%
-graded} whenever the $(i,j)$-homogeneous component $\Delta
_{C}^{i,j}:C^{i+j}\rightarrow C^{i}\otimes C^{j}$ of $\Delta _{C}$ is a
monomorphism for every $i,j\in \mathbb{N}$ (dually the notion of strongly $%
\mathbb{N}$-graded algebra can be introduced). Let $B$ be a graded braided
bialgebra and consider $B,S\left( B\right) ,S\left( S\left( B\right) \right) 
$ and so on. This yields a direct system whose direct limit is a graded
braided bialgebra which is strongly $\mathbb{N}$-graded as a coalgebra
(Theorem \ref{teo: S[inf] is strong}) and that coincides with $B^{0}\left[
B^{1}\right] ,$ the braided bialgebra of Type one associated to $B^{0}$ and $%
B^{1}$ (see \ref{claim: type one}), whenever $B$ is also strongly $\mathbb{N}
$-graded as an algebra (Corollary \ref{claim: type one}). We say that $B$%
\textbf{\ }\emph{has strongness degree }$n$ ($\mathrm{sdeg}\left( B\right)
=n $) if the direct system above is stationary exactly after $n$ steps (see
Definition \ref{def: strongness degree}). V. K. Kharchenko pointed out to
our attention that the notion of \emph{combinatorial rank} in \cite[%
Definition 5.4]{Kharchenko-SkewPrim} is essentially the same notion as the
strongness degree. Nevertheless we decided here to keep our terminology as
strongness degree is a measure of how far $B$ is to be strongly $%
\mathbb{N}
$-graded as a coalgebra. In Section \ref{sec: strongness degree} we also
provide criteria to compute $\mathrm{sdeg}\left( B\right) $ in some cases.
In appendix we collect some further results on strongness degree of graded
braided bialgebras. In particular we see that under suitable assumptions the
strongness degree of a graded braided bialgebra has an upper bound.

When we focus our attention on the case $B:=T\left( V,c\right) ,$ the tensor
algebra of a braided vector space $\left( V,c\right) $, we get $\emph{B}%
\left( V,c\right) =B^{0}\left[ B^{1}\right] ,S\left( V,c\right) =S\left(
B\right) ,U\left( V,c,b\right) =U\left( B,b\right) .$ The strongness degree
of a braided vector space is defined as the strongness degree of the
corresponding tensor algebra. Hence a braided vector space $\left(
V,c\right) $ is in $\mathcal{S}$ if and only if it has strongness degree at
most one.

Elsewhere, in the spirit of \cite[Section 6]{Kharchenko-SkewPrim}, we will
investigate a notion of Lie algebra and enveloping algebra for braided
vector spaces of arbitrary strongness degree.

\section{Preliminaries}

Throughout this paper $K$ will denote a field. All vector spaces will be
defined over $K$ and the tensor product over $K$ will be denoted by $\otimes$%
. \medskip\newline
In this section we define the main notions that we will deal with in the
paper.

\begin{definition}
Let $V$ be a vector space over a field $K$. A $K$-linear map $c:V\otimes
V\rightarrow V\otimes V$ is called a (quantum) \textbf{Yang-Baxter operator}
(or $YB$-operator) if it satisfies the quantum Yang-Baxter equation%
\begin{equation}
c_{1}c_{2}c_{1}=c_{2}c_{1}c_{2}  \label{ec: braided equation}
\end{equation}%
on $V\otimes V\otimes V$, where we set $c_{1}:=c\otimes V$ and $%
c_{2}:=V\otimes c.$ The pair $\left( V,c\right) $ will be called a\textbf{\
braided vector space} (or $YB$-space). A morphism of braided vector spaces $%
(V,c_{V})$ and $(W,c_{W})$ is a $K$-linear map $f:V\rightarrow W$ such that $%
c_{W}(f\otimes f)=(f\otimes f)c_{V}.$
\end{definition}

Note that, for every braided vector space $(V,c{})$ and every $\lambda \in
K, $ the pair $(V,\lambda c{})$ is a braided vector space too. A general
method for producing braided vector spaces is to take an arbitrary braided
category $(\mathcal{M},\otimes ,K,a,l,r,c),$ which is a monoidal subcategory
of the category of $K$-vector spaces. Hence any object $V\in \mathcal{M}$
can be regarded as a braided vector space with respect to $c:=c_{V,V}.$
Here, $c_{X,Y}:X\otimes Y\rightarrow Y\otimes X$ denotes the braiding in $%
\mathcal{M}{}.$ The category of comodules over a coquasitriangular Hopf
algebra and the category of Yetter-Drinfeld modules are examples of such
categories. More particularly, every bicharacter of a group $G$ induces a
braiding on the category of $G$-graded vector spaces.

\begin{definition}[Baez, \protect\cite{Ba}]
The quadruple $(A,m_{A},u_{A},c_{A})$ is called a \textbf{braided algebra} if

\begin{itemize}
\item $(A,m_{A},u_{A})$ is an associative unital algebra;

\item $(A,c_{A})$ is a braided vector space;

\item $m_{A}$ and $u_{A}$ commute with $c_{A}$, that is the following
conditions hold: 
\begin{gather}
c_{A}(m_{A}\otimes A)=(A\otimes m_{A})(c_{A}\otimes A)(A\otimes c_{A}),
\label{Br2} \\
c_{A}(A\otimes m_{A})=(m_{A}\otimes A)\left( A\otimes c_{A}\right)
(c_{A}\otimes A),  \label{Br3} \\
c_{A}(u_{A}\otimes A)=A\otimes u_{A},\qquad c_{A}(A\otimes
u_{A})=u_{A}\otimes A.  \label{Br4}
\end{gather}
\end{itemize}

A morphism of braided algebras is, by definition, a morphism of ordinary
algebras which, in addition, is a morphism of braided vector spaces.

The quadruple $(C,\Delta _{C},\varepsilon _{C},c_{C})$ is called a \textbf{%
braided coalgebra} if

\begin{itemize}
\item $(C,\Delta _{C},\varepsilon _{C})$ is a coassociative counital
coalgebra;

\item $(C,c_{C})$ is a braided vector space;

\item $\Delta _{C}$ and $\varepsilon _{C}$ commute with $c_{C}$, that is the
following relations hold: 
\begin{gather}
(\Delta _{C}\otimes C)c_{C}=(C\otimes c_{C})(c_{C}\otimes C)(C\otimes \Delta
_{C}),  \label{Br5} \\
(C\otimes \Delta _{C})c_{C}=(c_{C}\otimes C)(C\otimes c_{C})(\Delta
_{C}\otimes C),  \label{Br6} \\
(\varepsilon _{C}\otimes C)c_{C}=C\otimes \varepsilon _{C},\qquad (C\otimes
\varepsilon _{C})c_{C}=\varepsilon _{C}\otimes C.  \label{Br7}
\end{gather}
\end{itemize}

A morphism of braided coalgebras is, by definition, a morphism of ordinary
coalgebras which, in addition, is a morphism of braided vector spaces.

\cite[Definition 5.1]{Ta} A sextuple $(B,m_{B},u_{B},\Delta _{B},\varepsilon
_{B},c_{B})$ is a called a \textbf{braided bialgebra} if

\begin{itemize}
\item $(B,m_{B},u_{B},c_{B})$ is a braided algebra

\item $(B,\Delta _{B},\varepsilon _{B},c_{B}{})$ is a braided coalgebra

\item the following relations hold:%
\begin{equation}
\Delta _{B}m_{B}=(m_{B}\otimes m_{B})(B\otimes c_{B}\otimes B)(\Delta
_{B}\otimes \Delta _{B}).  \label{Br1}
\end{equation}
\end{itemize}
\end{definition}

Examples of the notions above are algebras, coalgebras and bialgebras in any
braided category $(\mathcal{M},\otimes ,K,a,l,r,c),$ which is a monoidal
subcategory of the category of $K$-vector spaces.

\begin{definition}
\label{def: graded bialg}We will need graded versions of braided algebras,
coalgebras and bialgebras. A \textbf{graded braided algebra }is a braided
algebra $(A,m_{A},u_{A},c_{A}{})$ such that $A=\bigoplus_{n\in \mathbb{N}%
}A^{n}$ and $m_{A}(A^{m}\otimes A^{n})\subseteq A^{m+n},$ for every $m,n\in 
\mathbb{N}$. The braiding $c_{A}$ is assumed to satisfy $c_{A}(A^{m}\otimes
A^{n})\subseteq A^{n}\otimes A^{m}.$ It is easy to see that $%
1_{A}=u_{A}\left( 1_{K}\right) \in A^{0}$. Therefore a graded braided
algebra is defined by maps $m_{A}^{m,n}:A^{m}\otimes A^{n}\rightarrow
A^{m+n} $ and $c_{A}^{m,n}:A^{m}\otimes A^{n}\rightarrow A^{n}\otimes A^{m}$%
, and by an element $1\in A^{0}$ such that: 
\begin{gather}
m_{A}^{n+m,p}(m_{A}^{n,m}\otimes A^{p})=m_{A}^{n,m+p}(A^{n}\otimes
m_{A}^{m,p}),\quad n,m,p\in \mathbb{N},  \label{gr1} \\
m_{A}^{0,n}(1\otimes a)=a=m_{A}^{n,0}(a\otimes 1),\quad \forall a\in
A^{n},\quad n\in \mathbb{N}.  \label{gr2} \\
c_{A}^{n+m,p}(m_{A}^{n,m}\otimes A^{p})=(A^{p}\otimes
m_{A}^{n,m})(c_{A}^{n,p}\otimes A^{m})(A^{n}\otimes c_{A}^{m,p}{}),
\label{gr3} \\
c_{A}^{n,m+p}(A^{n}\otimes m_{A}^{m,p})=(m_{A}^{m,p}\otimes
A^{n})({}A^{m}\otimes c_{A}^{n,p})(c_{A}^{n,m}\otimes A^{p}),  \label{gr4} \\
c_{A}^{0,n}(1\otimes a)=a\otimes 1\qquad \text{and\qquad }\mathfrak{\ }%
c_{A}^{n,0}(a\otimes 1)=1\otimes a,\qquad \forall a\in A^{n}.  \label{gr5}
\end{gather}%
The multiplication $m_{A}$ can be recovered from $(m_{A}^{n,m})_{n,m\in 
\mathbb{N}}$ as the unique $K$-linear map such that $m_{A}(x\otimes
y)=m_{A}^{p,q}(x\otimes y),\forall p,q\in \mathbb{N},\ \forall x\in A^{p},\
\forall y\in A^{q}.$ Analogously, the braiding $c_{A}$ is uniquely defined
by $c_{A}(x\otimes y)=c_{A}^{p,q}(x\otimes y),\forall p,q\in \mathbb{N},\
\forall x\in A^{p},\ \forall y\in A^{q}.$ We will say that $m_{A}^{n,m}$ and 
$c_{A}^{n,m}$ are the $(n,m)$-homogeneous components of $\nabla $ and $c_{A}$%
, respectively.

\textbf{Graded braided coalgebras} can by described in a similar way. By
definition a braided coalgebra $(C,\Delta _{C},\varepsilon _{C},c_{C})$ is
graded if $C=\bigoplus_{n\in \mathbb{N}}C^{n},$ $\Delta _{C}(C^{n})\subseteq
\sum_{p+q=n}C^{p}\otimes C^{q}$, $c_{C}{}(C^{n}\otimes C^{m})\subseteq
C^{m}\otimes C^{n}$ and $\varepsilon _{C}|_{C_{n}}=0$, for $n>0$ . If $\pi
^{p}$ denotes the projection onto $C^{p}$ then the comultiplication $\Delta
_{C}$ is uniquely defined by maps $\Delta _{C}^{p,q}:C^{p+q}\rightarrow
C^{p}\otimes C^{q}$, where $\Delta _{C}^{p,q}:=(\pi ^{p}\otimes \pi
^{q})\Delta _{C}|_{C^{p+q}}$. The counit is given by a map $\varepsilon
_{C}^{0}:C^{0}\rightarrow K,$ while the braiding $c_{C}{}$ is uniquely
determined by a family $(c_{C}^{n,m})_{n,m\in \mathbb{N}},$ as for braided
algebras. The families $(\Delta _{C}^{n,m})_{n,m\in \mathbb{N}},$ $%
(c_{C}^{n,m})_{n,m\in \mathbb{N}}$ and $\varepsilon _{C}^{0}$ has to satisfy
the relations that are dual to (\ref{gr1}) -- (\ref{gr5}), namely: 
\begin{gather}
(\Delta _{C}^{n,m}\otimes C^{p})\Delta _{C}^{n+m,p}=(C^{n}\otimes \Delta
_{C}^{m,p})\Delta _{C}^{n,m+p},\quad n,m,p\in \mathbb{N},  \label{c1} \\
(\varepsilon _{C}^{0}\otimes C^{n})\Delta _{C}^{0,n}(c)=c=(C^{n}\otimes
\varepsilon _{C}^{0})\Delta _{C}^{n,0}(c),\quad \forall c\in C^{n},\quad
n\in \mathbb{N}.  \label{c2} \\
{}(C^{p}\otimes \Delta _{C}^{n,m})c_{C}^{n+m,p}=(c_{C}^{n,p}\otimes
C^{m})(C^{n}\otimes c_{C}{}^{m,p})(\Delta _{C}^{n,m}\otimes C^{p}),
\label{c3} \\
(\Delta _{C}^{m,p}\otimes C^{n})c_{C}^{n,m+p}=(C^{m}\otimes
c_{C}^{n,p})(c_{C}^{n,m}\otimes C^{p})(C^{n}\otimes \Delta _{C}^{m,p}),
\label{c4} \\
(\varepsilon _{C}^{0}\otimes C)c_{C}{}(c\otimes d)=\varepsilon
_{C}^{0}(d)c=(C\otimes \varepsilon _{C}^{0})c_{C}{}(d\otimes c),\qquad
\forall c\in C^{n},\ \ \forall d\in C^{0}.  \label{c5}
\end{gather}%
We will say that $\Delta _{C}^{n,m}$ is the $(n,m)$-homogeneous component of 
$\Delta _{C}$.

A \textbf{graded braided bialgebra} is a braided bialgebra which is graded
both as an algebra and as a coalgebra.
\end{definition}

\begin{noname}
\label{nn: connected} Recall that a coalgebra $C$ is called \textbf{connected%
} if $C_{0}$, the coradical of $C$ (i.e the sum of simple subcoalgebras of $%
C $), is one dimensional. In this case there is a unique group-like element $%
1_{C}\in C$ such that $C_{0}=K1_{C}$. A morphism of connected coalgebras is
a coalgebra homomorphisms (clearly it preserves the grouplike elements).

By definition, a braided coalgebra $\left( C,c_{C}\right) $ is \textbf{%
connected} if $C_{0}=K1_{C}$ and, for any $x\in C$, 
\begin{equation}
c_{C}(x\otimes 1_{C})=1_{C}\otimes x,\qquad c_{C}(1_{C}\otimes x)=x\otimes g.
\label{de: connected}
\end{equation}
\end{noname}

\begin{definition}
Let $C=\bigoplus_{n\in \mathbb{N}}C^{n}$ be a graded braided coalgebra. We
say that $C$ is $0$\textbf{-connected}, whenever $C^{0}$ is a one
dimensional vector space.
\end{definition}

\begin{remark}
\label{re: graded connected} Let $C=\bigoplus_{n\in \mathbb{N}}C^{n}$ be a
graded braided coalgebra. By \cite[Proposition 11.1.1]{Sw}, if $%
(C_{n})_{n\in \mathbb{N}}$ is the coradical filtration, then $C_{n}\subseteq
\bigoplus_{m\leq n}C^{m}$. Therefore, $C$ is connected if it is $0$%
-connected.
\end{remark}

\begin{lemma}
\label{lem: gr(C) is braided}\cite[Lemma 1.1]{AMS-MM2} Let $(C,c)$ be a
connected braided coalgebra. Then $c$ induces a canonical braiding $c_{%
\mathrm{gr}\left( C\right) }$ on the graded coalgebra $\mathrm{gr}\left(
C\right) $ associated to the coradical filtration on $C$ such that $(\mathrm{%
gr}\left( C\right) ,c_{\mathrm{gr}\left( C\right) })$ is a $0$-connected
graded braided coalgebra.
\end{lemma}

\begin{lemma}
Let $\left( B,m_{B},u_{B},\Delta _{B},\varepsilon _{B},c_{B}\right) $ be a
graded braided bialgebra. Then%
\begin{equation}
\Delta _{B}^{a,b}m_{B}^{c,d}=\delta _{b+a,c+d}\sum\limits_{\substack{ 0\leq
u\leq c  \\ 0\leq u^{\prime }\leq d}}\left[ 
\begin{array}{c}
\delta _{a,u+u^{\prime }}\left( m_{B}^{u,u^{\prime }}\otimes
m_{B}^{c-u,d-u^{\prime }}\right) \\ 
\left( B^{u}\otimes c_{B}^{c-u,u^{\prime }}\otimes B^{d-u^{\prime }}\right)
\left( \Delta _{B}^{u,c-u}\otimes \Delta _{B}^{u^{\prime },d-u^{\prime
}}\right)%
\end{array}%
\right]  \label{form: delta mult gen}
\end{equation}
\end{lemma}

\begin{proof}
Since $\Delta _{B}^{a,b}m_{B}^{c,d}=\left( p_{B}^{a}\otimes p_{B}^{b}\right)
\Delta _{B}m_{B}\left( i_{B}^{c}\otimes i_{B}^{d}\right) ,$ using (\ref{Br1}%
) and the fact that $m_{B},c_{B},\Delta _{B}$ are graded maps one gets%
\begin{equation*}
\Delta _{B}^{a,b}m_{B}^{c,d}=\sum\limits_{\substack{ u+v=c  \\ u^{\prime
}+v^{\prime }=d}}\delta _{a,u+u^{\prime }}\delta _{b,v+v^{\prime }}\left(
m_{B}^{u,u^{\prime }}\otimes m_{B}^{v,v^{\prime }}\right) \left(
B^{u}\otimes c_{B}^{v,u^{\prime }}\otimes B^{v^{\prime }}\right) \left(
\Delta _{B}^{u,v}\otimes \Delta _{B}^{u^{\prime },v^{\prime }}\right) .
\end{equation*}%
It is easy to check that the last term of this equality coincide with the
second term of (\ref{form: delta mult gen}). 
\end{proof}

\section{Symmetric and enveloping algebras of graded braided bialgebras}

In this section we introduce and investigate the notion of symmetric algebra
of a graded braided bialgebra. We also introduce and study the notion of
universal enveloping algebra of a graded braided bialgebra endowed with a
bracket.

\begin{definitions}
\label{def: E(C)}Let $\left( C,\Delta _{C},\varepsilon _{C}\right) $ be a
graded coalgebra. Let%
\begin{equation*}
E\left( C\right) :=\bigoplus\limits_{n\in \mathbb{N}}E_{n}\left( C\right) 
\text{\quad where\quad }E_{n}\left( C\right) :=\left\{ 
\begin{tabular}{ll}
$0$ & $,$ if $n=0,1$ \\ 
$\bigcap\limits_{\substack{ a,b\geq 1  \\ a+b=n}}\ker \left( \Delta
_{C}^{a,b}\right) $ & $,$ if $n\geq 2.$%
\end{tabular}%
\right.
\end{equation*}%
Denote by $i_{C}^{n}:C^{n}\rightarrow C$ and $p_{C}^{n}:C\rightarrow C^{n}$
the canonical injection and projection respectively. Following \cite[1.2]{Ni}%
, we set%
\begin{equation*}
P_{1}\left( C\right) :=\ker \left[ \Delta _{C}-\left(
i_{C}^{0}p_{C}^{0}\otimes C\right) \Delta _{C}-\left( C\otimes
i_{C}^{0}p_{C}^{0}\right) \Delta _{C}\right] \subseteq C.
\end{equation*}
\end{definitions}

\begin{remark}
\label{rem: 0-connected}Note that $P_{1}\left( C\right) $ is just the space $%
P\left( C\right) $ of primitive elements in $C$ whenever $C$ is $0$%
-connected.
\end{remark}

The following result shows that the elements of $E\left( C\right) $ are
those of $P_{1}\left( C\right) $ of degree greater then $1$.

\begin{lemma}
\label{lem: Nichols}Let $\left( C,\Delta _{C},\varepsilon _{C}\right) $ be a
graded braided bialgebra. Then $P_{1}\left( C\right) =C^{1}\oplus \left[
\oplus _{n\geq 2}E_{n}\left( C\right) \right] .$ Moreover $C^{0}\wedge
_{C}C^{0}=C^{0}\oplus C^{1}\oplus \left[ \oplus _{n\geq 2}E_{n}\left(
C\right) \right] .$
\end{lemma}

\begin{proof}
Set $h:=\left[ \mathrm{Id}_{C\otimes C}-i_{C}^{0}p_{C}^{0}\otimes C-C\otimes
i_{C}^{0}p_{C}^{0}\right] \circ \Delta _{C}:C\rightarrow C\otimes C.$ For
every $a,b\in \mathbb{N},a+b\geq 1,$ define a map $h_{a,b}:C^{a+b}%
\rightarrow C^{a}\otimes C^{b}$ by 
\begin{equation*}
h_{a,b}:=\left\{ 
\begin{tabular}{ll}
$0,$ & if $a=0$ or $b=0,$ \\ 
$\Delta _{C}^{a,b},$ & if $a,b\geq 1.$%
\end{tabular}%
\right.
\end{equation*}%
Let $h_{n}:C^{n}\rightarrow \left( C\otimes C\right) ^{n}:=\oplus
_{a+b=n}C^{a}\otimes C^{b}$ be defined by $h_{0}:=-\Delta _{C}^{0,0}$ and $%
h_{n}:=\underset{_{a+b=n}}{\triangle }\left( h_{a,b}\right) ,$ for every $%
n\geq 1,$ where $\underset{_{a+b=n}}{\triangle }\left( h_{a,b}\right) $
denotes the diagonal morphism associated to the family $\left(
h_{a,b}\right) _{a+b=n}.$ Since $\left( C,\Delta _{C},\varepsilon
_{C}\right) $ is a graded coalgebra, we have%
\begin{eqnarray*}
h\circ i_{C}^{n} &=&\left[ \mathrm{Id}_{C\otimes
C}-i_{C}^{0}p_{C}^{0}\otimes C-C\otimes i_{C}^{0}p_{C}^{0}\right] \circ
\left( \sum\limits_{a+b=n}\left( i_{C}^{a}\otimes i_{C}^{b}\right) \circ
\Delta _{C}^{a,b}\right) \\
&=&\sum\limits_{a+b=n}\left( i_{C}^{a}\otimes i_{C}^{b}\right) \circ \Delta
_{C}^{a,b}-\left( i_{C}^{0}\otimes i_{C}^{n}\right) \circ \Delta
_{C}^{0,n}-\left( i_{C}^{n}\otimes i_{C}^{0}\right) \circ \Delta _{C}^{n,0}
\end{eqnarray*}%
Now, if $n=0,$ we obtain%
\begin{eqnarray*}
h\circ i_{C}^{n} &=&\sum\limits_{a+b=0}\left( i_{C}^{a}\otimes
i_{C}^{b}\right) \circ \Delta _{C}^{a,b}-\left( i_{C}^{0}\otimes
i_{C}^{0}\right) \circ \Delta _{C}^{0,0}-\left( i_{C}^{0}\otimes
i_{C}^{0}\right) \circ \Delta _{C}^{0,0} \\
&=&\left( i_{C}^{0}\otimes i_{C}^{0}\right) \circ \left( -\Delta
_{C}^{0,0}\right) =i_{C\otimes C}^{0}\circ h_{0}.
\end{eqnarray*}%
If $n\geq 1,$ we obtain%
\begin{eqnarray*}
h\circ i_{C}^{n} &=&\sum\limits_{\substack{ a,b\geq 1  \\ a+b=n}}\left(
i_{C}^{a}\otimes i_{C}^{b}\right) \circ \Delta
_{C}^{a,b}=\sum\limits_{a+b=n}\left( i_{C}^{a}\otimes i_{C}^{b}\right) \circ
h_{a,b} \\
&=&\sum\limits_{a+b=n}\left( i_{C}^{a}\otimes i_{C}^{b}\right) \circ
p_{\left( C\otimes C\right) ^{a+b}}^{a,b}\circ h_{a+b}=i_{C\otimes
C}^{n}\circ h_{n},
\end{eqnarray*}%
where $p_{\left( C\otimes C\right) ^{a+b}}^{a,b}:\left( C\otimes C\right)
^{a+b}\rightarrow C^{a}\otimes C^{b}$ denotes the canonical projection. This
entails that $h:C\rightarrow C\otimes C$ is a graded homomorphism so that%
\begin{eqnarray*}
P_{1}\left( C\right) &=&\ker \left( h\right) =\oplus _{n\in \mathbb{N}}\ker
\left( h_{n}\right) =\ker \left( h_{0}\right) \oplus \ker \left(
h_{1}\right) \oplus \left[ \oplus _{n\geq 2}\ker \left( h_{n}\right) \right]
\\
&=&C^{1}\oplus \left[ \oplus _{n\geq 2}E_{n}\left( C\right) \right] .
\end{eqnarray*}%
Now, by \cite[after Proposition 1.2.1]{Ni}, we have that $C^{0}\wedge
_{C}C^{0}=C^{0}\oplus P_{1}\left( C\right) $ so that%
\begin{equation*}
C^{0}\wedge _{C}C^{0}=C^{0}\oplus P_{1}\left( C\right) =C^{0}\oplus
C^{1}\oplus \left[ \oplus _{n\geq 2}E_{n}\left( C\right) \right] .
\end{equation*}
\end{proof}

\begin{definition}
\label{def: symmetric}Let $\left( B,m_{B},u_{B},\Delta _{B},\varepsilon
_{B},c_{B}\right) $ be a graded braided bialgebra. The \textbf{symmetric
algebra of }$B$ is defined to be%
\begin{equation}
S\left( B\right) =\frac{B}{\left( E\left( B\right) \right) }.
\label{form: S(B)}
\end{equation}%
Denote by $\pi _{S}:=\pi _{S}^{B}:B\rightarrow S\left( B\right) $ the
canonical projection.
\end{definition}

\begin{theorem}
\label{teo: S(B) graded}Let $\left( B,m_{B},u_{B},\Delta _{B},\varepsilon
_{B},c_{B}\right) $ be a graded braided bialgebra. Then $S\left( B\right) $
has a unique graded braided bialgebra structure such that the canonical
projection $\pi _{S}:B\rightarrow S\left( B\right) $ is a graded braided
bialgebra homomorphism.
\end{theorem}

\begin{proof}
Since $I:=\left( E\left( B\right) \right) $ is generated by homogeneous
elements, it is a graded ideal of $B$ 
with graded component $I_{n}=\left( E\left( B\right) \right) \cap B^{n}.$
Hence $S:=S\left( B\right) $ has a unique graded algebra structure such that 
$\pi _{S}:B\rightarrow S\left( B\right) $ is a graded algebra homomorphism.
Let us check it is also a graded coalgebra.

Since $B$ is a graded coalgebra and $E_{t}\left( B\right) =0$ if $n=0,1,$ we
have $\varepsilon _{B}\left( I\right) =0$ so that, in order to get that $I$
is a coideal, it remains to prove that $\Delta \left( I\right) \subseteq
I\otimes B+B\otimes I.$

Firstly, let us check that 
\begin{equation}
c_{B}\left( B^{u}\otimes E_{t}\left( B\right) \right) \subseteq E_{t}\left(
B\right) \otimes B^{u}\text{\qquad and\qquad }c_{B}\left( E_{t}\left(
B\right) \otimes B^{u}\right) \subseteq B^{u}\otimes E_{t}\left( B\right) .
\label{form: braided 1}
\end{equation}%
We focus on the second inclusion the first one having a similar proof. Since 
$B$ is a $c$-bracket, for every $z\in E_{t}\left( B\right) ,y\in B^{u}$ and $%
a,b\geq 1,a+b=t,$ we have 
\begin{equation*}
\left( B^{u}\otimes \Delta _{B}^{a,b}\right) c_{B}^{t,u}\left( z\otimes
y\right) \overset{\text{(\ref{c3})}}{=}\left( c_{B}^{a,u}\otimes
B^{b}\right) \left( B^{a}\otimes c_{B}^{b,u}\right) \left( \Delta
_{B}^{a,b}\otimes B^{u}\right) \left( z\otimes y\right) =0
\end{equation*}%
Then $c_{B}^{t,u}\left( z\otimes y\right) \in \ker \left( B^{u}\otimes
\Delta _{B}^{a,b}\right) =B^{u}\otimes \ker \left( \Delta _{B}^{a,b}\right)
, $ for every $a,b\geq 1,a+b=t$ so that 
\begin{equation*}
c_{B}^{t,u}\left( z\otimes y\right) \in \bigcap\limits_{\substack{ a,b\geq
1,  \\ a+b=t}}\left[ B^{u}\otimes \ker \left( \Delta _{B}^{a,b}\right) %
\right] =B^{u}\otimes \left( \bigcap\limits_{\substack{ a,b\geq 1,  \\ a+b=t 
}}\ker \left( \Delta _{B}^{a,b}\right) \right) =B^{u}\otimes E_{t}\left(
B\right) .
\end{equation*}%
Thus (\ref{form: braided 1}) holds.

Now, for every $z\in E_{t}\left( B\right) ,x\in B^{u},y\in B^{v},s\in B^{w},$
using (\ref{gr3}) and (\ref{form: braided 1}) it is straightforward to check
that $c_{B}^{u+t+v,w}\left( xzy\otimes s\right) \in \left( B^{w}\otimes
I_{u+t+v}\right) $%
so that 
\begin{equation}
c_{B}\left( B^{u}\otimes I_{t}\right) \subseteq I_{t}\otimes B^{u}\text{
\qquad and\qquad\ }c_{B}\left( I_{t}\otimes B^{u}\right) \subseteq
B^{u}\otimes I_{t}.  \label{form: braided 2}
\end{equation}%
Let $z\in E_{t}\left( B\right) .$ Thus, for every $a,b\geq 1$ such that $%
a+b=t,$ and we have $\left( B^{0}\otimes \Delta _{B}^{a,b}\right) \Delta
_{B}^{0,t}\left( z\right) =\left( \Delta _{B}^{0,a}\otimes B^{b}\right)
\Delta _{B}^{a,b}\left( z\right) =0$ so that%
\begin{equation*}
\Delta _{B}^{0,t}\left( z\right) \in \bigcap\limits_{\substack{ a,b\geq 1, 
\\ a+b=t}}\left[ B^{0}\otimes \ker \left( \Delta _{B}^{a,b}\right) \right]
=B^{0}\otimes \left[ \bigcap\limits_{\substack{ a,b\geq 1,  \\ a+b=t}}\left(
\ker \Delta _{B}^{0,a}\right) \right] =B^{0}\otimes E_{t}\left( B\right) .
\end{equation*}%
Therefore we have%
\begin{equation}
\Delta _{B}^{0,t}\left( E_{t}\left( B\right) \right) \subseteq B^{0}\otimes
E_{t}\left( B\right) \qquad \text{and}\qquad \Delta _{B}^{t,0}\left(
E_{t}\left( B\right) \right) \subseteq E_{t}\left( B\right) \otimes B^{0}
\label{form: braided 3}
\end{equation}%
where the second inclusion can be checked similarly.

Let $x\in B^{c}$ and $z\in E_{t}\left( B\right) .$ Using (\ref{form: delta
mult gen}), that $z\in E_{t}\left( B\right) $ and (\ref{form: braided 3}),
we get $\Delta _{B}^{a,b}\left( xz\right) \in B^{a}\otimes
I_{b}+I_{a}\otimes B^{b}$ %
so that $\Delta _{B}^{a,b}\left( I_{a+b}\right) \subseteq B^{a}\otimes
I_{b}+I_{a}\otimes B^{b},$ for every $a,b\in \mathbb{N} .$

Therefore $I$ is a graded coideal and hence $S:=S\left( B\right) $ carries a
unique graded coalgebra structure such that $\pi _{S}$ is a coalgebra
homomorphism. In fact there is a unique morphism $\Delta
_{S}^{a,b}:B^{a+b}/I_{a+b}\rightarrow B^{a}/I_{a}\otimes B^{b}/I_{b}$ such
that $\Delta _{S}^{a,b}\pi _{S}^{a+b}=\left( \pi _{S}^{a}\otimes \pi
_{S}^{b}\right) \Delta _{B}^{a,b},$ where $\pi _{S}^{n}:B^{n}\rightarrow
S^{n}:=B^{n}/I_{n}$ denotes the canonical projection. Using that $\pi
_{S}^{n}$ is an epimorphism for every $n\in \mathbb{N} $, one gets that $%
\left( \Delta _{S}^{a,b}\otimes S^{c}\right) \Delta _{S}^{a+b,c}=\left(
S^{a}\otimes \Delta _{S}^{b,c}\right) \Delta _{S}^{a,b+c}.$ On the other
hand $\pi _{S}^{0}$ is bijective so that we can set $\varepsilon
_{S}^{0}:=\varepsilon _{B}^{0}\left( \pi _{S}^{0}\right) ^{-1}$. One easily
checks that $\left( \varepsilon _{S}^{0}\otimes S^{c}\right) \Delta
_{S}^{0,c}=\mathrm{Id}_{S^{c}}=\left( S^{c}\otimes \varepsilon
_{S}^{0}\right) \Delta _{S}^{c,0}.$ By \cite[Proposition 2.5]{AM- Type One},
there is a unique map $\Delta _{S}:S\rightarrow S\otimes S$ such that $%
\Delta i_{S}^{n}=\sum\limits_{a+b=n}\left( i_{S}^{a}\otimes i_{S}^{b}\right)
\Delta _{S}^{a,b}$ where $i_{S}^{n}:S^{n}\rightarrow S$ denotes the
canonical injection. Furthermore $\left( S,\Delta _{S},\varepsilon
_{S}=\varepsilon _{S}^{0}p_{S}^{0}\right) $ is a graded coalgebra where $%
p_{S}^{0}:S\rightarrow S^{0}$ denotes the canonical projection.

Let us prove that $c_{B}$ factors trough a braiding $c_{S}$ of $S=S\left(
B\right) $ that makes $S$ a braided bialgebra.

By the foregoing we get%
\begin{equation*}
c_{B}^{a,b}\left( \ker \left( \pi _{S}^{a}\otimes \pi _{S}^{b}\right)
\right) =c_{B}^{a,b}\left( B^{a}\otimes I_{b}+I_{a}\otimes B^{b}\right) 
\overset{\text{(\ref{form: braided 2})}}{\subseteq }I_{b}\otimes
B^{a}+B^{b}\otimes I_{a}=\ker \left( \pi _{S}^{b}\otimes \pi _{S}^{a}\right)
\end{equation*}%
so that, since $\pi _{S}\otimes \pi _{S}$ is surjective, there exists a
unique $K$-linear map $c_{S}^{a,b}:S^{a}\otimes S^{b}\rightarrow
S^{b}\otimes S^{a}$ such that%
\begin{equation}
c_{S}^{a,b}\left( \pi _{S}^{b}\otimes \pi _{S}^{a}\right) =\left( \pi
_{S}^{a}\otimes \pi _{S}^{b}\right) c_{B}^{a,b}.  \label{form: braiding S(B)}
\end{equation}%
This relation can be used to prove that $\left( S,c_{S}\right) $ is a graded
braided bialgebra with structure induced by $\pi _{S}.$
\end{proof}

\begin{proposition}
\label{pro: S(f)}Let $f:\left( B,m_{B},u_{B},\Delta _{B},\varepsilon
_{B},c_{B}\right) \rightarrow \left( B^{\prime },m_{B^{\prime
}},u_{B^{\prime }},\Delta _{B^{\prime }},\varepsilon _{B^{\prime
}},c_{B^{\prime }}\right) $ be a graded braided bialgebra homomorphism. Then
there is a unique algebra homomorphism%
\begin{equation*}
S\left( f\right) :S\left( B\right) \rightarrow S\left( B^{\prime }\right)
\end{equation*}%
such that $S\left( f\right) \circ \pi _{S}=\pi _{S^{\prime }}\circ f.$
Moreover $S\left( f\right) $ is a graded braided bialgebra homomorphism.
\end{proposition}

\begin{proof}
Let $n\geq 2.$ Then, for every $a,b\geq 1,a+b=n$, we have 
\begin{equation*}
\Delta _{B^{\prime }}^{a,b}f_{n}\left( E_{n}\left( B\right) \right) =\left(
f_{a}\otimes f_{b}\right) \Delta _{B}^{a,b}\left( E_{n}\left( B\right)
\right) =0
\end{equation*}%
so that $f_{n}\left( E_{n}\left( B\right) \right) \subseteq \bigcap\limits 
_{\substack{ a,b\geq 1  \\ a+b=n}}\ker \left( \Delta _{B}^{a,b}\right)
=E_{n}\left( B^{\prime }\right) .$ Set $I:=E_{n}\left( B\right) $ and $%
I^{\prime }:=E_{n}\left( B^{\prime }\right) .$ Denote by $I_{n}$ and $%
I_{n}^{\prime }$ the graded components of these ideals. Since $f$ is a
graded algebra homomorphism, we get $f\left( I_{n}\right) \subseteq
I_{n}^{\prime }$ so that $f$ is a morphism of graded ideals. Hence there is
a unique algebra homomorphism $S\left( f\right) :S\left( B\right)
\rightarrow S\left( B^{\prime }\right) $ such that $S\left( f\right) \circ
\pi _{S}=\pi _{S^{\prime }}\circ f.$ Furthermore $S\left( f\right) $ is a
morphism of graded algebras with graded component $S\left( f^{n}\right)
:S^{n}\rightarrow \left( S^{\prime }\right) ^{n}$ uniquely defined by $%
S\left( f^{n}\right) \circ \pi _{S}^{n}=\pi _{S^{\prime }}^{n}\circ f^{n}.$
Since $\pi _{S}$ is an epimorphism, one gets that $S\left( f\right) $ is a
coalgebra homomorphism (and hence a graded coalgebra homomorphism.
\end{proof}

Next aim is to introduce a notion universal enveloping algebra of a graded
braided bialgebra endowed with a bracket and to see that it carries a
braided bialgebra structure.

\begin{claim}
Let $\left( B,m_{B},u_{B},\Delta _{B},\varepsilon _{B},c_{B}\right) $ be a
graded braided bialgebra. By (\ref{form: braided 1}), we have the the inclusions $c_{B}\left(
B^{u}\otimes E_{t}\left( B\right) \right) \subseteq E_{t}\left( B\right)
\otimes B^{u}$ and $c_{B}\left( E_{t}\left( B\right) \otimes B^{u}\right)
\subseteq B^{u}\otimes E_{t}\left( B\right) $ hold. Hence there exists a unique
morphism $c_{E_{t}\left( B\right) ,B^{u}}:E_{t}\left( B\right) \otimes
B^{u}\rightarrow B^{u}\otimes E_{t}\left( B\right) $ such that%
\begin{equation}
\left( B^{u}\otimes j_{B}^{t}\right) \circ c_{E_{t}\left( B\right)
,B^{u}}=c_{B}^{t,u}\circ \left( j_{B}^{t}\otimes B^{u}\right) ,
\label{form: c_E_t(B)}
\end{equation}%
where $j_{B}^{t}:E_{t}\left( B\right) \rightarrow B^{t}$ denotes the
canonical injection. Similarly one gets $c_{B^{u},E_{t}\left( B\right) }.$
Define now%
\begin{eqnarray*}
c_{E\left( B\right) ,B^{u}} &:&=\oplus _{t\in \mathbb{N}}c_{E_{t}\left(
B\right) ,B^{u}}:E\left( B\right) \otimes B^{u}\rightarrow B^{u}\otimes
E\left( B\right) , \\
c_{B^{u},E\left( B\right) } &:&=\oplus _{t\in \mathbb{N}}c_{B^{u},E_{t}%
\left( B\right) }:B^{u}\otimes E\left( B\right) \rightarrow E\left( B\right)
\otimes B^{u}.
\end{eqnarray*}
\end{claim}

\begin{definition}
\label{def: bracket}A \emph{bracket} for a graded braided bialgebra $\left(
B,m_{B},u_{B},\Delta _{B},\varepsilon _{B},c_{B}\right) $ is a $K$-linear
map $b=b_{B}:E\left( B\right) \rightarrow B^{1}$ such that, for every $t\in 
\mathbb{N} $ 
\begin{equation}
c_{B}^{1,t}\left( b\otimes B^{t}\right) =\left( B^{t}\otimes b\right)
c_{E\left( B\right) ,B^{t}}\qquad c_{B}^{t,1}\left( B^{t}\otimes b\right)
=\left( b\otimes B^{t}\right) c_{B^{t},E\left( B\right) }
\label{form: Bracket}
\end{equation}%
and%
\begin{equation}
\Delta _{B}^{1,0}b=\left( b\otimes B^{0}\right) \Delta _{B}^{t,0}\qquad
\Delta _{B}^{0,1}b=\left( B^{0}\otimes b\right) \Delta _{B}^{0,t}.
\label{form: Bracket2}
\end{equation}%
The restriction of $b$ to $E_{t}\left( B\right) $ will be denoted by $%
b^{t}:E_{t}\left( B\right) \rightarrow B^{1}.$ If $b$ is a bracket for $B,$
then we define \textbf{the universal enveloping algebra of }$\left(
B,b\right) $ to be 
\begin{equation*}
U\left( B,b\right) :=\frac{B}{\left( \left( \mathrm{Id}-b\right) \left[
E\left( B\right) \right] \right) }.
\end{equation*}%
We will denote by $i_{U}:B^{1}\rightarrow U\left( B,b\right) $ the
restriction to $B^{1}$ of the canonical projection $\pi _{U}:B\rightarrow
U\left( B,b\right) $. Note that $0$ is always a bracket for $B$ so that it
makes sense to consider%
\begin{equation*}
U\left( B,0\right) =\frac{B}{\left( E\left( B\right) \right) }\overset{\text{%
(\ref{form: S(B)})}}{=}S\left( B\right) .
\end{equation*}%
Let $f:B\rightarrow B^{\prime }$ be a graded braided bialgebra homomorphism.
Then one has $\Delta _{B^{\prime }}^{a,b}\circ f_{a+b}=\left( f_{a}\otimes
f_{b}\right) \circ \Delta _{B}^{a,b},$ for every $a,b\in \mathbb{N} $. This
entails that $f_{t}\left[ E_{t}\left( B\right) \right] \subseteq E_{t}\left(
B^{\prime }\right) $ so that there exists a unique map%
\begin{equation*}
E_{t}\left( f\right) :E_{t}\left( B\right) \rightarrow E_{t}\left( B^{\prime
}\right)
\end{equation*}%
such that $j_{t}^{B^{\prime }}\circ E_{t}\left( f\right) =f_{t}\circ
j_{t}^{B},$ for every $t\geq 2$, where $j_{t}^{B}:E_{t}\left( B\right)
\rightarrow B^{t}$ and $j_{t}^{B^{\prime }}:E_{t}\left( B^{\prime }\right)
\rightarrow \left( B^{\prime }\right) ^{t}$ are the canonical inclusions. We
set%
\begin{equation*}
E\left( f\right) :=\oplus _{n\in \mathbb{N} }E_{t}\left( f\right) .
\end{equation*}%
A \emph{morphism of brackets} $f:\left( B,b_{B}\right) \rightarrow \left(
B^{\prime },b_{B^{\prime }}\right) $ is a graded braided bialgebra
homomorphism $f:B\rightarrow B^{\prime }$ such that $f_{1}\circ
b_{B}^{t}=b_{B^{\prime }}^{t}\circ E_{t}\left( f\right) ,$ for every $t\in 
\mathbb{N} $.
\end{definition}

\begin{theorem}
\label{teo: U bialgebra}Let $\left( B,m_{B},u_{B},\Delta
_{B},\varepsilon _{B},c_{B}\right) $ be a graded braided bialgebra endowed with a
bracket $b:E\left( B\right) \rightarrow B^{1}$. Then $U\left( B,b\right) $ has a
unique braided bialgebra structure such that the canonical projection $\pi
_{U}:B\rightarrow U\left( B,b\right) $ is a braided bialgebra homomorphism.
\end{theorem}

\begin{proof}
Set $I=\left( \left( \mathrm{Id}-b\right) \left( E\left( B\right) \right)
\right) .$ Let us check that $I$ is a coideal i.e. that $\Delta \left(
I\right) \subseteq I\otimes B+B\otimes I$ and $\varepsilon _{B}\left(
I\right) =0.$ First, for every $z\in E_{t}\left( B\right) ,x\in B^{u},y\in
B^{v},s\in B^{w},$ using (\ref{gr3}), (\ref{form: Bracket}) and (\ref{form:
braided 1}) one can prove that $c_{B}\left( x\left( z-b\left( z\right)
\right) y\otimes s\right) \in B^{w}\otimes I$%
so that%
\begin{equation}
c_{B}\left( I\otimes B^{w}\right) \subseteq B^{w}\otimes I\qquad \text{and}%
\qquad c_{B}\left( B^{w}\otimes I\right) \subseteq I\otimes B^{w}
\label{form: braided 5}
\end{equation}%
where the second equality can be obtained similarly. Let $z\in E_{t}\left(
B\right) .$ By Lemma \ref{lem: Nichols}, one has 
\begin{equation*}
u=z-b\left( z\right) \in E_{t}\left( B\right) -B^{1}\subseteq B^{1}\oplus 
\left[ \oplus _{n\geq 2}E_{n}\left( B\right) \right] =P_{1}\left( C\right)
\end{equation*}%
so that $\Delta _{B}\left( u\right) =\left( i_{B}^{0}p_{B}^{0}\otimes
B\right) \Delta _{B}\left( u\right) +\left( B\otimes
i_{B}^{0}p_{B}^{0}\right) \Delta _{B}\left( u\right) .$ We have $\varepsilon
_{B}\left( I\right) \subseteq \varepsilon _{B}\left( P_{1}\left( C\right)
\right) =0.$ Moreover, using (\ref{form: Bracket2}) and (\ref{form: braided
3}) one easily gets that $\left( p_{B}^{0}\otimes B\right) \Delta _{B}\left(
u\right) \in B^{0}\otimes I$
and similarly $\left( B\otimes p_{B}^{0}\right) \Delta _{B}\left( u\right)
\in I\otimes B^{0}$ so that $\Delta _{B}\left( u\right) \in B\otimes
I+I\otimes B.\text{ } 
$ Let $x\in B^{c}.$ Using (\ref{Br1}), the last relation and (\ref{form:
braided 5}) we get that $\Delta _{B}\left( xu\right) \in B\otimes I+I\otimes
B. 
$ %
%
Using (\ref{Br1}), the last relation and (\ref{form: braided 5}), for every $%
y\in B^{d}\ $we obtain $\Delta _{B}\left( xuy\right) \in B\otimes I+I\otimes
B.$

In conclusion $\Delta \left( I\right) \subseteq I\otimes B+B\otimes I.$
Therefore $I$ is coideal and hence $U=U\left( B,b\right) $ carries a unique
coalgebra structure such that the canonical projection $\pi
_{U}:T\rightarrow U$ is a coalgebra homomorphism. In fact there are unique
morphisms $\Delta _{U}:B/I\rightarrow B/I\otimes B/I$ and $\varepsilon
_{U}:B/I\rightarrow K$ such that $\Delta _{U}\pi _{U}=\left( \pi _{U}\otimes
\pi _{U}\right) \Delta _{B}$ and $\varepsilon _{U}\pi _{U}=\varepsilon _{B}$%
. Using that $\pi _{U}$ is an epimorphism and that $B$ is a coalgebra one
gets that $\left( U,\Delta _{U},\varepsilon _{U}\right) $ is a coalgebra
too. Let us prove that $c_{B}$ factors trough a braiding $c_{U}$ of $U$ that
makes $U$ a braided bialgebra. We get%
\begin{equation*}
c_{B}\left( \ker \left( \pi _{U}\otimes \pi _{U}\right) \right) =c_{B}\left(
I\otimes T+T\otimes I\right) \overset{\text{(\ref{form: braided 5})}}{%
\subseteq }I\otimes T+T\otimes I=\ker \left( \pi _{U}\otimes \pi _{U}\right)
\end{equation*}%
so that, since $\pi _{U}\otimes \pi _{U}$ is surjective, there exists a
unique $K$-linear map $c_{U}:U\otimes U\rightarrow U\otimes U$ such that $%
c_{U}\left( \pi _{U}\otimes \pi _{U}\right) =\left( \pi _{U}\otimes \pi
_{U}\right) c_{T}.$ This relation can be used to prove that $\left(
U,m_{U},u_{U},\Delta _{U},\varepsilon _{U},c_{U}\right) $ is a braided
bialgebra with structure induced by $\pi _{U}.$ Here $m_{U}:U\otimes
U\rightarrow U$ and $u_{U}:K\rightarrow U$ denote the multiplication and the
unit of $U$ respectively.
\end{proof}

\section{Strongness degree\label{sec: strongness degree}}

In this section we introduce and study the notion of strongness degree of a
graded braided bialgebra $B$ which is somehow a measure of how far $B$ is to
be strongly $\mathbb{N}$-graded as a coalgebra. Later on we will focus on
braided bialgebras of strongness degree at most one. Further results on
strongness degree of graded braided bialgebras are included in appendix.

\begin{definition}
\label{def: strongly grAlg}\cite[Definition 3.5]{AM- Type One} Let $%
(A=\oplus _{n\in \mathbb{N} }A^{n},m,u)$ be a graded algebra in $\mathcal{M}$%
. In analogy with the group graded case, we say that $A$ is a \emph{strongly 
}$\mathbb{N} $\emph{-graded algebra} whenever $m_{A}^{i,j}:A^{i}\otimes
A^{j}\rightarrow A^{i+j}$ is an epimorphism for every $i,j\in \mathbb{N}$.

\label{def: strongly grCoalg}\cite[Definition 2.9]{AM- Type One} (see also 
\cite[Lemma 2.3]{MS-Pointed indec}) Let $(C=\oplus _{n\in \mathbb{N}
}C^{n},\Delta ,\varepsilon )$ be a graded coalgebra in $\mathcal{M}$. We say
that $C$ is a \emph{strongly }$\mathbb{N} $\emph{-graded coalgebra} whenever 
$\Delta _{C}^{i,j}:C^{i+j}\rightarrow C^{i}\otimes C^{j}$ is a monomorphism
for every $i,j\in \mathbb{N}$.
\end{definition}

\begin{claim}
Recall that, given a coalgebra $C$ and a $C$-bicomodule $M,$ one can
consider the cotensor coalgebra $\left( T^{c}:=T_{C}^{c}(M),\Delta
_{T^{c}},\varepsilon _{T^{c}}\right) $. It is defined as follows. As a
vector space $T^{c}=C\oplus M\oplus M^{\square _{C}2}\oplus M^{\square
_{C}3}\oplus \cdots ,$ where $\square _{C}$ denotes the cotensor product
over $C$. Furthermore, $\Delta _{T^{c}}^{i,j}$ is given by the comodule
structure map for $i=0$ or $j=0;$ when $i,j>0$ it is induced by the map $%
m_{1}\otimes \cdots \otimes m_{n}\mapsto \left( m_{1}\otimes \cdots \otimes
m_{i}\right) \otimes \left( m_{i+1}\otimes \cdots \otimes m_{n}\right) .$
The counit is zero on elements of positive degree. The cotensor coalgebra
fulfill a suitable universal property that resemble the universal property
of the tensor algebra (see \cite[Proposition 1.4.2]{Ni}). Note that $T^{c}$
is a strongly $%
\mathbb{N}
$-graded coalgebra.
\end{claim}

\begin{theorem}
\label{teo: strongness}Let $\left( B,m_{B},u_{B},\Delta _{B},\varepsilon
_{B},c_{B}\right) $ be a graded braided bialgebra. The following assertions
are equivalent.

\begin{enumerate}
\item The canonical projection $\pi _{S}:B\rightarrow S\left( B\right) $ is
an isomorphism.

\item $E_{n}\left( B\right) =0$ for every $n\geq 2$.

\item $B$ is strongly $%
\mathbb{N}
$-graded as a coalgebra.

\item $\Delta _{B}^{n,1}$ is injective for every $n\in 
\mathbb{N}
$.

\item The canonical map $\psi :B\rightarrow T_{B^{0}}^{c}\left( B^{1}\right) 
$ is injective.

\item $B^{0}\oplus B^{1}\oplus \cdots \oplus B^{n}=\underset{n+1}{%
\underbrace{B^{0}\wedge _{B}\cdots \wedge _{B}B^{0}}},$ for every $n\in 
\mathbb{N}
.$

\item $B^{0}\oplus B^{1}=B^{0}\wedge _{B}B^{0}.$
\end{enumerate}
\end{theorem}

\begin{proof}
By Theorem \ref{teo: S(B) graded}, $S\left( B\right) $ is a graded
coalgebra. Thus the equivalences $\left( 3\right) \Leftrightarrow \left(
4\right) \Leftrightarrow \left( 5\right) \Leftrightarrow \left( 6\right)
\Leftrightarrow \left( 7\right) $ follows by applying \cite[Theorem 2.22]%
{AM- Type One} to the monoidal category of vector spaces. $\left( 1\right)
\Leftrightarrow \left( 2\right) $ It is trivial. $\left( 2\right)
\Leftrightarrow \left( 7\right) $ By Lemma \ref{lem: Nichols}, we have $%
B^{0}\wedge _{B}B^{0}=B^{0}\oplus B^{1}\oplus \left[ \oplus _{n\geq
2}E_{n}\left( B\right) \right] .$
\end{proof}

\begin{definition}
\label{def: strongness degree}Let $\left( B,m_{B},u_{B},\Delta
_{B},\varepsilon _{B},c_{B}\right) $ be a graded braided bialgebra. Define
iteratively $S^{\left[ n\right] }\left( B\right) $ by $S^{\left[ 0\right]
}\left( B\right) :=B$ and $S^{\left[ n\right] }\left( B\right) :=S\left( S^{%
\left[ n-1\right] }\left( B\right) \right) ,$ for $n\geq 1.$ This yields a
direct system $\left( \left( S^{\left[ i\right] }\left( B\right) \right)
_{i\in 
\mathbb{N}
},\left( \pi _{i}^{j}\right) _{i,j\in 
\mathbb{N}
}\right) $%
\begin{equation*}
B=S^{\left[ 0\right] }\left( B\right) \overset{\pi _{0}^{1}}{\rightarrow }S^{%
\left[ 1\right] }\left( B\right) \overset{\pi _{1}^{2}}{\rightarrow }S^{%
\left[ 3\right] }\left( B\right) \rightarrow \cdots
\end{equation*}%
where $\pi _{i}^{j}$ is defined in an obvious way for every $j\geq i$.
Denote by $\left( S^{\left[ \infty \right] }\left( B\right) ,\pi
_{i}^{\infty }\right) =\underrightarrow{\lim }\left( S^{\left[ i\right]
}\left( B\right) \right) $ the direct limit of this direct system. We say
that $B$\textbf{\ has strongness degree }$n$ if the system above is
stationary and $n$ is the least $t\in 
\mathbb{N}
$ such that $\pi _{t}^{t+1}:S^{\left[ t\right] }\left( B\right) \rightarrow
S^{\left[ t+1\right] }\left( B\right) $ is an isomorphism. In this case we
will write $\mathrm{sdeg}\left( B\right) =n.$ If the system above is not
stationary we will write $\mathrm{sdeg}\left( B\right) =\infty .$
\end{definition}

\begin{remark}
\label{rem: combinatorial}V. K. Kharchenko pointed out to our attention that
the notion of \textbf{combinatorial rank} in \cite[Definition 5.4]%
{Kharchenko-SkewPrim} is essentially the same notion as the strongness
degree. Nevertheless we decided here to keep our terminology as, in view of
Theorem \ref{teo: strongness}, strongness degree is a measure of how far $B$
is to be strongly $%
\mathbb{N}
$-graded as a coalgebra.
\end{remark}

\begin{definition}
\label{def: Gamma}Let $\left( B,m_{B},u_{B},\Delta _{B},\varepsilon
_{B},c_{B}\right) $ be a graded braided bialgebra. For every $a,b,n\in 
\mathbb{N}
$, we set 
\begin{equation}
\mathrm{\Gamma }_{a,b}^{B}:=m_{B}^{a,b}\Delta _{B}^{a,b}\qquad \text{and}%
\qquad \mathrm{\Gamma }_{n}^{B}:=\left\{ 
\begin{tabular}{ll}
$\mathrm{Id}_{B^{0}}$ & if $n=0$ \\ 
$m_{B}^{n-1,1}\left( \mathrm{\Gamma }_{n-1}^{B}\otimes B^{1}\right) \Delta
_{B}^{n-1,1}$ & if $n\geq 1$%
\end{tabular}%
\right. \text{.}  \label{form: SS}
\end{equation}
\end{definition}

\begin{remark}
\label{rem: Gamma}Let $\left( V,c\right) $ be a braided vector space and let 
$T=T\left( V,c\right) $ be the associated tensor algebra. Then $m_{T}^{a,b}$
is just the juxtaposition. This entails that $\mathrm{\Gamma }_{a,b}^{T}=%
\mathrm{S}_{a,b}$ and $\mathrm{\Gamma }_{n}^{T}=\mathrm{S}_{n}$ where $%
\mathrm{S}_{a,b}$ and $\mathrm{S}_{n}$ are the morphisms defined in \cite[%
page 2815]{Schauenburg0}.
\end{remark}

\begin{lemma}
Let $\pi :\left( B,m_{B},u_{B},\Delta _{B},\varepsilon _{B},c_{B}\right)
\rightarrow \left( E,m_{E},u_{E},\Delta _{E},\varepsilon _{E},c_{E}\right) $
be a graded braided bialgebra homomorphism. Then 
\begin{equation}
\mathrm{\Gamma }_{a,b}^{E}\circ \pi ^{a+b}=\pi ^{a+b}\circ \mathrm{\Gamma }%
_{a,b}^{B},\text{ for every }a,b\in 
\mathbb{N}
,a+b\geq 1.  \label{form: S bar 2}
\end{equation}%
where, for every $n\in 
\mathbb{N}
$, $\pi ^{n}:B^{n}\rightarrow E^{n}$ denotes the $n$-th graded component of $%
\pi $.
\end{lemma}

\begin{proof}
We have $\mathrm{\Gamma }_{a,b}^{E}\pi ^{a+b}=m_{E}^{a,b}\Delta
_{E}^{a,b}\pi ^{a+b}=m_{E}^{a,b}\left( \pi ^{a}\otimes \pi ^{b}\right)
\Delta _{B}^{a,b}=\pi ^{a+b}m_{B}^{a,b}\Delta _{B}^{a,b}=\pi ^{a+b}\mathrm{%
\Gamma }_{a,b}^{B}.$
\end{proof}

\begin{lemma}
Let $\left( B,m_{B},u_{B},\Delta _{B},\varepsilon _{B},c_{B}\right) $ be a
graded braided bialgebra which is strongly $%
\mathbb{N}
$-graded as an algebra. Then%
\begin{equation}
\mathrm{\Gamma }_{a+b}^{B}=m_{B}^{a,b}\left( \mathrm{\Gamma }_{a}^{B}\otimes 
\mathrm{\Gamma }_{b}^{B}\right) \Delta _{B}^{a,b},\text{ for every }a,b\in 
\mathbb{N}
.  \label{form: SS bar}
\end{equation}
\end{lemma}

\begin{proof}
By \cite[Theorem 3.11]{AM- Type One}, there is a unique algebra homomorphism 
$\varphi _{B}:T:=T_{B^{0}}\left( B^{1}\right) \rightarrow B$ that restricted
to $B^{0}$ and $B^{1}$ gives the canonical inclusions. Furthermore each
graded component of $\varphi _{B}$ is an epimorphism as $B$ is strongly $%
\mathbb{N}
$-graded as an algebra. We have 
\begin{eqnarray*}
&&m_{B}^{a,b}\left( \mathrm{\Gamma }_{a}^{B}\otimes \mathrm{\Gamma }%
_{b}^{B}\right) \Delta _{B}^{a,b}\varphi _{B}^{a+b}=m_{B}^{a,b}\left( 
\mathrm{\Gamma }_{a}^{B}\varphi _{B}^{a}\otimes \mathrm{\Gamma }%
_{b}^{B}\varphi _{B}^{b}\right) \Delta _{T}^{a,b} \\
&\overset{\text{(\ref{form: S bar 2})}}{=}&m_{B}^{a,b}\left( \varphi _{B}^{a}%
\mathrm{\Gamma }_{a}^{T}\otimes \varphi _{B}^{b}\mathrm{\Gamma }%
_{b}^{T}\right) \Delta _{T}^{a,b} \\
&=&\varphi _{B}^{a+b}m_{T}^{a,b}\left( \mathrm{\Gamma }_{a}^{T}\otimes 
\mathrm{\Gamma }_{b}^{T}\right) \Delta _{T}^{a,b}\overset{\text{(*)}}{=}%
\varphi _{B}^{a+b}\mathrm{\Gamma }_{a+b}^{T}\overset{\text{(\ref{form: S bar
2})}}{=}\mathrm{\Gamma }_{a+b}^{B}\varphi _{B}^{a+b},
\end{eqnarray*}%
where in (*) we applied \cite[(1)]{Schauenburg0}. Since $\varphi _{B}^{a+b}$
is an epimorphism, we obtain (\ref{form: SS bar}).
\end{proof}

We now restrict our attention to connected graded braided bialgebras.

\begin{lemma}
Let $\left( B,m_{B},u_{B},\Delta _{B},\varepsilon _{B},c_{B}\right) $ be a
$0$-connected graded braided bialgebra. Then $\Delta _{B}^{n,0}\left( z\right)
=z\otimes 1$ and $\Delta _{B}^{0,n}=1\otimes z,$ for every $z\in B^{n}.$
Moreover%
\begin{equation}
\Delta _{B}^{n,1}m_{B}^{n,1}=\mathrm{Id}_{B^{n}\otimes B^{1}}+\left(
m_{B}^{n-1,1}\otimes B^{1}\right) \left( B^{n-1}\otimes c_{B}^{1,1}\right)
\left( \Delta _{B}^{n-1,1}\otimes B^{1}\right) .  \label{form: Delta Mult}
\end{equation}
\end{lemma}

\begin{proof}
The first assertion follows by \cite[Remark 1.16]{AMS-MM2}.

Since $B$ is a graded bialgebra we have%
\begin{eqnarray*}
&&\Delta _{B}^{n,1}m_{B}^{n,1}\overset{\text{(\ref{form: delta mult gen})}}{=%
}\left[ 
\begin{array}{c}
\left( m_{B}^{n,0}\otimes m_{B}^{0,1}\right) \left( B^{n}\otimes
c_{B}^{0,0}\otimes B^{1}\right) \left( \Delta _{B}^{n,0}\otimes \Delta
_{B}^{0,1}\right) + \\ 
+\left( m_{B}^{n-1,1}\otimes m_{B}^{1,0}\right) \left( B^{n-1}\otimes
c_{B}^{1,1}\otimes B^{0}\right) \left( \Delta _{B}^{n-1,1}\otimes \Delta
_{B}^{1,0}\right)%
\end{array}%
\right] \\
&=&\mathrm{Id}_{B^{n}\otimes B^{1}}+\left( m_{B}^{n-1,1}\otimes B^{1}\right)
\left( B^{n-1}\otimes c_{B}^{1,1}\right) \left( \Delta _{B}^{n-1,1}\otimes
B^{1}\right) .
\end{eqnarray*}
\end{proof}

Let $\left( B,m_{B},u_{B},\Delta _{B},\varepsilon _{B},c_{B}\right) $ be a
$0$-connected graded braided bialgebra. Next aim is to provide an upper bound of 
$\mathrm{sdeg}\left( B\right) $ under suitable assumptions on $c_{B}^{1,1}$.
These results will be needed in the last sections of the paper.

\begin{definition}
\label{def: regular element}For every $n\geq 1$, we set 
\begin{equation*}
\left( n\right) _{X}:=\frac{X^{n}-1}{X-1}=1+X+\cdots +X^{n-1},\quad \left(
n\right) _{X}!:=\left( 1\right) _{X}\cdot \left( 2\right) _{X}\cdot \cdots
\cdot \left( n\right) _{X}.
\end{equation*}%
An element $q\in K$ is called \textbf{regular} whenever $\left( n\right)
_{q}\neq 0$, for all $n\geq 2.$
\end{definition}

\begin{theorem}
\label{teo: sdeg(B)=0}Let $\left( B,m_{B},u_{B},\Delta _{B},\varepsilon
_{B},c_{B}\right) $ be a $0$-connected graded braided bialgebra. Assume that

\begin{itemize}
\item $B$ is strongly $%
\mathbb{N}
$-graded as an algebra;

\item $m_{B}^{1,1}\left( c_{B}^{1,1}-q\mathrm{Id}_{B^{2}}\right) =0$ for
some $q\in K$ (e.g. $c_{B}^{1,1}$ has minimal polynomial of degree $1$);

\item $q$ is regular.
\end{itemize}

Then the canonical projection $\pi _{S}:B\rightarrow S\left( B\right) $ is
an isomorphism so that $\mathrm{sdeg}\left( B\right) =0$.
\end{theorem}

\begin{proof}
By using the definition of $\mathrm{\Gamma }_{n,1}^{B},$ (\ref{form: Delta
Mult}), (\ref{gr1}), and the condition $m_{B}^{1,1}c_{B}^{1,1}=qm_{B}^{1,1},$
we arrive at%
\begin{equation}
\mathrm{\Gamma }_{n,1}^{B}m_{B}^{n,1}=m_{B}^{n,1}+qm_{B}^{n,1}\left( \mathrm{%
\Gamma }_{n-1,1}^{B}\otimes B^{1}\right) .  \label{form: grado 1}
\end{equation}%
Inductively, using (\ref{form: grado 1}) and the surjectivity of $%
m_{B}^{n,1} $ we obtain%
\begin{equation}
\mathrm{\Gamma }_{n,1}^{B}=\left( n+1\right) _{q}\mathrm{Id}_{B^{n+1}},\text{
for every }n\geq 1\text{.}  \label{form: grado 1 reduced}
\end{equation}%
%
%
%
%
%
%
%
%
%
%
%
%
%
%
%
%
%
%
%
%
%
%
%
%
%
%
%
%
%
%
%
%
%
%
%
%
%
%
%
%
%
%
%
%
%
%
%
%
%
%
%
%
%
%
%
%
%
%
%
%
%
%
%
%
%
%
Since $q$ is regular i.e. $\left( n\right) _{q}\neq 0,\forall n\geq 2,$ we
get that $\mathrm{\Gamma }_{n,1}^{B}\ $is bijective for every $n\geq 1$.
Since, by definition, $\mathrm{\Gamma }_{n,1}^{B}:=m_{B}^{n,1}\Delta
_{B}^{n,1},$ we find that $\Delta _{B}^{n,1}$ is injective for every $n\geq
1 $ and hence for every $n\in 
\mathbb{N}
$ ($\left( \varepsilon _{B}^{0}\otimes B^{1}\right) \Delta
_{B}^{0,1}=l_{B^{1}}$ is bijective). By Theorem \ref{teo: strongness}, the
canonical projection $\pi _{S}:B\rightarrow S\left( B\right) $ is an
isomorphism so that $\mathrm{sdeg}\left( B\right) =0.$
\end{proof}

\begin{theorem}
\label{teo: sdeg(B)=n-1}Let $\left( B,m_{B},u_{B},\Delta _{B},\varepsilon
_{B},c_{B}\right) $ be a connected graded braided bialgebra. Assume that

\begin{itemize}
\item $m_{B}^{1,1}$ is an isomorphism.

\item $c_{B}^{1,1}$ has minimal polynomial $f\left( X\right) \in K\left[ X%
\right] $ where

\item $f\left( -1\right) =0$ i.e. $f\left( X\right) =\left( X+1\right)
h\left( X\right) $ for some $h\in K\left[ X\right] .$
\end{itemize}

Then $\mathrm{sdeg}\left( B\right) =\mathrm{sdeg}\left( S\left( B\right)
\right) +1$ and $m_{S\left( B\right) }^{1,1}h\left( c_{S\left( B\right)
}^{1,1}\right) =0.$
\end{theorem}

\begin{proof}
Note that, if $\mathrm{sdeg}\left( B\right) =0,$ then $\Delta _{B}^{1,1}$ is
injective and hence $c_{B}^{1,1}+\mathrm{Id}_{B^{2}}=\Delta
_{B}^{1,1}m_{B}^{1,1}$ is injective too. In this case $f\left(
c_{B}^{1,1}\right) =0$ imply $h\left( c_{B}^{1,1}\right) =0$ so that $f$ is
not the minimal polynomial for $c_{B}^{1,1},$ a contradiction. We have so
proved that $\mathrm{sdeg}\left( B\right) \neq 0.$ Thus $\mathrm{sdeg}\left(
B\right) =\mathrm{sdeg}\left( S\left( B\right) \right) -1$.

By (\ref{form: braiding S(B)}), we have $c_{S}^{1,1}\left( \pi
_{S}^{1}\otimes \pi _{S}^{1}\right) =\left( \pi _{S}^{1}\otimes \pi
_{S}^{1}\right) c_{B}^{1,1}.$ Since $\pi _{S}^{1}$ is bijective, we infer
that $c_{S}^{1,1}$ has minimal polynomial $f.$ Thus $0=f\left(
c_{S}^{1,1}\right) =\left( c_{S}^{1,1}+\mathrm{Id}_{S^{2}}\right) h\left(
c_{S}^{1,1}\right) =\Delta _{S}^{1,1}m_{S}^{1,1}h\left( c_{S}^{1,1}\right) .$
By Proposition \ref{pro: Delta S[n]}, $\Delta _{S}^{1,1}$ is injective
whence $m_{S}^{1,1}h\left( c_{S}^{1,1}\right) =0.$
\end{proof}

\section{Braided Lie algebras}

Although many result in this section could be obtained for a general graded
braided bialgebra $B$, we focus here our attention on the case $B=T\left(
V,c\right) ,$ the tensor algebra of a braided vector space $\left(
V,c\right) .$

\begin{claim}
\label{claim: Gamma}Let $\left( V,c\right) $ be a braided vector space and
let $T=T\left( V,c\right) $. By the universal property of the tensor algebra
there is a unique algebra homomorphism 
\begin{equation*}
\mathrm{\Gamma }^{T}:T\left( V,c\right) \rightarrow T^{c}\left( V,c\right)
\end{equation*}%
such that $\mathrm{\Gamma }_{\mid V}^{T}=\mathrm{Id}_{V},$ where $%
T^{c}\left( V,c\right) $ denotes the quantum shuffle algebra (see \ref%
{claim: type one}). This is a morphism of graded braided bialgebras. By \cite%
[page 2815]{Schauenburg0} (see also \cite[page 25-26]{AG}), it is clear that
the $n$-th graded component of $\mathrm{\Gamma }^{T}$ is the map $\mathrm{%
\Gamma }_{n}^{T}:V^{\otimes n}\rightarrow V^{\otimes n}$ in the sense of
Definition \ref{def: Gamma} (see Remark \ref{rem: Gamma}). The bialgebra of
Type one generated by $V$ over $K$ (or Nichols algebra) is by definition%
\begin{equation*}
\emph{B}\left( V,c\right) =\mathrm{Im}\left( \mathrm{\Gamma }^{T}\right)
\simeq \frac{T\left( V,c\right) }{\ker \left( \mathrm{\Gamma }^{T}\right) }.
\end{equation*}%
We set%
\begin{equation*}
E_{n}\left( V,c\right) :=E_{n}\left( T\left( V,c\right) \right) \qquad \text{%
and}\qquad E\left( V,c\right) :=E\left( T\left( V,c\right) \right)
\end{equation*}%
where we used the notations of Definition \ref{def: E(C)}.
\end{claim}

\begin{definition}
A \textbf{(braided) bracket} on a braided vector space $\left( V,c\right) $
is a $K$-linear map $b:E\left( V,c\right) \rightarrow V$ such that 
\begin{equation}
c\left( b\otimes V\right) =\left( V\otimes b\right) c_{E\left( V,c\right)
,V}\qquad c\left( V\otimes b\right) =\left( b\otimes V\right) c_{V,E\left(
V,c\right) }.  \label{form: c-bracket}
\end{equation}%
The restriction of $b$ to $E_{t}\left( V,c\right) $ will be denoted by $%
b^{t}:E_{t}\left( V,c\right) \rightarrow V.$ A \textbf{morphism of (braided)
brackets} $f:\left( V,c_{V,V},b_{V}\right) \rightarrow \left(
W,c_{W,W},b_{W}\right) $ is a morphism of braided vector spaces $f:\left(
V,c_{V,V}\right) \rightarrow \left( W,c_{W,W}\right) $ such that%
\begin{equation}
f\circ b_{V}^{t}=b_{W}^{t}\circ E_{t}\left( f\right) .
\label{form: morph brackets}
\end{equation}
\end{definition}

\begin{lemma}
\label{lem: bracket VS bracket}Let $\left( V,c\right) $ be a braided vector
space. The following assertions are equivalent for a $K$-linear map $%
b:E\left( V,c\right) \rightarrow V:$

\begin{enumerate}
\item[$\left( i\right) $] $b$ is a bracket for the graded braided bialgebra $%
T\left( V,c\right) $;

\item[$\left( ii\right) $] $b$ is a bracket on $\left( V,c\right) $.
\end{enumerate}
\end{lemma}

\begin{proof}
$\left( i\right) \Rightarrow \left( ii\right) $ It is trivial.

$\left( ii\right) \Rightarrow \left( i\right) $ Set $T:=T\left( V,c\right) .$
Let us prove, by induction on $t\geq 1,$ that $b$ fulfills (\ref{form:
Bracket}) i.e. that 
\begin{equation*}
c_{T}^{1,t}\left( b\otimes V^{\otimes t}\right) =\left( V^{\otimes t}\otimes
b\right) c_{E\left( V,c\right) ,V^{\otimes t}}\qquad c_{T}^{t,1}\left(
V^{\otimes t}\otimes b\right) =\left( b\otimes V^{\otimes t}\right)
c_{V^{\otimes t},E\left( V,c\right) }.
\end{equation*}%
For $t=1$ there is nothing to prove as $c_{T}^{1,1}=c$. Assume that the
formulas are true for $t-1.$ Then, by construction of $c_{T}$, we have 
\begin{eqnarray*}
c_{T}^{1,t}\left( b\otimes V^{\otimes t}\right) &=&\left( V^{\otimes
t-1}\otimes c_{T}^{1,1}\right) \left( c_{T}^{1,t-1}\otimes V\right) \left(
b\otimes V^{\otimes t-1}\otimes V\right) \\
&=&\left( V^{\otimes t-1}\otimes c_{T}^{1,1}\right) \left( V^{\otimes
t-1}\otimes b\otimes V\right) \left( c_{E\left( V,c\right) ,V^{\otimes
t-1}}\otimes V\right) \\
&=&\left( V^{\otimes t-1}\otimes V\otimes b\right) \left( V^{\otimes
t-1}\otimes c_{E\left( V,c\right) ,V}\right) \left( c_{E\left( V,c\right)
,V^{\otimes t-1}}\otimes V\right)\\&=&\left( V^{\otimes t}\otimes b\right)
c_{E\left( V,c\right) ,V^{\otimes t}}
\end{eqnarray*}%
and similarly we get the second equality. It remains to prove that $b$
fulfills (\ref{form: Bracket2}) i.e. that $\Delta _{T}^{1,0}b=\left(
b\otimes K\right) \Delta _{T}^{t,0}$ and $\Delta _{T}^{0,1}b=\left( K\otimes
b\right) \Delta _{T}^{0,t}$ but these equalities are trivially true by the
definition of $\Delta _{T}.$
\end{proof}

\begin{definition}
Let $b$ be a bracket on a braided vector space $\left( V,c\right) $. In view
of Lemma \ref{lem: bracket VS bracket}, it makes sense to define \textbf{the
universal enveloping algebra of }$\left( V,c,b\right) $ to be 
\begin{equation}
U\left( V,c,b\right) :=U\left( T\left( V,c\right) ,b\right) =\frac{T\left(
V,c\right) }{\left( \left( \mathrm{Id}-b\right) \left[ E\left( V,c\right) %
\right] \right) }  \label{form: U(V,c,b)}
\end{equation}%
which, by Theorem \ref{teo: U bialgebra}, is a braided bialgebra quotient of
the tensor algebra $T\left( V,c\right) .$ The \textbf{symmetric algebra of }$%
\left( V,c\right) $ is then 
\begin{equation*}
S\left( V,c\right) :=U\left( V,c,0\right) =U\left( T\left( V,c\right)
,0\right) =S\left( T\left( V,c\right) \right)
\end{equation*}%
(see (\ref{form: S(B)})).

We will denote by $\pi _{U}:T\left( V,c\right) \rightarrow U\left(
V,c,b\right) $ and $\pi _{S}:T\left( V,c\right) \rightarrow S\left(
V,c\right) $ the canonical projections.
\end{definition}

\begin{remark}
Since $\pi _{U}$ (resp. $\pi _{S}$) is an epimorphism, then $U\left(
V,c,b\right) $ (resp. $S\left( V,c\right) $) is a connected coalgebra (cf. 
\cite[Corollary 5.3.5]{Montgomery}).
\end{remark}

\begin{definition}
\label{def: c-Lie alg}We say that $\left( V,c,b\right) $ is a \textbf{%
braided Lie algebra} whenever

\begin{itemize}
\item $\left( V,c\right) $ is a braided vector space;

\item $b:E\left( V,c\right) \rightarrow V$ is a bracket on $\left(
V,c\right) $;

\item the canonical $K$-linear map $i_{U}:V\rightarrow U\left( V,c,b\right) $
is injective i.e. $V\cap \ker \left( \pi _{U}\right) =\ker \left(
i_{U}\right) =0.$
\end{itemize}
\end{definition}

Injectivity of $i_{U}$ plays here the role of an implicit Jacobi identity.
On the other hand, antisymmetry of the bracket is encoded in the choice of
the domain of $b$. This becomes clearer when $c$ is of Hecke type as in
Theorem \ref{teo: Hecke Lie}. There it is shown that $b$ can be substituted
with the $c$-antisymmetric map $\left[ -\right] _{b}:V\otimes V\rightarrow V$
defined by setting $\left[ z\right] _{b}=b\left( c-q\mathrm{Id}_{V\otimes
V}\right) \left( z\right) $, for every $z\in V\otimes V$. \medskip\newline
A first example of braided Lie algebra is given in the following
Proposition. A large class of examples with not necessarily trivial brackets
will be included in Proposition \ref{pro: alg is Lie} or Theorem \ref{teo:
univ U}.

\begin{proposition}
\label{pro: (V,c,0) Lie}Let $\left( V,c\right) $ is a braided vector space.
Then $\left( V,c,0\right) $ braided Lie algebra.
\end{proposition}

\begin{proof}
We already observed that $0$ is always a braided bracket for $\left(
V,c\right) $ whence it remains to prove that $i_{U}$ is injective. Set $A:=%
\frac{T\left( V,c\right) }{\left( V\otimes V\right) }$ where $\left(
V\otimes V\right) $ is the two sided ideal of $T\left( V,c\right) $
generated by $V\otimes V$. Denote by $\pi :T\left( V,c\right) \rightarrow A$
and $\pi _{S}:T\left( V,c\right) \rightarrow S\left( V,c\right) $ the
canonical projections. Since $E\left( V,c\right) \subseteq \left( V\otimes
V\right) $, then $\pi :T\left( V,c\right) \rightarrow A$ quotients to an
algebra homomorphism $\overline{\pi }:S\left( V,c\right) =U\left(
V,c,0\right) \rightarrow A\ $such that $\overline{\pi }\circ \pi _{S}=\pi $.
Hence $\overline{\pi }\circ i_{U}=\overline{\pi }\circ \pi _{S}\circ
i_{T}=\pi \circ i_{T}$ where $i_{T}:V\rightarrow T\left( V,c\right) $ is the
canonical injection. Note that as a vector space $A\simeq K\oplus V$ and
through this isomorphism, $\pi \circ i_{T}$ identifies with the canonical
injection of $V$ in $K\oplus V$ whence $i_{U}$ is injective too.
\end{proof}

\begin{proposition}
\label{pro: alg is Lie}Let $\left( A,c_{A}\right) $ be a braided algebra.
Then $\left( A,c_{A},b_{A}\right) $ is a braided Lie algebra, where $%
b_{A}\left( z\right) :=m_{A}^{t-1}\left( z\right) ,$ for every $z\in
E_{t}\left( A,c_{A}\right) .$
\end{proposition}

\begin{proof}
Let us check that $b_{A}$ is a $c_{A}$-bracket on the braided vector space $%
\left( A,c_{A}\right) $. For every $z\in E_{t}\left( A,c_{A}\right) $ and $%
x\in A,$ we have%
\begin{eqnarray*}
c_{A}\left( b_{A}^{t}\otimes A\right) \left( z\otimes x\right)
&=&c_{A}\left( m_{A}^{t-1}\otimes A\right) \left( z\otimes x\right) \\
&=&\left( A\otimes m_{A}^{t-1}\right) c_{T}^{t,1}\left( z\otimes x\right)
=\left( A\otimes b_{A}^{t}\right) c_{E_{t}\left( A,c_{A}\right) ,A}\left(
z\otimes x\right) .
\end{eqnarray*}%
so that the left-hand side of (\ref{form: c-bracket}) is proved. Similarly
one proves also the right-hand side. By the universal property of $T\left(
A,c_{A}\right) $ (see \cite[Theorem 1.17]{AMS-MM2}), $Id_{A}$ can be lifted
to a unique braided algebra homomorphism $\varphi _{A}:T\left(
A,c_{A}\right) \rightarrow A.$

For every $z\in E_{t}\left( A,c_{A}\right) ,$ we have $\varphi
_{A}b_{A}\left( z\right) =\varphi _{A}m_{A}^{t-1}\left( z\right)
=m_{A}^{t-1}\left( z\right) =\varphi _{A}\left( z\right) $ so that $\varphi
_{A}$ quotients to a morphism $\overline{\varphi _{A}}:U\left(
A,c_{A},b_{A}\right) \rightarrow A.$ If we denote by $i_{U}:P\rightarrow
U\left( A,c_{A},b_{A}\right) $ the canonical map, we get that $\overline{%
\varphi _{A}}\circ i_{U}=Id_{A}.$ In particular $i_{U}$ is injective whence $%
\left( A,c_{A},b_{A}\right) $ is a braided Lie algebra.
\end{proof}

\begin{lemma}
\label{lem: subLie}Let $\left( V,c_{V},b_{V}\right) $ a braided Lie algebra
and let $h:\left( W,c_{W}\right) \rightarrow \left( V,c_{V,V}\right) $ be an
injective morphism of braided vector spaces. If there exists a $K$-linear
map $b_{W}:E\left( W,c_{W}\right) \rightarrow W$ such that $h\circ
b_{W}=b_{V}\circ E\left( h\right) ,$ then $\left( W,c_{W},b_{W}\right) $ is
a braided Lie algebra too.
\end{lemma}

\begin{proof}
We have%
\begin{eqnarray*}
\left( h\otimes h\right) c_{W}\left( b_{W}\otimes W\right) &=&c_{V}\left(
h\otimes h\right) \left( b_{W}\otimes W\right) =c_{V}\left( b_{V}\otimes
V\right) \left( E\left( h\right) \otimes h\right) \\
&\overset{\text{(\ref{form: c-bracket})}}{=}& \left( V\otimes b_{V}\right)
c_{E\left( V,c_{V}\right) ,V}\left( E\left( h\right) \otimes h\right)
=\left( V\otimes b_{V}\right) \left( h\otimes E\left( h\right) \right)
c_{E\left( W,c_{W}\right) ,W} \\
&=&\left( h\otimes h\right) \left( W\otimes b_{W}\right) c_{E\left(
W,c_{W}\right) ,W}.
\end{eqnarray*}%
Since $h$ is injective, we deduce that $c_{W}\left( b_{W}\otimes W\right)
=\left( W\otimes b_{W}\right) c_{E\left( W,c_{W}\right) ,W}$ so that the
left-hand side of (\ref{form: c-bracket}) is proved for $\left(
W,c_{W}\right) $. Similarly one proves also the right-hand side. If we
denote by $i_{U}^{W}:W\rightarrow U\left( W,c_{W},b_{W}\right) $ and $%
i_{U}^{V}:V\rightarrow U\left( V,c_{V},b_{V}\right) $ the canonical maps, by
the universal property of the universal enveloping algebra, there is a
unique algebra homomorphism $U\left( h\right) :U\left( W,c_{W}\right)
\rightarrow U\left( V,c_{V,V}\right) $ such that $U\left( h\right) \circ
i_{U}^{W}=i_{U}^{V}\circ h.$ Since both $i_{U}^{V}$ and $h$ are injective,
we deduce that $i_{U}^{W}$ is injective too. Therefore $\left(
W,c_{W},b_{W}\right) $ is a braided Lie algebra.
\end{proof}

\begin{theorem}
\label{teo: univ U}Let $\left( A,c_{A}\right) $ be a connected braided
bialgebra. Let $P=P\left( A\right) $ be the space of primitive elements of $%
A.$ Then%
\begin{equation*}
c_{A}\left( P\otimes P\right) \subseteq P\otimes P.
\end{equation*}%
Let $c_{P}:P\otimes P\rightarrow P\otimes P$ be the restriction of $c_{A}$
to $P\otimes P.$ Then

\begin{enumerate}
\item[i)] $m_{A}^{t-1}\left( E_{t}\left( P,c_{P}\right) \right) \subseteq P$
for every $t\geq 2,$ so we can define $b_{P}:E\left( P,c_{P}\right)
\rightarrow P$ by $b_{P}\left( z\right) :=m_{A}^{t-1}\left( z\right) ,$ for
every $z\in E_{t}\left( P,c_{P}\right) .$

\item[ii)] $\left( P,c_{P},b_{P}\right) $ is a braided Lie algebra.

\item[iii)] Every morphism of braided brackets $f:\left(
V,c_{V,V},b_{V}\right) \rightarrow \left( P,c_{P},b_{P}\right) $ can be
lifted to a morphism of braided bialgebras $\overline{f}:U\left(
V,c_{V,V},b_{V}\right) \rightarrow A.$
\end{enumerate}
\end{theorem}

\begin{proof}
The first assertion follows by \cite[Remark 1.12]{AMS-MM2}. Set $T:=T\left(
P,c_{P}\right) .$

i) Since $T$ is connected, by Lemma \ref{lem: Nichols}, we have $P\left(
T\right) =P_{1}\left( T\right) =T^{1}\oplus \left[ \oplus _{n\geq
2}E_{n}\left( P,c_{P}\right) \right] .$ Thus $E_{t}\left( P,c_{P}\right)
\subseteq P\left( T\right) $ for every $t\geq 2$. Let $\varphi
_{A}:T\rightarrow A$ be the natural braided bialgebra homomorphism arising
from the universal property of the tensor algebra. Then, for every $t\geq 2$%
, we have $m_{A}^{t-1}\left( E_{t}\left( P,c_{P}\right) \right) =\varphi
_{A}\left( E_{t}\left( P,c_{P}\right) \right) \subseteq \varphi _{A}\left(
P\left( T\right) \right) \subseteq P\left( A\right) =P.$

ii) Let $b_{P}$ be defined as in the statement and denote by $%
i_{A}:P\rightarrow A$ the canonical inclusion. Then $i_{A}\circ
b_{P}=b_{A}\circ E\left( i_{A}\right) ,$ where $b_{A}$ is the map defined in
Proposition \ref{pro: alg is Lie}. By this proposition and Lemma \ref{lem:
subLie} we have that $\left( P,c_{P},b_{P}\right) $ is a braided Lie algebra.

iii) By the universal property of $T\left( V,c_{V,V}\right) $ (see \cite[%
Theorem 1.17]{AMS-MM2}), $f$ can be lifted to a unique braided bialgebra
homomorphism $f^{\prime }:T\left( V,c_{V,V}\right) \rightarrow A.$ We have%
\begin{equation*}
f^{\prime }b_{V}\left( z\right) =fb_{V}^{t}\left( z\right) \overset{\text{(%
\ref{form: morph brackets})}}{=}b_{P}^{t}E_{t}\left( f\right) \left(
z\right) =m_{A}^{t-1}f^{\otimes t}\left( z\right) =f^{\prime }\left( z\right)
\end{equation*}%
for every $z\in E_{t}\left( V,c_{V,V}\right) .$ Therefore $f^{\prime }$
quotients to a morphism $\overline{f}:U\left( V,c_{V,V},b_{V}\right)
\rightarrow A.$
\end{proof}

\begin{definition}
With the same assumptions and notations of Theorem \ref{teo: univ U}, $%
\left( P,c_{P},b_{P}\right) $ will be called the \textbf{infinitesimal
braided Lie algebra of }$A.$
\end{definition}

\begin{theorem}
\label{teo: univ S}Let $\left( A,c_{A}\right) $ be a connected braided
bialgebra and let $\left( P,c_{P},b_{P}\right) $ be its infinitesimal
braided Lie algebra. Then, every morphism of braided vector spaces $f:\left(
V,c_{V,V}\right) \rightarrow \left( P,c_{P}\right) $ such that $%
m_{A}^{t-1}E_{t}\left( f\right) =0,$ for every $t\geq 2,$ can be lifted to a
morphism of braided bialgebras $\overline{f}:S\left( V,c_{V,V}\right)
\rightarrow A$ which respects the gradings on $S\left( V,c_{V,V}\right) $
and $A$ whenever $A$ is graded.
\end{theorem}

\begin{proof}
By Proposition \ref{pro: (V,c,0) Lie}, we know that $\left(
V,c_{V,V},b_{V}\right) $ is a braided Lie algebra for $b_{V}=0$. Let us
prove that $f:\left( V,c_{V,V},b_{V}\right) \rightarrow \left(
P,c_{P},b_{P}\right) $ is a morphism of braided Lie algebras. Since $%
b_{P}^{t}$ acts like $m_{A}^{t-1},$ by hypothesis we get $b_{P}^{t}\circ
E_{t}\left( f\right) =0=f\circ b_{V}^{t}.$ Hence, we can apply Theorem \ref%
{teo: univ U}.
\end{proof}

\begin{claim}
\label{claim: graded}Let $\left( V,c\right) $ be a braided vector space and
let $b:E\left( V,c\right) \rightarrow V$ be a $c$-bracket. Set%
\begin{equation*}
T:=T\left( V,c\right) \qquad \text{and}\qquad U:=U\left( V,c,b\right) .
\end{equation*}%
Let $\pi _{U}:T\rightarrow U$ be the canonical projection. By construction $%
\pi _{U}$ is a morphism of braided bialgebras. Mimicking \cite[4.11]{AMS-MM}%
, set $T_{\left( n\right) }:=\oplus _{0\leq t\leq n}V^{\otimes t}$ and 
\begin{equation*}
U_{n}^{\prime }:=\pi _{U}\left( T_{\left( n\right) }\right) .
\end{equation*}%
Then $\left( U_{n}^{\prime }\right) _{n\in 
\mathbb{N}
}$ is both an algebra and a coalgebra filtration on $U$ which is called the 
\emph{standard filtration on }$U$. Note that this filtration is not the
coradical filtration $\left( U_{n}\right) _{n\in \mathbb{N}}$ in general. Still one has $U_{n}^{\prime }=\pi _{U}\left( T_{\left(
n\right) }\right) \subseteq \pi _{U}\left( T_{n}\right) \subseteq U_{n}$
where $T_{n}$ and $U_{n}$ denote the $n$-th terms of the coradical
filtration of $T$ and $U$ respectively. Denote by%
\begin{equation*}
\mathrm{gr}^{\prime }\left( U\right) :=\oplus _{n\in \mathbb{N}}\frac{U_{n}^{\prime }}{U_{n-1}^{\prime }}.
\end{equation*}%
the graded coalgebra associated to the standard filtration (see \cite[page
228]{Sw}).

If $b=0$, then $S\left( V,c\right) =U\left( V,c,0\right) $ is a graded
bialgebra $S\left( V,c\right) =\oplus _{n\in \mathbb{N}}S^{n}\left( V,c\right) .$ The standard filtration on $S\left( V,c\right) $
is the filtration associated to this grading.
\end{claim}

The following result is inspired to \cite[Proposition 4.19]{AMS-MM}.

\begin{proposition}
\label{pro: theta}Let $\left( V,c\right) $ be a braided vector space and let 
$b:E\left( V,c\right) \rightarrow V$ be a $c$-bracket. Then $\mathrm{gr}%
^{\prime }\left( U\left( V,c,b\right) \right) $ is a graded braided
bialgebra and there is a canonical morphism of graded braided bialgebras $%
\theta :S\left( V,c\right) \rightarrow \mathrm{gr}^{\prime }\left( U\left(
V,c,b\right) \right) $ which is surjective and lifts the map $\theta
_{1}:V\rightarrow U_{1}^{\prime }/U_{0}^{\prime }:v\mapsto \pi _{U}\left(
v\right) +U_{0}^{\prime }.$
\end{proposition}

\begin{proof}
Set $T:=T\left( V,c\right) $, $U:=U\left( V,c,b\right) $, $G:=\mathrm{gr}%
^{\prime }\left( U\left( V,c,b\right) \right) ,$ $G^{n}:=U_{n}^{\prime
}/U_{n-1}^{\prime }$ and let $p_{n}:U_{n}^{\prime }\rightarrow G^{n}$ be the
canonical projection, for every $n\in 
\mathbb{N}
$. We have%
\begin{eqnarray*}
c_{U}\left( U_{a}^{\prime }\otimes U_{b}^{\prime }\right) &=&c_{U}\left( \pi
_{U}\otimes \pi _{U}\right) \left( T_{\left( a\right) }\otimes T_{\left(
b\right) }\right) =\left( \pi _{U}\otimes \pi _{U}\right) c_{T}\left(
T_{\left( a\right) }\otimes T_{\left( b\right) }\right) \\
&\subseteq &\left( \pi _{U}\otimes \pi _{U}\right) \left( T_{\left( b\right)
}\otimes T_{\left( a\right) }\right) =U_{b}^{\prime }\otimes U_{a}^{\prime }.
\end{eqnarray*}%
Hence $c_{U}$ induces a braiding $c_{G}:G\otimes G\rightarrow G\otimes G.$
Now, we have already observed that $G$ carries a coalgebra structure $\left(
G,\Delta _{G},\varepsilon _{G}\right) $. Since $\left( U_{n}^{\prime
}\right) _{n\in 
\mathbb{N}
}$ is also an algebra filtration on $U$, the multiplication of $U$ induces,
for every $a,b\in 
\mathbb{N}
$, maps $m_{G}^{a,b}:G^{a}\otimes G^{b}\rightarrow G^{a+b}$ and the unit of $%
U$ induces a map $u_{G}^{0}:K\rightarrow G^{0}$ such that (\ref{gr1}) and (%
\ref{gr2}) hold for $"A"=G.$ Thus there is a unique map $m_{G}:G\otimes
G\rightarrow G$ such that $m_{G}\left( i_{G}^{a}\otimes i_{G}^{b}\right)
=i_{G}^{a+b}m_{G}^{a,b},$ for every $a,b\in 
\mathbb{N}
$, where $i_{G}^{t}:G^{t}\rightarrow G$ is the canonical map. Moreover $%
\left( G,m_{G},u_{G}=i_{G}^{0}u_{G}^{0}\right) $ is a graded algebra (see
e.g. \cite[Proposition 3.4]{AM- Type One}). It is straightforward to prove
that $\left( G,m_{G},u_{G},\Delta _{G},\varepsilon _{G},c_{G}\right) $ is
indeed a graded braided bialgebra.

Set $P:=P\left( G\right) .$ Since $G$ is a connected graded coalgebra, it is
clear that $\mathrm{Im}\left( \theta _{1}\right) \subseteq G^{1}\subseteq P.$
Moreover $\theta _{1}:\left( V,c\right) \rightarrow \left( P,c_{P}\right) $
is a morphism of braided vector spaces as, for every $u,v\in V,$ we have%
\begin{eqnarray*}
c_{G}\left( \theta _{1}\otimes \theta _{1}\right) \left( u\otimes v\right)
&=&c_{G}\left[ \left( \pi _{U}\left( u\right) +U_{0}^{\prime }\right)
\otimes \left( \pi _{U}\left( v\right) +U_{0}^{\prime }\right) \right] \\
&=&\left( p_{1}\otimes p_{1}\right) c_{U}\left( \pi _{U}\left( u\right)
\otimes \pi _{U}\left( v\right) \right) \\
&=&\left( p_{1}\otimes p_{1}\right) c_{U}\left( \pi _{U}\otimes \pi
_{U}\right) \left( u\otimes v\right) \\
&=&\left( p_{1}\otimes p_{1}\right) \left( \pi _{U}\otimes \pi _{U}\right)
c_{T}\left( u\otimes v\right) =\left( \theta _{1}\otimes \theta _{1}\right)
c\left( u\otimes v\right) .
\end{eqnarray*}%
Let us prove that $m_{G}^{t-1}E_{t}\left( \theta _{1}\right) =0,$ for every $%
t\geq 2.$ For every $t\geq 2$ and $z\in E_{t}\left( V,c\right) ,$ we have 
\begin{eqnarray*}
&&m_{G}^{t-1}E_{t}\left( \theta _{1}\right) \left( z\right)
=m_{G}^{t-1}\theta _{1}^{\otimes t}\left( z\right) =p_{t}m_{U}^{t-1}\pi
_{U}^{\otimes t}\left( z\right) =p_{t}\pi _{U}m_{T}^{t-1}\left( z\right) \\
&=&p_{t}\pi _{U}\left( z\right) =\pi _{U}\left( z\right) +U_{t-1}^{\prime }%
\overset{\text{def. }U}{=}\pi _{U}b^{t}\left( z\right) +U_{t-1}^{\prime }=0
\end{eqnarray*}%
so that $m_{G}^{t-1}\circ E_{t}\left( \theta _{1}\right) =0.$ By Theorem \ref%
{teo: univ S}, there is a canonical morphism of graded braided bialgebras $%
\theta :S\left( V,c\right) \rightarrow G$ lifting $\theta _{1}.$ By
construction $\theta $ lifts the canonical map $T\left( V,c_{V,V}\right)
\rightarrow G$ which is surjective as $\theta _{1}$ is surjective and $G$ is
generated as an $K$-algebra by $G^{1}.$ Thus $\theta $ is surjective too.
\end{proof}

\begin{definition}
\label{def: PBW}Let $\left( V,c\right) $ be a braided vector space and let $%
b:E\left( V,c\right) \rightarrow V$ be a $c$-bracket. Following \cite[%
Definition, page 316]{Braverman-Gaitsgory}, we will say that $U\left(
V,c,b\right) $ is \textbf{Poincar\'{e}-Birkhoff-Witt (PBW) type} whenever
the projection $\theta :S\left( V,c\right) \rightarrow \mathrm{gr}^{\prime
}\left( U\left( V,c,b\right) \right) $ of Proposition \ref{pro: theta} is an
isomorphism (compare with \cite[page 92]{Humphreys} for justifying this
terminology).
\end{definition}

\section{Strongness degree of a braided vector space}

This section mainly concerns the characterization of braided vector spaces
of strongness degree not greater then $1$ (see the following definition). As
we will see in the sequel, braided vector spaces with braiding of Hecke type
of regular mark $q$ (Theorem \ref{teo: sdeg Hecke}) and braided vector
spaces whose Nichols algebra is a quadratic algebra (Theorem \ref{teo:
quadratic => sdeg 1}) are examples of such a braided vector space.

\begin{definition}
\label{def: sdeg(V,c)}Let $\left( V,c\right) $ be braided vector space. The 
\textbf{strongness degree} of $\left( V,c\right) $ is defined to be%
\begin{equation*}
\mathrm{sdeg}\left( V,c\right) :=\mathrm{sdeg}\left( T\left( V,c\right)
\right) .
\end{equation*}
\end{definition}

The following result shows how the investigation of braided vector spaces
with finite strongness degree fits into the classification of finite
dimensional Hopf algebras problem. In fact, in this direction, one of the
steps of the celebrated lifting method by Andruskiewitsch and Schneider (see
e.g. \cite{AS- Lifting}) is the classification of finite dimensional Nichols
algebras.

\begin{theorem}
Let $\left( V,c\right) $ be a braided vector space such that the Nicholas
algebra $\mathcal{B}\left( V,c\right) $ is finite dimensional. Then $\mathrm{%
sdeg}\left( V,c\right) \leq \deg \left( \mathfrak{h}\left( \mathcal{B}\left(
V,c\right) \right) \right) -1\leq \dim _{K}\mathcal{B}\left( V,c\right) -1.$
\end{theorem}

\begin{proof}
Observe that the tensor algebra $T:=T_{K}\left( V,c\right) $ is strongly $%
\mathbb{N}
$-graded as an algebra. Since $\mathcal{B}\left( V,c\right) $ identifies
with $T^{0}\left[ T^{1}\right] $ i.e. the braided bialgebra of type
associated to $T^{0}$ and $T^{1}$, the conclusion follows by applying
Corollary \ref{coro: f.g. Type1} to the case $B=T.$
\end{proof}


\begin{theorem}
\label{teo: magnum}Let $\left( V,c\right) $ be a braided vector space and
let $b:E\left( V,c\right) \rightarrow V$ be a braided bracket. Assume that $%
\mathrm{sdeg}\left( V,c\right) \leq 1$. Then, the following assertions are
equivalent.

\begin{enumerate}
\item[(i)] $\left( V,c,b\right) $ is a braided Lie algebra in the sense of
Definition \ref{def: c-Lie alg}.

\item[(ii)] $V\cap \ker \left( \pi _{U}\right) =0.$

\item[(iii)] The map $\theta _{1}:V\rightarrow U_{1}^{\prime }/U_{0}^{\prime
}:v\mapsto \pi _{U}\left( v\right) +U_{0}^{\prime }$ is injective.

\item[(iv)] $U\left( V,c,b\right) $ is of PBW type in the sense of
Definition \ref{def: PBW}.

\item[(v)] $i_{U}:V\rightarrow U\left( V,c,b\right) $ induces an isomorphism
between $V$ and $P\left( U\left( V,c,b\right) \right) .$
\end{enumerate}
\end{theorem}

\begin{proof}
Set $T:=T\left( V,c\right) ,$ $S:=S\left( V,c\right) $ and $U:=U\left(
V,c,b\right) .$ Let $p_{t}:U_{t}^{\prime }\rightarrow U_{t}^{\prime
}/U_{t-1}^{\prime }$ be the canonical projection and let $i:=V\rightarrow
U_{1}^{\prime }$ be the corestriction of $i_{U}.$ Let $\theta
_{1}:V\rightarrow U_{1}^{\prime }/U_{0}^{\prime }:v\mapsto \pi _{U}\left(
v\right) +U_{0}^{\prime }$. Then $p_{1}\circ i=\theta _{1}.$ Let us prove
that 
\begin{equation}
\ker \left( \theta _{1}\right) =V\cap \ker \left( \pi _{U}\right) .
\label{form: ker theta1}
\end{equation}%
Let $x\in \ker \left( \theta _{1}\right) .$ Then $\pi _{U}\left( x\right)
\in U_{0}^{\prime }=\pi _{U}\left( K\right) \ $so that there exists $k\in K$
such that $\pi _{U}\left( x\right) =\pi _{U}\left( k\right) $. From this we
get $k=\varepsilon _{U}\pi _{U}\left( k\right) =\varepsilon _{U}\pi
_{U}\left( x\right) =\varepsilon _{T}\left( x\right) =0$ where the last
equality holds as $V\subseteq P\left( T\right) .$ Thus $\pi _{U}\left(
x\right) =0\ $and hence $x\in V\cap \ker \left( \pi _{U}\right) .$ The other
inclusion is trivial.

$\left( i\right) \Leftrightarrow \left( ii\right) $ It follows by Definition %
\ref{def: c-Lie alg}.

$\left( ii\right) \Leftrightarrow \left( iii\right) $ It follows trivially
by (\ref{form: ker theta1}).

$\left( iii\right) \Rightarrow \left( iv\right) $ Since $\mathrm{sdeg}\left(
V,c\right) \leq 1$, we have that $S\left( V,c\right) $ is strongly $%
\mathbb{N}
$-graded as a coalgebra so that (cf. Theorem \ref{teo: strongness}) $P\left(
S\right) =V$. Therefore, by \cite[Lemma 5.3.3]{Montgomery}, we obtain that $%
\theta $ is injective as it is injective its restriction $\theta _{1}$ to $%
P\left( S\left( V,c\right) \right) $. In conclusion the surjective map $%
\theta $ is indeed bijective.

$\left( iv\right) \Rightarrow \left( v\right) $ Let $z\in P\left( U\right) .$
Since $z\in U,$ there is $n\in 
\mathbb{N}
$ such that $z\in U_{n}^{\prime }\backslash U_{n-1}^{\prime }.$Let $%
\overline{z}:=z+U_{n-1}^{\prime }\in U_{n}^{\prime }/U_{n-1}^{\prime }$. By
definition of the bialgebra structure of $G:=\mathrm{gr}^{\prime }\left(
U\left( V,c,b\right) \right) $ one has that $\overline{z}\in P\left(
G\right) .$ By $\left( iv\right) $, we have $P\left( G\right) =\theta \left(
P\left( S\right) \right) =\theta \left( V\right) =U_{1}^{\prime
}/U_{0}^{\prime }.$ Thus $\overline{z}\in \left( U_{n}^{\prime
}/U_{n-1}^{\prime }\right) \cap \left( U_{1}^{\prime }/U_{0}^{\prime
}\right) .$ Since $\overline{z}\neq 0,$ we get $n=1.$ Hence $z\in
U_{n}^{\prime }=U_{1}^{\prime }=\pi _{U}\left( T_{\left( 1\right) }\right)
=\pi _{U}\left( K\oplus V\right) .$ Therefore, there are $k\in K$ and $v\in
V $ such that $z=\pi _{U}\left( k\right) +\pi _{U}\left( v\right) $ so that $%
k=\varepsilon _{U}\pi _{U}\left( k\right) =\varepsilon _{U}\left( z\right)
-\varepsilon _{U}\pi _{U}\left( v\right) =\varepsilon _{U}\left( z\right)
-\varepsilon _{T}\left( v\right) =0 $ where the last equality holds as $z\in
P\left( U\right) $ and $v\in P\left( T\right) .$ Then $z=\pi _{U}\left(
k\right) +\pi _{U}\left( v\right) =\pi _{U}\left( v\right) =i_{U}\left(
v\right) .$ We have so proved that $P\left( U\right) \subseteq \mathrm{Im}%
\left( i_{U}\right) .$ The other inclusion is trivial so that $\mathrm{Im}
\left( i_{U}\right) =P\left( U\right) .$

It remains to prove that $i_{U}$ is injective or, equivalently, that $i$ is
injective. Recall that, by Proposition \ref{pro: theta}, the map $\theta $
lifts the map $\theta _{1}.$

Since the canonical map $V\rightarrow S$ is injective (cf. \ref{pro: (V,c,0)
Lie}) and since, by hypothesis, $\theta $ is injective we infer that $\theta
_{1}$ and hence $i$ is injective too.

$\left( v\right) \Rightarrow \left( i\right) $ By hypothesis, $i_{U}$ is
injective.
\end{proof}

\begin{definition}
A connected braided bialgebra $\left( A,c_{A}\right) $ is called \textbf{%
primitively generated} if it is generated as a $K$-algebra by its space $%
P\left( A\right) $ of primitive elements.
\end{definition}

The following is one of the main results of this section.

\begin{theorem}
\label{teo: generated}Let $\left( A,c_{A}\right) $ be a primitively
generated connected braided bialgebra and let $\left( P,c_{P},b_{P}\right) $
be its infinitesimal braided Lie algebra. If $\mathrm{sdeg}\left(
P,c_{P}\right) \leq 1,$ then $A$ is isomorphic to $U\left(
P,c_{P},b_{P}\right) $ as a braided bialgebra.
\end{theorem}

\begin{proof}
Set $U:=U\left( P,c_{P},b_{P}\right) .$ By Theorem \ref{teo: univ U}, the
identity map of $P$ can be lifted to a morphism of braided bialgebras $%
\alpha :U\rightarrow A.$ Since $P$ generates $A$ as a $K$-algebra, this
morphism is surjective. On the other hand, by Theorem \ref{teo: magnum}, $%
i_{U}:P\rightarrow U$ induces an isomorphism between $V$ and $P\left(
U\right) .$ Hence the restriction of $\alpha $ on $P\left( U\right) $ is
injective. By \cite[Lemma 5.3.3]{Montgomery} $\alpha $ is injective.
\end{proof}

\section{Some braided vector space of strongness degree not greater the one 
\label{section: some braided}}

In this section meaningful examples of braided vector space of strongness
degree not greater the one are given. We also illustrate an example of a
braided vector space of strongness degree two.

\begin{theorem}
\label{teo: useful}Let $L\subseteq E\left( V,c\right) \ $and let $R=T\left(
V,c\right) /\left( L\right) .$ If $R$ inherits from $T\left( V,c\right) $ a
braided bialgebra structure such that the space $P\left( R\right) $ of
primitive elements of $R$ identifies with the image of $V$ in $R$, then $%
S\left( V,c\right) =R$ and $\mathrm{sdeg}\left( V,c\right) \leq 1.$
\end{theorem}

\begin{proof}
Denote by $\tau :R\rightarrow S:=S\left( V,c\right) ,$ $\pi _{S}:T\left(
V,c\right) \rightarrow S$ and $\pi _{R}:T\left( V,c\right) \rightarrow R$
the canonical projections. Since $\tau \pi _{R}=\pi _{S}$ and both $\pi _{R}$
and $\pi _{S}$ are coalgebra homomorphism, then $\tau $ is a coalgebra
homomorphism too. By \cite[Lemma 5.3.3]{Montgomery}, we obtain that $\tau $
is injective as its restriction to $P\left( R\right) $ is injective (note
that the map $i_{S}:V\rightarrow S\left( V,c\right) $ is injective by
Proposition \ref{pro: (V,c,0) Lie}). Hence $\tau $ is an isomorphism so that 
$S\left( V,c\right) =R.$ Therefore $P\left( S\right) $ identifies with $V\ $%
whence $S$ is strongly $%
\mathbb{N}
$-graded as a coalgebra which means $\mathrm{sdeg}\left( V,c\right) \leq 1.$
\end{proof}

\begin{definition}
\label{def: diagonal type}Recall that a braided vector space $\left(
V,c\right) $ is of \textbf{diagonal type} is there is a basis $x_{1},\ldots
,x_{n}$ of $V$ over $K$ and a matrix $\left( q_{i,j}\right) ,$ $q_{i,j}\in
K\backslash \left\{ 0\right\} ,$ such that $c\left( x_{i}\otimes
x_{j}\right) =q_{i,j}x_{j}\otimes x_{i},1\leq i,j\leq n.$
\end{definition}

\begin{theorem}
\label{teo: GK dim}Let $\left( V,c\right) $ be a braided vector space of
diagonal type such that $\emph{B}\left( V,c\right) $ is a domain and its
Gelfand-Kirillov dimension is finite. Then $\mathrm{sdeg}\left( V,c\right)
\leq 1.$
\end{theorem}

\begin{proof}
By hypothesis there is a basis $x_{1},\ldots ,x_{n}$ of $V$ over $K$ and a
matrix $\left( q_{i,j}\right) ,$ $q_{i,j}\in K\backslash \left\{ 0\right\} ,$
such that $c\left( x_{i}\otimes x_{j}\right) =q_{i,j}x_{j}\otimes
x_{i},1\leq i,j\leq n.$ By \cite[Corollary 3.11]{Andruskiewitsch-Angiono}
(which was obtained from results of Lusztig and Rosso), $R:=\emph{B}\left(
V,c\right) $ as a graded braided bialgebra divides out $T\left( V,c\right) $
by the ideal generated by set $L$ whose elements are $z_{i,j}:=\left(
ad_{c}x_{i}\right) ^{m_{i,j}+1}\left( x_{j}\right) ,$ $1\leq i,j\leq n,$
where $i\neq j$ and $m_{i,j}\ $is a suitable non-negative integer. Since $%
z_{i,j}$ is primitive in $T\left( V,c\right) $ (see e.g. \cite[Theorem 6.1]%
{Kharchenko-AnAlgSkew}, where $z_{i,j}$ is denoted by $W\left(
x_{j},x_{i}\right) $ in formulae (18) and (19), or see \cite[Lemma A.1]{AS-
FiniteQuantCartan}) and of degree at least two, we get $L\subseteq E\left(
V,c\right) .$ By Theorem \ref{teo: useful}, $\mathrm{sdeg}\left( V,c\right)
\leq 1.$
\end{proof}

\begin{example}
Assume $K$ is algebraically closed of characteristic $0$. Let $V$ be a $n$%
-dimensional vector space with basis $x_{1},\ldots ,x_{n}.$ Let $q_{i,j}\in
K,i,j\in \left\{ 1,\ldots ,n\right\} ,$ be such that $q_{i,i}\neq 1$ are $%
N_{i}$-roots of unity ($N_{i}>1$) and $q_{i,j}q_{j,i}=1$ for $i\neq j.$
Define $c:V\otimes V\rightarrow V\otimes V$ by setting $c\left( x_{i}\otimes
x_{j}\right) :=q_{i,j}x_{j}\otimes x_{i}$ (whence $c$ is of diagonal type).
One checks that $c$ fulfills (\ref{ec: braided equation}). By means of the
quantum binomial formula (see \cite[Lemma 3.6]{AS}) one gets that $%
x_{i}\otimes x_{j}-q_{i,j}x_{j}\otimes x_{i}\in E_{2}\left( V,c\right) $ for 
$i\neq j$ and $x_{i}^{\otimes N_{i}}\in E_{N_{i}}\left( V,c\right) .$
Consider the \textbf{quantum linear space} $R=K<x_{1},\ldots ,x_{n}\mid
x_{1}^{N_{1}}=0;\ldots ;x_{n}^{N_{n}}=0;x_{i}x_{j}=q_{i,j}x_{j}x_{i}>.$ By
Theorem \ref{teo: useful}, $S\left( V,c\right) =R$ and $\mathrm{sdeg}\left(
V,c\right) \leq 1.$ Since $S\left( V,c\right) \neq T\left( V,c\right) $ it
is clear that $\mathrm{sdeg}\left( V,c\right) =1.$
\end{example}

\begin{remark}
\label{rem: mutlilinear}Let $\left( V,c\right) $ be a braided vector space
of diagonal type as in Definition \ref{def: diagonal type}. By \cite[Theorem
1]{Kharchenko-AnExistence}, the space $E_{m}\left( V,c\right) $ contains an
element of the form $\sum \limits_{\sigma \in S_{m}}\lambda _{\sigma
}x_{\sigma \left( t_{1}\right) }\otimes \cdots \otimes x_{\sigma \left(
t_{m}\right) }$ with $0\leq t_{1}<\ldots <t_{m}\leq n,$ $\lambda _{\sigma
}\in K$, if and only if 
\begin{equation*}
\prod_{1\leq i\neq j\leq n}q_{t_{i},t_{j}}=1.
\end{equation*}%
Furthermore, if this condition holds, then $\left( n-2\right) !\leq \dim
_{K}\left( E_{n}\left( V,c\right) \right) .$

For instance, in the previous example, condition $q_{i,j}q_{j,i}=1$ is
necessary to have that $x_{i}\otimes x_{j}-q_{i,j}x_{j}\otimes x_{i}\in
E_{2}\left( V,c\right) $ for $i\neq j.$
\end{remark}

\begin{definition}
\cite[Definition 1.6]{AS} Let $\left( V,c\right) $ be a braided vector
space. Recall that $\left( V,c\right) $ is of \textbf{group type} if $V$ has
a basis $x_{1},\ldots ,x_{n}$ such that $c\left( x_{i}\otimes x_{j}\right)
=\gamma _{i}\left( x_{j}\right) \otimes x_{i}$ for some $\gamma _{1},\ldots
,\gamma _{n}\in GL\left( V\right) .$ If in addition the subgroup $G\left(
V,c\right) $ of $GL\left( V\right) $ generated by $\gamma _{1},\ldots
,\gamma _{n}$ is abelian we will say that $\left( V,c\right) $ is of \textbf{%
abelian group type}.
\end{definition}

\begin{proposition}
\label{pro: abelian grp} Let $\left( V,c\right) $ be a two dimensional
braided vector space of abelian group type. Suppose that $\dim _{K}S\left(
V,c\right) <32.$ Then $\mathrm{sdeg}\left( V,c\right) \leq 1.$
\end{proposition}

\begin{proof}
We have two graded braided bialgebra homomorphisms $\pi _{S}:T\left(
V,c\right) \twoheadrightarrow S\left( V,c\right) $ and $\pi _{1}^{\infty
}:S\left( V,c\right) \twoheadrightarrow \emph{B}\left( V,c\right) .$ By a
famous Takeuchi's result (see \cite[Lemma 5.2.10]{Montgomery}), these
bialgebras are indeed Hopf algebras whence $\pi _{S}$ and $\pi _{1}^{\infty
} $ are in fact graded braided Hopf algebra homomorphisms. It is well known
that $V$ can be regarded as an object in the monoidal category $\mathcal{M}$
of Yetter-Drinfeld modules over $G\left( V,c\right) $. Furthermore $c$
coincides with the braiding of $V$ in the category. Clearly $T\left(
V,c\right) \ $and $\emph{B}\left( V,c\right) $ are Hopf algebras in $%
\mathcal{M}$ and the same is true for $S\left( V,c\right) $ as $E\left(
V,c\right) $ is an object in $\mathcal{M}$. Furthermore $\pi _{S}$ and $\pi
_{1}^{\infty }$ are morphism in $\mathcal{M}$. Since $\emph{B}\left(
V,c\right) $ is a quotient of $S\left( V,c\right) ,$ one has $\dim _{K}\emph{%
B}\left( V,c\right) \leq \dim _{K}S\left( V,c\right) <32.$ By \cite[Lemma 6.2%
]{Grana- low dim}, we have that $\pi _{1}^{\infty }$ is an isomorphism i.e. $%
\mathrm{sdeg}\left( V,c\right) \leq 1.$
\end{proof}

\begin{lemma}
\label{lem: S is graded}Let $(V,c)$ be a braided vector space. Then $%
I=\left( E\left( V,c\right) \right) $ is a graded ideal of $T\left(
V,c\right) $ with graded component $I_{n}:=\left( E\left( V,c\right) \right)
\cap V^{\otimes n}.$ Moreover 
\begin{equation}
I_{n}=I_{n-1}\otimes V+V\otimes I_{n-1}+E_{n}\left( V,c\right) \text{, for
every }n\in 
\mathbb{N}
\text{.}  \label{form: graded S}
\end{equation}
\end{lemma}

\begin{proof}
For the first part see the proof of Theorem \ref{teo: S(B) graded}. The last
part follows by induction on $n\in 
\mathbb{N}
.$ 
\end{proof}

\begin{theorem}
\label{teo: minimal poly deg 1}Let $V\neq \left\{ 0\right\} $ be a $K$%
-vector space, let $q\in K\backslash \left\{ 0\right\} $ and set $c=q\mathrm{%
Id}_{V\otimes V}.$

Then $\left( V,c\right) $ is a braided vector space and

\begin{enumerate}
\item[i)] $\mathrm{sdeg}\left( V,c\right) =0$, if $q$ is regular.

\item[ii)] $\mathrm{sdeg}\left( V,c\right) =1$, otherwise. In this case $q$
is a primitive $N$-th root of unity for some $N\geq 2$ and $S\left(
V,c\right) =T\left( V,c\right) /\left( V^{N}\right) $.
\end{enumerate}
\end{theorem}

\begin{proof}
Clearly $c:V\otimes V\rightarrow V\otimes V$ is a Yang-Baxter operator over $%
V$.

i) It follows by Theorem \ref{teo: sdeg(B)=0}.

ii) Set $T:=T\left( V,c\right) $ and $S:=S\left( V,c\right) $. By (\ref%
{form: grado 1 reduced}), which holds for a not necessarily regular $q\in K,$
we have $\mathrm{\Gamma }_{n,1}^{T}=\left( n+1\right) _{q}\mathrm{Id}%
_{B^{n+1}},$ for every $n\geq 1.$ Using this equality and $\mathrm{\Gamma }%
_{n}^{T}=\left( \mathrm{\Gamma }_{n-1}^{T}\otimes V\right) \mathrm{\Gamma }%
_{n-1,1}^{T}$, inductively one proves that 
\begin{equation}
\mathrm{\Gamma }_{n}^{T}=\left( n\right) _{q}!\mathrm{Id}_{V^{\otimes n}},%
\text{ for every }n\geq 2\text{.}  \label{form: grado n}
\end{equation}

%
Now, since $q$ is not regular, there is a minimal $N\geq 2$ such that $%
\mathrm{\Gamma }_{N}^{T}=0$. This entails that $\left( N\right) _{q}!=0$ but 
$\left( N-1\right) _{q}!\neq 0$ i.e. $q$ is a primitive $N$-th root of unity.

Let $I=\left( E\left( V,c\right) \right) $ and let $I_{n}:=\left( E\left(
V,c\right) \right) \cap V^{\otimes n}.$ Assume that 
\begin{equation}
\ker \left( \mathrm{\Gamma }_{n}^{T}\right) \subseteq I_{n}\text{ for every }%
n\in 
\mathbb{N}
\text{.}  \label{form: Gamma I}
\end{equation}%
Since $\mathrm{\Gamma }^{T}:T\left( V,c\right) \rightarrow T^{c}\left(
V,c\right) $ is a graded algebra homomorphism with $n$-th graded component $%
\mathrm{\Gamma }_{n}^{T}$ (see \ref{claim: Gamma}), we get $\mathrm{\Gamma }%
_{n}^{T}\left( I_{n}\right) =\mathrm{\Gamma }^{T}\left[ \left( E\left(
V,c\right) \right) \cap V^{\otimes n}\right] \subseteq \mathrm{\Gamma }^{T}%
\left[ \left( E\left( V,c\right) \right) \right] \subseteq \left( \mathrm{%
\Gamma }^{T}\left[ E\left( V,c\right) \right] \right) =0$ so that $%
I_{n}\subseteq \ker \left( \mathrm{\Gamma }_{n}^{T}\right) $ whence $%
I_{n}=\ker \left( \mathrm{\Gamma }_{n}^{T}\right) .$ In this case 
\begin{equation*}
S=\oplus _{n\in 
\mathbb{N}
}\frac{V^{\otimes n}}{I_{n}}=\oplus _{n\in 
\mathbb{N}
}\frac{V^{\otimes n}}{\ker \left( \mathrm{\Gamma }_{n}^{T}\right) }=\emph{B}%
\left( V,c\right) .
\end{equation*}%
Moreover, in view of (\ref{form: grado n}), we deduce 
\begin{equation*}
\ker \left( \mathrm{\Gamma }_{n}^{T}\right) =\left\{ 
\begin{tabular}{ll}
$0,$ & for every $0\leq n\leq N-1,$ \\ 
$V^{n},$ & for every $N\leq n$%
\end{tabular}%
\right.
\end{equation*}%
whence $\emph{B}\left( V,c\right) =T\left( V,c\right) /\left( V^{N}\right) $%
. Since $\emph{B}\left( V,c\right) $ is strongly $\mathbb{N}$-graded as a coalgebra, we have that $\mathrm{sdeg}\left( V,c\right) =1$. It
remains to prove (\ref{form: Gamma I}). Hence it remains to check that $%
I_{n}=V^{n},$ for every $N\leq n.$ In view of (\ref{form: graded S}) it
suffices to prove that $E_{N}\left( V,c\right) =V^{N}$. For every $a,b\geq
1,a+b=N,$ we have $0=\mathrm{\Gamma }_{N}^{T}=\mathrm{\Gamma }_{a+b}^{T}%
\overset{\text{(\ref{form: SS})}}{=}\left( \mathrm{\Gamma }_{a}^{T}\otimes 
\mathrm{\Gamma }_{b}^{T}\right) \mathrm{\Gamma }_{a,b}^{T}.$ Since $a,b\leq
N-1,$ then $\mathrm{\Gamma }_{a}^{T}$ and $\mathrm{\Gamma }_{b}^{T}$ are
injective so that $\mathrm{\Gamma }_{a,b}^{T}=0.$ Since $m_{T}^{a,b}$ is an
isomorphism, we get $E_{N}\left( V,c\right) =V^{N}$.
\end{proof}

Next we give some examples of braided vector spaces of strongness degree at
least too. The most interesting of them are due to the work of V. K.
Kharchenko on combinatorial rank (see Remark \ref{rem: combinatorial}).

\begin{example}
\label{es: AG}From \cite[Example 3.3.22]{AG}, we quote the following
example. Take $H=K\mathbb{D}_{4},$ where $\mathbb{D}_{4}=\left\{ \sigma
,\rho \mid o\left( \sigma \right) =2,o\left( \rho \right) =4,\rho \sigma
=\sigma \rho ^{-1}\right\} .$ Consider the module $V$ in ${_{H}^{H}\mathcal{%
YD}}$ with basis $\left\{ z_{0},z_{1},z_{2},z_{3}\right\} $, with the
structure given by%
\begin{equation*}
\delta \left( z_{i}\right) =\rho ^{i}\sigma \oplus z_{i},\qquad \rho
^{i}\vartriangleright z_{i}=z_{i+2j},\qquad \sigma \vartriangleright
z_{i}=-z_{-i}
\end{equation*}%
where we take subindexes of the $z_{i}$ to be on $\mathbb{Z}/4\mathbb{Z}.$
The braiding is then given by $c\left( z_{i}\otimes z_{j}\right)
=-z_{2i-j}\otimes z_{i}.$ Let $T:=T\left( V,c\right) $, let $%
a:=z_{1}z_{2}+z_{0}z_{1}$, $b:=z_{1}z_{0}+z_{2}z_{1}$ and let $Z$ be the set
consisting of the following elements:%
\begin{gather}
z_{0}^{2},\quad z_{1}^{2},\quad z_{2}^{2},\quad z_{3}^{2},\quad
z_{0}z_{2}+z_{2}z_{0},\quad z_{1}z_{3}+z_{3}z_{1},  \label{form: AG1} \\
z_{0}z_{1}+z_{1}z_{2}+z_{2}z_{3}+z_{3}z_{0},\quad
z_{0}z_{3}+z_{1}z_{0}+z_{2}z_{1}+z_{3}z_{2}  \label{form: AG2} \\
a^{2},\quad b^{2},\quad ab+ba.  \label{form: AG3}
\end{gather}%
Then the Nichols algebra $\mathcal{B}\left( V,c\right) $ comes out to be $%
T/\left( Z\right) .$ One can check that:

\begin{itemize}
\item $E_{2}\left( V,c\right) \ $is generated over $K$ by the elements in (%
\ref{form: AG1}) and (\ref{form: AG2});

\item $E_{3}\left( V,c\right) $ is generated over $K$ by the elements $%
z_{1}^{2}z_{3}-z_{3}z_{1}^{2},$ $z_{1}z_{3}^{2}-z_{3}^{2}z_{1},$ $%
z_{2}^{2}z_{0}-z_{0}z_{2}^{2},$ $z_{2}z_{0}^{2}-z_{0}^{2}z_{2},$ $%
z_{1}\left( z_{0}z_{2}+z_{2}z_{0}\right) -\left(
z_{0}z_{2}+z_{2}z_{0}\right) z_{1},$ $z_{3}\left(
z_{0}z_{2}+z_{2}z_{0}\right) -\left( z_{0}z_{2}+z_{2}z_{0}\right) z_{3},$ $%
z_{2}\left( z_{1}z_{3}+z_{3}z_{1}\right) -\left(
z_{1}z_{3}+z_{3}z_{1}\right) z_{2}\in E_{2}\left( V,c\right) \otimes
V+V\otimes E_{2}\left( V,c\right) ;$

\item $E_{4}\left( V,c\right) =\left\{ 0\right\} .$
\end{itemize}

Thus, by (\ref{form: graded S}) we obtain $\left( E\left( V,c\right) \right)
\cap V^{\otimes 4}\subseteq E_{2}\left( V,c\right) \otimes V\otimes
V+V\otimes E_{2}\left( V,c\right) \otimes V+V\otimes V\otimes E_{2}\left(
V,c\right) .$ Since $%
a^{2}=z_{1}z_{2}z_{1}z_{2}+z_{1}z_{2}z_{0}z_{1}+z_{0}z_{1}^{2}z_{2}+z_{0}z_{1}z_{0}z_{1} 
$ is not an element in the ideal generated by $E_{2}\left( V,c\right) $,
then $\mathrm{sdeg}\left( V,c\right) =\mathrm{sdeg}\left( T\right) >1.$ In
fact the elements in (\ref{form: AG3}) belongs to $E_{4}\left( S\left(
V,c\right) \right) $ so that we obtain $\mathrm{sdeg}\left( V,c\right) =2.$
\end{example}

\begin{example}
\label{ex: twodim sdeg 2}At the end of \cite{Kharchenko-SkewPrim}, an
example of a two dimensional braided vector space $\left( V,c\right) $ of
strongness degree $2$ is given. The braiding $c$ is of diagonal type of the
form%
\begin{equation*}
c\left( x_{i}\otimes x_{j}\right) =q_{i,j}x_{j}\otimes x_{i},1\leq i,j\leq 2.
\end{equation*}%
where $x_{1},x_{2}$ is a basis of $V$ over $K$, $q_{1,2}=1\neq -1$ and $%
q_{i,j}=-1$ for all $\left( i,j\right) \neq \left( 1,2\right) .$
\end{example}

\begin{example}
\label{ex: Cartan}In \cite{Kharchenko-Andrade}, one can find the strongness
degree of a braided vector space with diagonal braiding of Cartan type $%
A_{n}.$ That is, if $c\left( x_{i}\otimes x_{j}\right) =q_{i,j}x_{j}\otimes
x_{i}$ with 
\begin{equation*}
q_{i,i}=q_{jj}=q,\qquad q_{i,i+1}q_{i+1,i}=q^{-1},\qquad
q_{i,j}q_{j,i}=1,i-j>1
\end{equation*}%
where $q$ is a primitive $t$-th root of unity, $t>2,$ then the strongness
degree of the related braided vector space equals $\lceil 1+log_{2}(n)\rceil
.$ If $q$ is not a root of unity, then the strongness degree is one.
\end{example}

\section{Braidings of Hecke type \label{section: Hecke}}

\begin{definition}
\label{def: Hecke type} (see e.g. \cite[Definition 3.3]{AS} or \cite[%
Definition 3.1.1]{Abella-Andruskiewitsch}) Let $\left( V,c\right) $ be a
braided vector space. We say that $c$ is a braiding of \textbf{Hecke type
(or a Hecke symmetry)} \textbf{with mark }$q$ if $c$ satisfies the equation $%
\left( c+\mathrm{Id}_{V\otimes V}\right) \left( c-q\mathrm{Id}_{V\otimes
V}\right) =0$ for some $q\in K\backslash \left\{ 0\right\} .$
\end{definition}

\begin{theorem}
\label{teo: sdeg(B)=1}Let $\left( B,m_{B},u_{B},\Delta _{B},\varepsilon
_{B},c_{B}\right) $ be a $0$-connected graded braided bialgebra. Assume that

\begin{itemize}
\item $B$ is strongly $%
\mathbb{N}
$-graded as an algebra and $m_{B}^{1,1}$ is an isomorphism.

\item $c_{B}^{1,1}$ has minimal polynomial $f\left( X\right) =\left(
X+1\right) \left( X-q\right) ,$ for some regular $q\in K$.
\end{itemize}

Then $\mathrm{sdeg}\left( B\right) =1$.
\end{theorem}

\begin{proof}
By Theorem \ref{teo: sdeg(B)=n-1}, we have $\mathrm{sdeg}\left( B\right) =%
\mathrm{sdeg}\left( S\left( B\right) \right) +1$ and $m_{S\left( B\right)
}^{1,1}\left( c_{S\left( B\right) }^{1,1}-q\right) =0.$ Therefore, by
Theorem \ref{teo: sdeg(B)=0}, we get that $\mathrm{sdeg}\left( S\left(
B\right) \right) =0$ whence $\mathrm{sdeg}\left( B\right) =1.$
\end{proof}

The following result explains how braided vector spaces with braiding of
Hecke type fits into our study of braided vector spaces of strongness degree
at most $1.$

\begin{theorem}
\label{teo: sdeg Hecke}Let $\left( V,c\right) $ be braided vector space such
that $c$ is a braiding of Hecke type with regular mark $q.$ Then $\mathrm{%
sdeg}\left( V,c\right) \leq 1$.
\end{theorem}

\begin{proof}
By hypothesis, $c$ satisfies the equation $\left( c+\mathrm{Id}_{V\otimes
V}\right) \left( c-q\mathrm{Id}_{V\otimes V}\right) =0$ so that $c$ has
minimal polynomial a divisor $f$ of the polynomial $\left( X+1\right) \left(
X-q\right) .$

Now, the tensor algebra $T:=T\left( V,c\right) $ is strongly $%
\mathbb{N}
$-graded as an algebra and $m_{T}^{1,1}$ is an isomorphism. Moreover $%
c_{T}^{1,1}=c$. Thus, if $f=\left( X+1\right) \left( X-q\right) $, by
Theorem \ref{teo: sdeg(B)=1}, we have $\mathrm{sdeg}\left( V,c\right) :=%
\mathrm{sdeg}\left( T\left( V,c\right) \right) =1.$ On the other hand, by
Theorem \ref{teo: minimal poly deg 1}, if $f=X+1$, we get $\mathrm{sdeg}%
\left( V,c\right) =1\ $while, if $f=X-q$, we obtain $\mathrm{sdeg}\left(
V,c\right) =0.$
\end{proof}

Using the results in \cite{AMS-MM2} and \cite{AMS-MM}, we will now give a simple description of
the universal enveloping algebra of a braided Lie algebra such that the
braiding is of Hecke type with regular mark. We point out that this
description does not depend on the bracket in the non-symmetric case
whenever $\mathrm{char}\left( K\right) \neq 2.$

\begin{theorem}
\label{teo: Hecke Lie}Let $\left( V,c,b\right) $ be a braided Lie algebra
such that $c$ is of Hecke type with regular mark $q$. Then%
\begin{equation*}
U\left( V,c,b\right) =\frac{T\left( V,c\right) }{\left( c\left( z\right) -qz-%
\left[ z\right] _{b}\mid z\in V\otimes V\right) }
\end{equation*}%
where $\left[ z\right] _{b}=b\left( c\left( z\right) -qz\right) $, for every 
$z\in V\otimes V$. Moreover if $q\neq 1$ and $\mathrm{char}\left( K\right)
\neq 2,$ then $\left[ -\right] _{b}=0.$
\end{theorem}

\begin{proof}
We set $T:=T\left( V,c\right) $. Let us prove that the following equality 
\begin{equation}
\mathrm{Im}\left( c-q\mathrm{Id}_{V\otimes V}\right) =E_{2}\left( V,c\right)
\label{form: remains 2}
\end{equation}%
holds. Since $\left( c+\mathrm{Id}_{V\otimes V}\right) \left( c-q\mathrm{Id}%
_{V\otimes V}\right) =0,$ it is clear that $\mathrm{Im}\left( c-q\mathrm{Id}%
_{V\otimes V}\right) \subseteq \ker \left( c+\mathrm{Id}_{V\otimes V}\right)
.$ Let $z\in \ker \left( c+\mathrm{Id}_{V\otimes V}\right) .$ Then $\left(
c-q\mathrm{Id}_{V\otimes V}\right) \left( z\right) =c\left( z\right)
-qz=-\left( 1+q\right) z=-\left( 2\right) _{q}z.$ Since $q$ is regular, we
have $\left( 2\right) _{q}\neq 0$, so that $z=-\left( 2\right)
_{q}^{-1}\left( c-q\mathrm{Id}_{V\otimes V}\right) \left( z\right) \in \text{%
Im}\left( c-q\mathrm{Id}_{V\otimes V}\right) .$ Hence $\mathrm{Im}\left( c-q%
\mathrm{Id}_{V\otimes V}\right) =\ker \left( c+\mathrm{Id}_{V\otimes
V}\right) =\ker \left( \Delta _{T}^{1,1}\right) =E_{2}\left( V,c\right) .$
We have so proved (\ref{form: remains 2}). By (\ref{form: remains 2}), we
have that $\left( c-q\mathrm{Id}_{V\otimes V}\right) \left( z\right) $ lies
in the domain of $b$ for every $z\in V\otimes V.$ Hence we can define a map $%
\left[ -\right] _{b}:V\otimes V\rightarrow V$ by setting $\left[ z\right]
_{b}=b\left( c-q\mathrm{Id}_{V\otimes V}\right) \left( z\right) $, for every 
$z\in V\otimes V$.

Let us check that $\left[ -\right] _{b}$ is a $c$-bracket in the sense of 
\cite[Definition 4.1]{AMS-MM} i.e. that $c\left( \left[ -\right] _{b}\otimes
V\right) =\left( V\otimes \left[ -\right] _{b}\right) c_{1}c_{2}$ and $%
c\left( V\otimes \left[ -\right] _{b}\right) =\left( \left[ -\right]
_{b}\otimes V\right) c_{2}c_{1}.$ By (\ref{form: c-bracket}), $b$ is a
bracket on $\left( V,c\right) $ if $c\left( b\otimes V\right) =\left(
V\otimes b\right) c_{E\left( V,c\right) ,V} $ and $c\left( V\otimes b\right)
=\left( b\otimes V\right) c_{V,E\left( V,c\right) }.$ By (\ref{form: remains
2}), the previous formulas rereads as follows%
\begin{eqnarray*}
c\left( b\otimes V\right) \left[ \left( c-q\mathrm{Id}_{V\otimes V}\right)
\otimes V\right] \left( z\right) &=&\left( V\otimes b\right) c_{1}c_{2}\left[
\left( c-q\mathrm{Id}_{V\otimes V}\right) \otimes V\right] \left( z\right)
\qquad \text{and} \\
c\left( V\otimes b\right) \left[ V\otimes \left( c-q\mathrm{Id}_{V\otimes
V}\right) \right] \left( z\right) &=&\left( b\otimes V\right) c_{2}c_{1}%
\left[ V\otimes \left( c-q\mathrm{Id}_{V\otimes V}\right) \right] \left(
z\right) ,
\end{eqnarray*}%
for every $z\in V^{\otimes 3}.$ Thus $\left[ -\right] _{b}$ is a $c$-bracket.

Set $U_{H}\left( V,c,\left[ -\right] _{b}\right) :=T\left( V,c\right)
/\left( c\left( z\right) -qz-\left[ z\right] _{b}\mid z\in V\otimes V\right)
$. This is the enveloping algebra as defined in \cite[Definition 2.1]{AMS-MM2}. Let us prove that $U_{H}\left( V,c,\left[ -\right] _{b}\right) =U\left(
V,c,b\right) .$ 

By definition of $\left[ -\right] _{b}$ and (\ref{form: remains 2}), we have%
\begin{equation*}
\left( c\left( z\right) -qz-\left[ z\right] _{b}\mid z\in V\otimes V\right)
=\left( \left( \mathrm{Id}_{V\otimes V}-b\right) \left[ E_{2}\left(
V,c\right) \right] \right) .
\end{equation*}%
Therefore there is a canonical projection $\gamma :U_{H}\left( V,c,\left[ -%
\right] _{b}\right) \rightarrow U\left( V,c,b\right) .$ Since $%
i_{U}:V\rightarrow U\left( V,c,b\right) $ is injective and $\gamma \circ
i_{U_{H}}=i_{U}$, it is clear that the canonical map $i_{U_{H}}:V\rightarrow
U_{H}\left( V,c,\left[ -\right] _{b}\right) $ is injective too. By \cite[%
Theorem 4.20 and Corollary 4.23]{AMS-MM}, the injectivity of $i_{U_{H}}$
imply that $U_{H}\left( V,c,\left[ -\right] _{b}\right) $ is strictly graded
as a coalgebra. Therefore $P(U_{H}\left( V,c,\left[ -\right] _{b}\right) )$
identifies with $V\ $so that the restriction of the bialgebra homomorphism $%
\gamma $ to the space of primitive elements of $U_{H}\left( V,c,\left[ -%
\right] _{b}\right) $ is injective. By \cite[Lemma 5.3.3]{Montgomery} this
entails that $\gamma $ itself is injective and hence bijective. We have so
proved that $U_{H}\left( V,c,\left[ -\right] _{b}\right) =U\left(
V,c,b\right) .$

Assume now $q\neq 1$ and $\mathrm{char}\left( K\right) \neq 2.$ Then, by 
\cite[Theorem 4.3]{AMS-MM2}, since $q$ is regular, the map $\left[ -\right]
_{b}$ is necessarily zero.
\end{proof}

\begin{remark}
\label{rem: Gu}Note that, when $q=1,$ then regularity of $q$ means $\mathrm{%
char}\left( K\right) =0.$ Moreover $\left( V,c,\left[ -\right] _{b}\right) $
becomes a Lie algebra in the sense of \cite{Gu- Gen Trans Lie} and $%
U_{H}\left( V,c,\left[ -\right] _{b}\right) $ is the corresponding
enveloping algebra. In this case $\left[ -\right] _{b}$ is not necessarily
trivial in general.
\end{remark}

\section{Quadratic algebras \label{section: Quadratic}}

\begin{definition}
\label{def: quadratic}Recall that a \textbf{quadratic algebra} \cite[page 19]%
{Manin} is an associative graded $K$-algebra $A=\oplus _{n\in 
\mathbb{N}
}A^{n}$ such that:

\begin{enumerate}
\item[1)] $A^{0}=K$;

\item[2)] $A$ is generated as a $K$-algebra by $A^{1}$;

\item[3)] the ideal of relations among elements of $A^{1}$ is generated by
the subspace of all quadratic relations $R(A)\subseteq A^{1}\otimes A^{1}.$
\end{enumerate}

Equivalently $A$ is a graded $K$-algebra such that the natural map $\pi
:T_{K}(A^{1})\rightarrow A$ from the tensor algebra generated by $A^{1}$ is
surjective and $\mathrm{\ker }(\pi )$ is generated as a two sided ideal in $%
T_{K}(A^{1})$ by $\mathrm{\ker }(\pi )\cap (A^{1}\otimes A^{1}).$
\end{definition}

\begin{corollary}
\label{coro: quadratic}Let $\left( B,m_{B},u_{B},\Delta _{B},\varepsilon
_{B},c_{B}\right) $ be a graded braided bialgebra which is strongly $%
\mathbb{N}
$-graded as an algebra. Then $\mathrm{sdeg}\left( B\right) \leq 1$ whenever $%
B^{0}\left[ B^{1}\right] $ is a quadratic algebra.
\end{corollary}

\begin{proof}
By Definition \ref{def: quadratic}, $B^{0}=K$ and the ideal of relations
among elements of $B^{1}$ is generated by the subspace of all quadratic
relations $R(B)\subseteq B^{1}\otimes B^{1}.$ By Corollary \ref{coro: f.g.
Type1}, we get $\mathrm{sdeg}\left( B\right) \leq 1.$
\end{proof}

\begin{theorem}
\label{teo: quadratic => sdeg 1}Let $\left( V,c\right) $ be braided vector
space. Assume that $\emph{B}\left( V,c\right) $ is a quadratic algebra. Then 
$\mathrm{sdeg}\left( V,c\right) \leq 1$.
\end{theorem}

\begin{proof}
Apply Corollary \ref{coro: quadratic} to the case $B:=T\left( V,c\right) $
once observed that $B^{0}\left[ B^{1}\right] =\emph{B}\left( V,c\right) .$
\end{proof}

\begin{remark}
\label{rem: B(V) quadratic}Examples of braided vector spaces $\left(
V,c\right) $ such that $\emph{B}\left( V,c\right) $ is a quadratic algebra
can be found e.g. in \cite{MS-Pointed indec} and in \cite{AG- Racks}. By 
\cite[Proposition 3.4]{AS}, another example is given by braided vector
spaces of Hecke-type with regular mark. This gives a different proof of
Theorem \ref{teo: sdeg Hecke}.
\end{remark}

\begin{theorem}
\label{teo: quadratic}Let $\left( A,c_{A}\right) $ be a connected braided
bialgebra such that the graded coalgebra associated to the coradical
filtration is a quadratic algebra with respect to its natural braided
bialgebra structure. Let $\left( P,c_{P},b_{P}\right) $ be the infinitesimal
braided Lie algebra of $A$. Then $\emph{B}\left( P,c_{P}\right) $ is a
quadratic algebra and $A$ is isomorphic to $U\left( P,c_{P},b_{P}\right) $
as a braided bialgebra.
\end{theorem}

\begin{proof}
Since $A$ is connected, then graded coalgebra $B:=\mathrm{gr}\left( A\right) 
$ associated to the coradical filtration $\left( A_{n}\right) _{n\in 
\mathbb{N}
}$ of $A$ is indeed a braided bialgebra. Since, by assumption, it is a
quadratic algebra (see Definition \ref{def: quadratic}), then $B$ is
generated as a $K$-algebra by $B^{1}=A_{1}/A_{0}$ so that $B$ is strongly $%
\mathbb{N}
$-graded as an algebra. Since, by construction, it is also strongly $%
\mathbb{N}
$-graded as a coalgebra, we get that $B\simeq B_{0}\left[ B_{1}\right]
\simeq \emph{B}\left( P,c_{P}\right) .$ Hence $\emph{B}\left( P,c_{P}\right) 
$ is a quadratic algebra. By Theorem \ref{teo: quadratic => sdeg 1}, $%
\mathrm{sdeg}\left( P,c_{P}\right) \leq 1.$ On the other hand $B$ is
generated as a $K$-algebra by $B^{1}=A_{1}/A_{0}$ implies that $P$ generates 
$A$ as a $K$-algebra. Therefore, by Theorem \ref{teo: generated}, $A$ is
isomorphic to $U\left( P,c_{P},b_{P}\right) $ as a braided bialgebra.
\end{proof}

\begin{example}
\label{ex: gu} We are going to apply the previous result to a braided
algebra $A$ that appeared in \cite[page 325]{Gurevich: Hecke sym}. Let us
recall its definition. Let $K$ be a field of characteristic $0$ and let $q,\alpha ,\beta \in K^{\ast }$ be such that $%
q\neq 1$ is regular and $\left( \alpha /\beta \right) ^{2}=q$. Set $m=-\alpha /\beta \ $%
and consider the vector space $V$ with basis $\left\{
e_{0},e_{1},e_{2}\right\} $ and braiding given by $c\left( e_{i}\otimes
e_{j}\right) =e_{j}\otimes e_{i}$ for $i=0$ or $j=0$, $c\left( e_{i}\otimes
e_{i}\right) =qe_{i}\otimes e_{i}$ for $i=1,2$, $c\left( e_{2}\otimes
e_{1}\right) =me_{1}\otimes e_{2}+\left( q-1\right) e_{2}\otimes e_{1}$ and $%
c\left( e_{1}\otimes e_{2}\right) =qm^{-1}e_{2}\otimes e_{1}.$ One can check
that $E_{2}\left( V,c\right) \ $is generated over $K$ by the elements $%
e_{1}\otimes e_{0}-e_{0}\otimes e_{1},e_{2}\otimes e_{0}-e_{0}\otimes e_{2}$
and $e_{2}\otimes e_{1}-me_{1}\otimes e_{2}.$ Consider the $K$-linear map $%
b_{2}:E_{2}\left( V,c\right) \rightarrow V$ defined, on generators, by 
\begin{equation*}
b_{2}\left( e_{1}\otimes e_{0}-e_{0}\otimes e_{1}\right) =e_{1},\qquad
b_{2}\left( e_{2}\otimes e_{0}-e_{0}\otimes e_{2}\right) =e_{2},\qquad
b_{2}\left( e_{2}\otimes e_{1}-me_{1}\otimes e_{2}\right) =0.
\end{equation*}%
Then, by definition, $A$ is given by%
\begin{equation*}
A:=\frac{T\left( V,c\right) }{\left( \left( \mathrm{Id}-b_{2}\right) \left[
E_{2}\left( V,c\right) \right] \right) }.
\end{equation*}%
Explicitly $A$ is generated over $K$ by $e_{0},e_{1},e_{2}$ with relations%
\begin{equation*}
e_{1}\otimes e_{0}-e_{0}\otimes e_{1}=e_{1},\qquad e_{2}\otimes
e_{0}-e_{0}\otimes e_{2}=e_{2}\qquad e_{2}\otimes e_{1}-m e_{1}\otimes e_{2}=0.
\end{equation*}
Now, $A$ is a braided bialgebra quotient of $T\left( V,c\right) $ (to prove
this one can adapt the proof of Theorem \ref{teo: U bialgebra} by taking $%
E_{2}\left( V,c\right) $ instead of $E\left( V,c\right) $ and using the
equalities $c\left( b_{2}\otimes V\right) =\left( V\otimes b_{2}\right)
c_{E_{2}\left( V,c\right) ,V}$ and $c\left( V\otimes b_{2}\right) =\left(
b_{2}\otimes V\right) c_{V,E_{2}\left( V,c\right) }$). Thus, by \cite[%
Corollary 5.3.5]{Montgomery}, $A$ is connected as $T\left( V,c\right) $ is
connected. Let $\left( P,c_{P},b_{P}\right) $ is the infinitesimal braided Lie algebra of $%
A. $ Then $A$ has basis $\{e_0^{n_0}e_1^{n_1}e_2^{n_2}| n_0,n_1,n_2\in \mathbb{N}\}$ and, using regularity of $q$, one gets (by the same argument as in the
proof of \cite[Lemma 3.3]{AS- Lifting}) that $(P,c_P)$ identifies with $(V,c)$. Thus, the graded coalgebra associated to the coradical
filtration of $A$ is $\emph{B}\left( V,c\right)$, which is isomorphic to $T\left( V,c\right) / \left(E_{2}\left( V,c\right)\right) $ which is quadratic. By Theorem \ref{teo: quadratic}, $A$ is
isomorphic to $U\left( P,c_{P},b_{P}\right) $ as a braided bialgebra. Note that, in this example, $U\left( P,c_{P},b_{P}\right) $ is not
trivial in the sense that it differs from $S\left( P,c_{P}\right).$
\end{example}

In this paper we are not interested in a classification of braided Lie
algebras. Thus, for instance, we have not tried to list all possible
brackets giving a braided Lie algebra structure on the braided vector space $%
\left( V,c\right) $ of Example \ref{ex: gu}. We just intended to provide an
example having both not trivial bracket and not Hecke type braiding. On the
other hand, an attempt of classification would require to rewrite the third
condition in Definition \ref{def: c-Lie alg} in a Jacobi identity type form.
This will be concerned in \cite{AS: Ardi-Stumbo} using results in \cite%
{Braverman-Gaitsgory}. Therein, other examples of braided Lie algebras whose
associated Nichols algebra is quadratic will be investigated.

\section{Pareigis-Lie algebras \label{section: Pareigis}}

In this section we will investigate the relation between our notions of Lie
algebra and universal enveloping algebra and the ones introduced by Pareigis
in \cite{Pareigis-On Lie}. Although these constructions were performed for
braided vector spaces which are in addition objects inside the category of
Yetter-Drinfeld modules, we will not take care of this extra structure.

\begin{claim}
Let $\mathcal{B}_{n}$ be the Artin braid group with generators $\tau
_{i},i\in \left\{ 1,\ldots ,n-1\right\} $ and relations%
\begin{equation*}
\tau _{i}\tau _{i+1}\tau _{i}=\tau _{i+1}\tau _{i}\tau _{i+1}\qquad \text{and%
}\qquad \tau _{i}\tau _{j}=\tau _{j}\tau _{i}\qquad \text{if }\left\vert
i-j\right\vert \geq 2.
\end{equation*}%
Let $\nu _{n}:\mathcal{B}_{n}\rightarrow \mathcal{S}_{n}:\tau \mapsto 
\widetilde{\tau }$ be the canonical quotient homomorphism from the braid
group onto the symmetric group and set $s_{i}:=\widetilde{\tau _{i}}$ for
every $i\in \left\{ 1,\ldots ,n-1\right\} .$

Let $\sigma \in \mathcal{S}_{n}$. The length $l(\sigma )$ of $\sigma $ is
the smallest number $m$ such that there is a decomposition $\sigma
=s_{i_{1}}s_{i_{2}}\cdots s_{i_{m}}$. Such a decomposition will be called
reduced. A permutation $\sigma $ may have several reduced decompositions,
but if $\sigma =s_{i_{1}}s_{i_{2}}\cdots s_{i_{m}}$ is one of them, then the
element $\tau _{i_{1}}\tau _{i_{2}}\cdots \tau _{i_{m}}$ is uniquely
determined in $\mathcal{B}_{n}$.

It is well-known that the $\nu _{n}$ has a canonical section $\iota _{n}:%
\mathcal{S}_{n}\rightarrow \mathcal{B}_{n}$. By definition, 
\begin{equation*}
\iota _{n}(\sigma ):=\tau _{i_{1}}\tau _{i_{2}}\cdots \tau _{i_{m}},
\end{equation*}%
where $\sigma =s_{i_{1}}s_{i_{2}}\cdots s_{i_{m}}$ is a reduced
decomposition of $\sigma $. Note that $\iota _{n}$ is just a map, but 
\begin{equation*}
\iota _{n}(\sigma \tau )=\iota _{n}(\sigma )\iota _{n}(\tau ),
\end{equation*}%
for any $\sigma ,\tau \in \mathcal{S}_{n}$ such that $l(\sigma \tau
)=l(\sigma )+l(\tau ).$

Recall that $\sigma \in \mathcal{S}_{p+q}$ is a $(p,q)$-shuffle if $\sigma
(1)<\dots <\sigma (p)$ and $\sigma (p+1)<\dots <\sigma (p+q).$ The set of $%
(p,q)$-shuffles will be denoted by $\left( p\mid q\right) .$
\end{claim}

\begin{claim}
\label{claim: braid action} Let $(V,c)$ be a braided vector space. For $n\in 
\mathbb{N}^{\ast }$, there is a canonical linear representation $\rho _{n}:%
\mathcal{B}_{n}\rightarrow \mathrm{Aut}_{K}(V^{\otimes n})$ such that $\rho
_{n}(\tau _{i})=c_{i},$ where 
\begin{equation*}
c_{i}=V^{\otimes (i-1)}\otimes c\otimes V^{\otimes (n-i-1)}.
\end{equation*}%
The induced action will be denoted by 
\begin{equation*}
\vartriangleright :\mathcal{B}_{n}\times V^{\otimes n}\rightarrow V^{\otimes
n}.
\end{equation*}%
On generators we have $\tau _{i}\vartriangleright x=c_{i}\left( x\right) $
for each $x\in V^{\otimes n}.$ Note that, since $\iota _{n}$ is not a group
homomorphism, this does not restrict in a natural way to an action of $%
\mathcal{S}_{n}$ on $V^{\otimes n}.$

Let $n\geq 2$ and let $\zeta \in K\backslash \left\{ 0\right\} $. Following 
\cite[Definition 2.3]{Pareigis-On Lie}, we set%
\begin{equation*}
V^{\otimes n}\left( \zeta \right) :=\left\{ x\in V^{\otimes n}\mid \left(
\varphi ^{-1}\tau _{i}^{2}\varphi \right) \vartriangleright x=\zeta ^{2}x,%
\text{ for every }\varphi \in \mathcal{B}_{n},i\in \left\{ 1,\ldots
,n-1\right\} \right\} .
\end{equation*}%
There is an action $\blacktriangleright :\mathcal{S}_{n}\times V^{\otimes
n}\left( \zeta \right) \rightarrow V^{\otimes n}\left( \zeta \right) $
defined on generators by 
\begin{equation}
s_{i}\blacktriangleright x=\zeta ^{-1}\tau _{i}\vartriangleright x.
\label{form: (3) Pareigis}
\end{equation}
\end{claim}

\begin{lemma}
\label{lem Pi}Let $(V,c)$ be a braided vector space. For every $n\geq 2$
consider the map%
\begin{equation*}
\Pi _{\zeta }^{n}:V^{\otimes n}\left( \zeta \right) \rightarrow V^{\otimes
n}\left( \zeta \right) ,\Pi _{\zeta }^{n}\left( x\right) =\sum_{\sigma \in 
\mathcal{S}_{n}}\sigma \blacktriangleright x
\end{equation*}%
for $\zeta \in K\backslash \left\{ 0\right\} .$ If $\zeta $ is a primitive $%
n $-th root of unity then $\Pi _{\zeta }^{n}\left( x\right) \in E_{n}\left(
V,c\right) $ for each $x\in V^{\otimes n}\left( \zeta \right) .$
\end{lemma}

\begin{proof}
We set $T:=T\left( V,c\right) .$ By \cite[Corollary 1.22]{AMS-MM}, for every 
$0\leq i\leq n$ and for every $x\in V^{\otimes n},$ we have $\Delta
_{T}^{i,n-i}\left( x\right) =\sum_{\sigma \in \left( i\mid n-i\right) }\iota
_{n}(\sigma ^{-1})\vartriangleright x. $ Note that, for every $y\in
V^{\otimes n}\left( \zeta \right) $ we have 
\begin{equation*}
\sigma ^{-1}\blacktriangleright y\overset{\text{(\ref{form: (3) Pareigis})}}{%
=}\zeta ^{-l(\sigma ^{-1})}\iota _{n}(\sigma ^{-1})\vartriangleright y=\zeta
^{-l(\sigma )}\iota _{n}(\sigma ^{-1})\vartriangleright y
\end{equation*}%
as $l(\sigma ^{-1})=l(\sigma )$ for every $\sigma \in \mathcal{S}_{n}.$
Therefore, if $x\in V^{\otimes n}\left( \zeta \right) $ we get%
\begin{eqnarray*}
\Delta _{T}^{i,n-i}\left[ \Pi _{\zeta }^{n}\left( x\right) \right]
&=&\sum_{\sigma \in \left( i\mid n-i\right) }\iota _{n}(\sigma
^{-1})\vartriangleright \Pi _{\zeta }^{n}\left( x\right) =\sum_{\sigma \in
\left( i\mid n-i\right) }\sum_{\gamma \in \mathcal{S}_{n}}\iota _{n}(\sigma
^{-1})\vartriangleright \left( \gamma \blacktriangleright x\right) \\
&=&\sum_{\sigma \in \left( i\mid n-i\right) }\sum_{\gamma \in \mathcal{S}%
_{n}}\zeta ^{l(\sigma )}\sigma ^{-1}\blacktriangleright \left( \gamma
\blacktriangleright x\right) =\sum_{\sigma \in \left( i\mid n-i\right)
}\sum_{\gamma \in \mathcal{S}_{n}}\zeta ^{l(\sigma )}\left( \sigma
^{-1}\gamma \right) \blacktriangleright x \\
&=&\sum_{\sigma \in \left( i\mid n-i\right) }\sum_{\gamma \in \mathcal{S}%
_{n}}\zeta ^{l(\sigma )}\gamma \blacktriangleright x=\left( \sum_{\sigma \in
\left( i\mid n-i\right) }\zeta ^{l(\sigma )}\right) \Pi _{\zeta }^{n}\left(
x\right) =\binom{n}{i}_{\zeta }\Pi _{\zeta }^{n}\left( x\right) .
\end{eqnarray*}%
where the last equality holds in view of \cite[1.29]{AMS-MM}. If $\zeta $ is
a primitive $n$-th root of unity, one has that $\binom{n}{i}_{\zeta }=0$
unless $i=0$ or $i=n.$ We have so proved that $\Pi _{\zeta }^{n}\left(
x\right) \in E_{n}\left( V,c\right) $ for each $x\in V^{\otimes n}\left(
\zeta \right) .$

We point out that the calculations above are inspired by part of the proof
of \cite[Theorem 5.3]{Pareigis-On Lie}.
\end{proof}

\begin{definition}
\label{def: Pareigis-Lie alg}(cf. \cite[Definition 4.1]{Pareigis-On Lie}) We
will call a \textbf{Pareigis-Lie algebra} any braided vector space $(V,c)$
together with a $K$-linear map%
\begin{equation*}
\left[ -\right] :\bigoplus\limits_{n\in 
\mathbb{N}
,\zeta \in \mathbb{P}_{n}}V^{\otimes n}\left( \zeta \right) \rightarrow V,
\end{equation*}%
uniquely defined by its restriction $\left[ -\right] _{\zeta
}^{n}:V^{\otimes n}\left( \zeta \right) \rightarrow V$ to $V^{\otimes
n}\left( \zeta \right) ,$ such that the following identities hold:%
\begin{gather}
\left[ \sigma \blacktriangleright x\right] _{\zeta }^{n}=\left[ x\right]
_{\zeta }^{n},\text{ for every }\sigma \in \mathcal{S}_{n},x\in V^{\otimes
n}\left( \zeta \right) ,  \label{form: PL1} \\
\sum_{i=1}^{n+1}\left\{ \left[ -\right] _{-1}^{2}\circ \left( V\otimes \left[
-\right] _{\zeta }^{n}\right) \right\} \left( \left( 1,\ldots ,i\right)
\blacktriangleright x\right) =0,\text{ for every }x\in V^{\otimes n+1}\left(
\zeta \right) ,  \label{form: PL2} \\
\left\{ \left[ -\right] _{-1}^{2}\circ \left( V\otimes \left[ -\right]
_{\zeta }^{n}\right) \right\} \left( x\right) =\sum_{i=1}^{n}\left\{ \left[ -%
\right] _{\zeta }^{n}\circ \left( V^{\otimes i-1}\otimes \left[ -\right]
_{-1}^{2}\otimes V^{\otimes n-i}\right) \right\} \left( \left( \tau
_{i-1}\cdots \tau _{1}\right) \vartriangleright x\right) ,  \label{form: PL3}
\\
\text{for every }x\in V^{\otimes n+1}\left( -1,\zeta \right)  \notag
\end{gather}%
where%
\begin{equation*}
V^{\otimes n+1}\left( -1,\zeta \right) :=\left\{ x\in V\otimes V^{\otimes
n}\left( \zeta \right) \mid \left( 1\otimes \varphi \right)
^{-1}\blacktriangleright \left( \tau _{1}^{2}\vartriangleright \left( \left(
1\otimes \varphi \right) \blacktriangleright x\right) \right) =x,\forall
\varphi \in \mathcal{S}_{n}\right\} .
\end{equation*}%
Given a Pareigis-Lie algebra $\left( V,c,\left[ -\right] \right) $ one can
define (see \cite{Pareigis-On Lie}) \textbf{the universal enveloping algebra}%
\begin{equation*}
U_{P}\left( V,c,\left[ -\right] \right) :=\frac{T\left( V,c\right) }{\left(
\Pi _{\zeta }^{n}\left( z\right) -\left[ z\right] _{\zeta }^{n}\mid n\in 
\mathbb{N}
,\zeta \in \mathbb{P}_{n},z\in V^{\otimes n}\left( \zeta \right) \right) }
\end{equation*}%
not to be mixed with (\ref{form: U(V,c,b)}).
\end{definition}

\begin{remark}
\label{rem: Pareigis-Khar}V. K. Kharchenko drew our attention to \cite[%
Corollary 7.5]{Kharchenko-Multilin} where it is proven there exist at least $%
\left( n-2\right) !$ possibilities to define symmetric Pareigis-Lie algebras
in the case when the underline braided vector space is an object in the
category of Yetter-Drinfeld modules over an abelian group algebra.
\end{remark}

Next result associates a Pareigis-Lie algebra to any braided Lie algebra.

\begin{theorem}
\label{teo: BrLie => PLie}Let $\left( V,c,b\right) $ be a braided Lie
algebra. Then $\left( V,c,\left[ -\right] \right) $ is a Pareigis-Lie
algebra where $\left[ x\right] _{\zeta }^{n}:=b_{n}\Pi _{\zeta }^{n}\left(
x\right) $ for every $x\in V^{\otimes n}\left( \zeta \right) .$ Moreover
there is a canonical braided bialgebra projection $U_{P}\left( V,c,\left[ -%
\right] \right) \rightarrow U\left( V,c,b\right) $ which is the identity map
whenever%
\begin{equation}
\sum_{\zeta \in \mathbb{P}_{n}}\mathrm{Im}\left( \Pi _{\zeta }^{n}\right)
=E_{n}\left( V,c\right) ,\text{ for every }n\in 
\mathbb{N}
\label{form: Pi su}
\end{equation}%
is fulfilled.
\end{theorem}

\begin{proof}
We keep many of the notations used in \cite{Pareigis-On Lie} and apply some
properties proved therein. First, observe that, by Lemma \ref{lem Pi}, $\Pi
_{\zeta }^{n}\left( x\right) $ lies in the domain of $b_{n},$ for every $%
x\in V^{\otimes n}\left( \zeta \right) $ whence $\left[ -\right] _{\zeta
}^{n}$ is well defined. For every $\sigma \in \mathcal{S}_{n},x\in
V^{\otimes n}\left( \zeta \right) ,$ we have%
\begin{eqnarray*}
\left[ \sigma \blacktriangleright x\right] _{\zeta }^{n} &=&b_{n}\Pi _{\zeta
}^{n}\left( \sigma \blacktriangleright x\right) =b_{n}\left[ \sum_{\gamma
\in \mathcal{S}_{n}}\gamma \blacktriangleright \left( \sigma
\blacktriangleright x\right) \right] =b_{n}\left[ \sum_{\gamma \in \mathcal{S%
}_{n}}\left( \gamma \sigma \right) \blacktriangleright x\right] \\
&=&b_{n}\left[ \sum_{\gamma \in \mathcal{S}_{n}}\gamma \blacktriangleright x%
\right] =b_{n}\Pi _{\zeta }^{n}\left( x\right) =\left[ x\right] _{\zeta }^{n}
\end{eqnarray*}%
so that (\ref{form: PL1}) holds. In order to obtain (\ref{form: PL2}) we
adapt the proof of \cite[Theorem 3.4]{Pareigis-On Lie} as follows. For every 
$x\in V^{\otimes n+1}\left( \zeta \right) $ and $\zeta \in \mathbb{P}_{n},$
we have%
\begin{equation*}
\sum_{i=1}^{n+1}\left\{ \left[ -\right] _{-1}^{2}\circ \left( V\otimes \left[
-\right] _{\zeta }^{n}\right) \right\} \left( \left( 1,\ldots ,i\right)
\blacktriangleright x\right) =b_{2}\left( \sum_{i=1}^{n+1}\left( \mathrm{Id}%
_{V^{\otimes 2}}-c\right) \left\{ \left( V\otimes \left[ -\right] _{\zeta
}^{n}\right) \left( \left( 1,\ldots ,i\right) \blacktriangleright x\right)
\right\} \right) .
\end{equation*}%
Let us focus on the argument of $b_{2}$ that we denote by $a.$ Note that $%
a\in E_{2}\left( V,c\right) \ $and observe that%
\begin{equation*}
c\left( V\otimes \left[ -\right] _{\zeta }^{n}\right) \left( z\right) 
\overset{\text{(\ref{form: c-bracket})}}{=}\left( b_{n}\otimes V\right)
c_{V,E_{n}\left( V,c\right) }\left( V\otimes \Pi _{\zeta }^{n}\right) \left(
z\right) =\left( \left[ -\right] _{\zeta }^{n}\otimes V\right) \left\{
\left( \tau _{n}\cdots \tau _{1}\right) \vartriangleright z\right\}
\end{equation*}%
so that%
\begin{equation}
c\left( V\otimes \left[ -\right] _{\zeta }^{n}\right) \left( z\right)
=\left( \left[ -\right] _{\zeta }^{n}\otimes V\right) \left\{ \left( \tau
_{n}\cdots \tau _{1}\right) \vartriangleright z\right\} ,\text{ for every }%
z\in V\otimes V^{\otimes n}\left( \zeta \right) .  \label{form: nat [-]}
\end{equation}%
Note also that $V^{\otimes n+1}\left( \zeta \right) \subseteq V\otimes
V^{\otimes n}\left( \zeta \right) $ (see \cite[Proof of Proposition 3.1]%
{Pareigis-On Lie}). We get 
\begin{eqnarray*}
a &=&\sum_{i=1}^{n+1}\left( V\otimes \left[ -\right] _{\zeta }^{n}\right)
\left( \left( 1,\ldots ,i\right) \blacktriangleright x\right)
-\sum_{i=1}^{n+1}c\left( V\otimes \left[ -\right] _{\zeta }^{n}\right)
\left( \left( 1,\ldots ,i\right) \blacktriangleright x\right) \\
&\overset{\text{(\ref{form: nat [-]})}}{=}&\sum_{i=1}^{n+1}\left( V\otimes %
\left[ -\right] _{\zeta }^{n}\right) \left( \left( 1,\ldots ,i\right)
\blacktriangleright x\right) -\sum_{i=1}^{n+1}\left( \left[ -\right] _{\zeta
}^{n}\otimes V\right) \left\{ \left( \tau _{n}\cdots \tau _{1}\right)
\vartriangleright \left( \left( 1,\ldots ,i\right) \blacktriangleright
x\right) \right\}
\end{eqnarray*}%
\begin{eqnarray*}
&=&\sum_{i=1}^{n+1}\left( V\otimes \left[ -\right] _{\zeta }^{n}\right)
\left( \left( 1,\ldots ,i\right) \blacktriangleright x\right)
-\sum_{i=1}^{n+1}\left( \left[ -\right] _{\zeta }^{n}\otimes V\right)
\left\{ \left( n+1,\ldots ,1\right) \blacktriangleright \left( 1,\ldots
,i\right) \blacktriangleright x\right\} \\
&=&\sum_{i=1}^{n+1}\left( V\otimes \left[ -\right] _{\zeta }^{n}\right)
\left( \left( 1,\ldots ,i\right) \blacktriangleright x\right)
-\sum_{i=1}^{n+1}\left( \left[ -\right] _{\zeta }^{n}\otimes V\right)
\left\{ \left( n+1,\ldots ,i\right) \blacktriangleright x\right\} \\
&=&\sum_{i=1}^{n+1}\left( V\otimes b_{n}\Pi _{\zeta }^{n}\right) \left(
\left( 1,\ldots ,i\right) \blacktriangleright x\right)
-\sum_{i=1}^{n+1}\left( b_{n}\Pi _{\zeta }^{n}\otimes V\right) \left\{
\left( n+1,\ldots ,i\right) \blacktriangleright x\right\}
\end{eqnarray*}%
\begin{eqnarray*}
&=&\left( 
\begin{array}{c}
\left( V\otimes b_{n}\right) \left( \sum_{i=1}^{n+1}\sum_{\sigma \in 
\mathcal{S}_{n}}\left\{ \left( 1\otimes \sigma \right) \blacktriangleright
\left( 1,\ldots ,i\right) \blacktriangleright x\right\} \right) + \\ 
-\left( b_{n}\otimes V\right) \left( \sum_{i=1}^{n+1}\sum_{\sigma \in 
\mathcal{S}_{n}}\left\{ \left( \sigma \otimes 1\right) \blacktriangleright
\left( n+1,\ldots ,i\right) \blacktriangleright x\right\} \right)%
\end{array}%
\right) \\
&=&\left( V\otimes b_{n}\right) \left( \sum_{\sigma \in \mathcal{S}%
_{n+1}}\sigma \blacktriangleright x\right) -\left( b_{n}\otimes V\right)
\left( \sum_{\sigma \in \mathcal{S}_{n+1}}\sigma \blacktriangleright x\right)
\\
&=&\left( V\otimes b_{n}\right) \left( \Pi _{\zeta }^{n+1}\left( x\right)
\right) -\left( b_{n}\otimes V\right) \left( \Pi _{\zeta }^{n+1}\left(
x\right) \right)
\end{eqnarray*}%
Thus%
\begin{equation*}
\sum_{i=1}^{n+1}\left\{ \left[ -\right] _{-1}^{2}\circ \left( V\otimes \left[
-\right] _{\zeta }^{n}\right) \right\} \left( \left( 1,\ldots ,i\right)
\blacktriangleright x\right) =b_{2}\left\{ \left( V\otimes b_{n}\right)
\left( \Pi _{\zeta }^{n+1}\left( x\right) \right) -\left( b_{n}\otimes
V\right) \left( \Pi _{\zeta }^{n+1}\left( x\right) \right) \right\} .
\end{equation*}%
Note that it must be true that $\Pi _{\zeta }^{n+1}\left( x\right) \in
\left( V\otimes E_{n}\left( V\right) \right) \cap \left( E_{n}\left(
V\right) \otimes V\right) $.

In order to prove (\ref{form: PL2}), it suffices to check that 
\begin{equation*}
b_{2}\left\{ \left( V\otimes b_{n}\right) \left( \Pi _{\zeta }^{n+1}\left(
x\right) \right) -\left( b_{n}\otimes V\right) \left( \Pi _{\zeta
}^{n+1}\left( x\right) \right) \right\} =0.
\end{equation*}%
Let $F:=\left( \mathrm{Id}-b\right) \left[ E\left( V,c\right) \right] .$
Since $a\in E_{2}\left( V,c\right) $ and $\Pi _{\zeta }^{n+1}\left( x\right)
\in \left( V\otimes E_{n}\left( V\right) \right) \cap \left( E_{n}\left(
V\right) \otimes V\right) ,$ we get that%
\begin{eqnarray*}
u &:&=b_{2}\left\{ \left( V\otimes b_{n}\right) \left( \Pi _{\zeta
}^{n+1}\left( x\right) \right) -\left( b_{n}\otimes V\right) \left( \Pi
_{\zeta }^{n+1}\left( x\right) \right) \right\} \\
&=&\left( b_{2}-\mathrm{Id}_{V^{\otimes 2}}\right) \left( a\right) + \\
&&+\left( V\otimes \left( b_{n}-\mathrm{Id}_{V^{\otimes n}}\right) \right)
\left( \Pi _{\zeta }^{n+1}\left( x\right) \right) +\left( \left( \mathrm{Id}%
_{V^{\otimes 2}}-b_{n}\right) \otimes V\right) \left( \Pi _{\zeta
}^{n+1}\left( x\right) \right) \\
&\in &F+V\otimes F+F\otimes V
\end{eqnarray*}%
so that $u\in V\cap \left( F\right) =V\cap \ker \left( \pi _{U}\right) $
which is zero by Definition \ref{def: c-Lie alg}. Hence (\ref{form: PL2}) is
proved.

Let us check that (\ref{form: PL3}) is true. We will adapt the proof of \cite%
[Theorem 3.5]{Pareigis-On Lie}. We get 
\begin{eqnarray*}
&&\sum_{i=1}^{n}\left\{ \left[ -\right] _{\zeta }^{n}\circ \left( V^{\otimes
i-1}\otimes \left[ -\right] _{-1}^{2}\otimes V^{\otimes n-i}\right) \right\}
\left( \left( \tau _{i-1}\cdots \tau _{1}\right) \vartriangleright x\right)
\\
&=&b_{n}\left( \sum_{i=1}^{n}\Pi _{\zeta }^{n}\left( V^{\otimes i-1}\otimes 
\left[ -\right] _{-1}^{2}\otimes V^{\otimes n-i}\right) \left( \left( \tau
_{i-1}\cdots \tau _{1}\right) \vartriangleright x\right) \right) \\
&=&b_{n}\left( \sum_{i=1}^{n}\Pi _{\zeta }^{n}\left( V^{\otimes i-1}\otimes
b_{2}\left( \mathrm{Id}_{V^{\otimes 2}}-c\right) \otimes V^{\otimes
n-i}\right) \left( \left( \tau _{i-1}\cdots \tau _{1}\right)
\vartriangleright x\right) \right) \\
&=&b_{n}\left( 
\begin{array}{c}
\sum_{i=1}^{n}\Pi _{\zeta }^{n}\left( V^{\otimes i-1}\otimes b_{2}\otimes
V^{\otimes n-i}\right) \left( \left( \tau _{i-1}\cdots \tau _{1}\right)
\vartriangleright x\right) + \\ 
-\sum_{j=1}^{n}\Pi _{\zeta }^{n}\left( V^{\otimes j-1}\otimes b_{2}\otimes
V^{\otimes n-j}\right) \left( V^{\otimes j-1}\otimes c\otimes V^{\otimes
n-j}\right) \left( \left( \tau _{j-1}\cdots \tau _{1}\right)
\vartriangleright x\right)%
\end{array}%
\right) \\
&=&b_{n}\left( 
\begin{array}{c}
\sum_{i=1}^{n}\sum_{\sigma \in \mathcal{S}_{n}}\sigma \blacktriangleright
\left( V^{\otimes i-1}\otimes b_{2}\otimes V^{\otimes n-i}\right) \left(
\left( \tau _{i-1}\cdots \tau _{1}\right) \vartriangleright x\right) + \\ 
-\sum_{j=1}^{n}\sum_{\sigma \in \mathcal{S}_{n}}\sigma \blacktriangleright
\left( V^{\otimes j-1}\otimes b_{2}\otimes V^{\otimes n-j}\right) \left(
\left( \tau _{j}\tau _{j-1}\cdots \tau _{1}\right) \vartriangleright x\right)%
\end{array}%
\right)
\end{eqnarray*}%
\begin{eqnarray*}
&=&b_{n}\left( 
\begin{array}{c}
\sum_{i=1}^{n}\sum_{\varphi \in \mathrm{Im}(\iota _{n})}\zeta ^{-l\left( 
\widetilde{\varphi }\right) }\varphi \vartriangleright \left( V^{\otimes
i-1}\otimes b_{2}\otimes V^{\otimes n-i}\right) \left\{ \left( \tau
_{i-1}\cdots \tau _{1}\right) \vartriangleright x\right\} + \\ 
-\sum_{j=1}^{n}\sum_{\varphi \in \mathrm{Im}(\iota _{n})}\zeta ^{-l\left( 
\widetilde{\varphi }\right) }\varphi \vartriangleright \left( V^{\otimes
j-1}\otimes b_{2}\otimes V^{\otimes n-i}\right) \left\{ \left( \tau _{j}\tau
_{j-1}\cdots \tau _{1}\right) \vartriangleright x\right\}%
\end{array}%
\right) \\
&\overset{\text{(*)}}{=}&b_{n}\left( 
\begin{array}{c}
\sum_{i=1}^{n}\sum_{\varphi \in \mathrm{Im}(\iota _{n})}\zeta ^{-l\left( 
\widetilde{\varphi }\right) }\left( V^{\otimes \widetilde{\varphi }\left(
i\right) -1}\otimes b_{2}\otimes V^{\otimes n-\widetilde{\varphi }\left(
i\right) }\right) \left\{ \left( \varphi _{\left( i\right) }\tau
_{i-1}\cdots \tau _{1}\right) \vartriangleright x\right\} + \\ 
-\sum_{j=1}^{n}\sum_{\varphi \in \mathrm{Im}(\iota _{n})}\zeta ^{-l\left( 
\widetilde{\varphi }\right) }\left( V^{\otimes \widetilde{\varphi }\left(
j\right) -1}\otimes b_{2}\otimes V^{\otimes n-\widetilde{\varphi }\left(
j\right) }\right) \left\{ \left( \varphi _{\left( j\right) }\tau _{j}\tau
_{j-1}\cdots \tau _{1}\right) \vartriangleright x\right\}%
\end{array}%
\right)
\end{eqnarray*}%
\begin{equation*}
\overset{\text{(**)}}{=}b_{n}\left( 
\begin{array}{c}
\sum_{i=1}^{n}\sum_{\varphi \in \mathrm{Im}(\iota _{n})}\zeta ^{-l\left( 
\widetilde{\varphi }\right) }\left( V^{\otimes \widetilde{\varphi }\left(
i\right) -1}\otimes b_{2}\otimes V^{\otimes n-\widetilde{\varphi }\left(
i\right) }\right) \left\{ \left( \tau _{\widetilde{\varphi }\left( i\right)
-1}\cdots \tau _{1}\left( 1\otimes \varphi \right) \right) \vartriangleright
x\right\} + \\ 
-\sum_{j=1}^{n}\sum_{\varphi \in \mathrm{Im}(\iota _{n})}\zeta ^{-l\left( 
\widetilde{\varphi }\right) }\left( V^{\otimes \widetilde{\varphi }\left(
j\right) -1}\otimes b_{2}\otimes V^{\otimes n-\widetilde{\varphi }\left(
j\right) }\right) \left\{ \left( \tau _{\widetilde{\varphi }\left( j\right)
}\tau _{\widetilde{\varphi }\left( j\right) -1}\cdots \tau _{1}\left(
1\otimes \varphi \right) \right) \vartriangleright x\right\}%
\end{array}%
\right)
\end{equation*}%
\begin{eqnarray*}
&=&b_{n}\left( 
\begin{array}{c}
\left( b_{2}\otimes V^{\otimes n-1}\right) \left( \sum_{\varphi \in \mathrm{%
Im}(\iota _{n})}\zeta ^{-l\left( \widetilde{\varphi }\right) }\left(
1\otimes \varphi \right) \vartriangleright x\right) + \\ 
-\left( V^{\otimes n-1}\otimes b_{2}\right) \left( \sum_{\varphi \in \mathrm{%
Im}(\iota _{n})}\zeta ^{-l\left( \widetilde{\varphi }\right) }\left( \left(
\tau _{n}\tau _{n-1}\cdots \tau _{1}\right) \left( 1\otimes \varphi \right)
\right) \vartriangleright x\right)%
\end{array}%
\right) \\
&=&b_{n}\left( 
\begin{array}{c}
\left( b_{2}\otimes V^{\otimes n-1}\right) \left( \sum_{\sigma \in \mathcal{S%
}_{n}}\left( 1\otimes \sigma \right) \blacktriangleright x\right) + \\ 
-\left( V^{\otimes n-1}\otimes b_{2}\right) \left( \sum_{\sigma \in \mathcal{%
S}_{n}}\left( \tau _{n}\tau _{n-1}\cdots \tau _{1}\right) \vartriangleright
\left( \left( 1\otimes \sigma \right) \blacktriangleright x\right) \right)%
\end{array}%
\right) \\
&=&b_{n}\left( \left( b_{2}\otimes V^{\otimes n-1}\right) \left( V\otimes
\Pi _{\zeta }^{n}\right) \left( x\right) -\left( V^{\otimes n-1}\otimes
b_{2}\right) c_{V,V^{\otimes n}}\left( V\otimes \Pi _{\zeta }^{n}\right)
\left( x\right) \right) \\
&=&b_{n}\left( \left( b_{2}\otimes V^{\otimes n-1}\right) -\left( V^{\otimes
n-1}\otimes b_{2}\right) c_{V,V^{\otimes n}}\right) \left( V\otimes \Pi
_{\zeta }^{n}\right) \left( x\right) ,
\end{eqnarray*}%
where in (*)\ we used that, for every $x^{\prime }\in V^{\otimes i-1}\otimes
E_{2}\left( V,c\right) \otimes V^{\otimes n-i},$%
\begin{equation*}
\varphi \vartriangleright \left( V^{\otimes i-1}\otimes b_{2}\otimes
V^{\otimes n-i}\right) \left( x^{\prime }\right) =\left( V^{\otimes 
\widetilde{\varphi }\left( i\right) -1}\otimes b_{2}\otimes V^{\otimes n-%
\widetilde{\varphi }\left( i\right) }\right) \varphi _{\left( i\right)
}\left( x^{\prime }\right) ,
\end{equation*}%
where $\varphi _{\left( i\right) }\in \mathcal{B}_{n+1}$ is a suitable braid
depending on $\varphi $ and that was defined in \cite[Appendix]{Pareigis-On
Lie}; in (**)\ we applied \cite[Proposition 8.1]{Pareigis-On Lie}.

On the other hand, we have%
\begin{eqnarray*}
&&\left[ -\right] _{-1}^{2}\left( V\otimes \left[ -\right] _{\zeta
}^{n}\right) \left( x\right) \\
&=&b_{2}\left( \mathrm{Id}_{V^{\otimes 2}}-c\right) \left( V\otimes \left[ -%
\right] _{\zeta }^{n}\right) \left( x\right) \\
&=&b_{2}\left\{ \left( V\otimes \left[ -\right] _{\zeta }^{n}\right) \left(
x\right) -c\left( V\otimes \left[ -\right] _{\zeta }^{n}\right) \left(
x\right) \right\} \\
&=&b_{2}\left\{ \left( V\otimes b_{n}\Pi _{\zeta }^{n}\right) \left(
x\right) -c\left( V\otimes b_{n}\Pi _{\zeta }^{n}\right) \left( x\right)
\right\} \\
&\overset{\text{(\ref{form: c-bracket})}}{=}&b_{2}\left\{ \left( V\otimes
b_{n}\right) \left( V\otimes \Pi _{\zeta }^{n}\right) \left( x\right)
-\left( b_{n}\otimes V\right) c_{V,E_{n}\left( V,c\right) }\left( V\otimes
\Pi _{\zeta }^{n}\right) \left( x\right) \right\} \\
&=&b_{2}\left\{ \left( V\otimes b_{n}\right) -\left( b_{n}\otimes V\right)
c_{V,E_{n}\left( V,c\right) }\right\} \left( V\otimes \Pi _{\zeta
}^{n}\right) \left( x\right) .
\end{eqnarray*}%
Let $w:=\left( V\otimes \Pi _{\zeta }^{n}\right) \left( x\right) .$ Hence it
remains to check that 
\begin{equation*}
b_{n}\left( \left( b_{2}\otimes V^{\otimes n-1}\right) -\left( V^{\otimes
n-1}\otimes b_{2}\right) c_{V,V^{\otimes n}}\right) \left( w\right)
=b_{2}\left\{ \left( V\otimes b_{n}\right) -\left( b_{n}\otimes V\right)
c_{V,E_{n}\left( V,c\right) }\right\} \left( w\right) .
\end{equation*}%
We have%
\begin{eqnarray*}
s &:&=\left( 
\begin{array}{c}
b_{n}\left( \left( b_{2}\otimes V^{\otimes n-1}\right) -\left( V^{\otimes
n-1}\otimes b_{2}\right) c_{V,V^{\otimes n}}\right) \left( w\right) + \\ 
-b_{2}\left\{ \left( V\otimes b_{n}\right) -\left( b_{n}\otimes V\right)
c_{V,E_{n}\left( V,c\right) }\right\} \left( w\right)%
\end{array}%
\right) \\
&=&\left( 
\begin{array}{c}
\left( b_{n}-\mathrm{Id}\right) \left( \left( b_{2}\otimes V^{\otimes
n-1}\right) -\left( V^{\otimes n-1}\otimes b_{2}\right) c_{V,V^{\otimes
n}}\right) \left( w\right) + \\ 
+\left( \left( \left( b_{2}-\mathrm{Id}\right) \otimes V^{\otimes
n-1}\right) -\left( V^{\otimes n-1}\otimes \left( b_{2}-\mathrm{Id}\right)
\right) c_{V,V^{\otimes n}}\right) \left( w\right) + \\ 
\left( \mathrm{Id}-b_{2}\right) \left\{ \left( V\otimes b_{n}\right) -\left(
b_{n}\otimes V\right) c_{V,E_{n}\left( V,c\right) }\right\} \left( w\right) +
\\ 
-\left\{ \left( V\otimes \left( b_{n}-\mathrm{Id}\right) \right) -\left(
\left( b_{n}-\mathrm{Id}\right) \otimes V\right) c_{V,E_{n}\left( V,c\right)
}\right\} \left( w\right)%
\end{array}%
\right) .
\end{eqnarray*}%
so that $s\in V\cap \left( F\right) =V\cap \ker \left( \pi _{U}\right) $
which is zero by Definition \ref{def: c-Lie alg}. Hence (\ref{form: PL3}) is
proved.

Thus $\left( V,c,\left[ -\right] \right) $ is a Pareigis-Lie algebra. By the
universal property of its universal enveloping algebra (see \cite[Section 6]%
{Pareigis-On Lie}), there is a unique braided bialgebra homomorphism $%
U_{P}\left( V,c,\left[ -\right] \right) \rightarrow U\left( V,c,b\right) $
that lifts the map $i_{U}:V\rightarrow U\left( V,c,b\right) .$ Assume that (%
\ref{form: Pi su}) is fulfilled. Then%
\begin{eqnarray*}
U\left( V,c,b\right) &=&\frac{T\left( V,c\right) }{\left( \left( \mathrm{Id}%
-b\right) \left\{ E\left( V,c\right) \right\} \right) }=\frac{T\left(
V,c\right) }{\left( \left( \mathrm{Id}-b_{n}\right) \left\{ \sum_{\zeta \in 
\mathbb{P}_{n}}\mathrm{Im}\left( \Pi _{\zeta }^{n}\right) \right\} \mid n\in 
\mathbb{N}\right) } \\
&=&\frac{T\left( V,c\right) }{\left( \left( \mathrm{Id}-b_{n}\right) \left\{
\Pi _{\zeta }^{n}\left( z\right) \right\} \mid n\in \mathbb{N},\zeta \in 
\mathbb{P}_{n},z\in V^{\otimes n}\left( \zeta \right) \right) } \\
&=&\frac{T\left( V,c\right) }{\left( \Pi _{\zeta }^{n}\left( z\right) -\left[
z\right] _{\zeta }^{n}\mid n\in \mathbb{N},\zeta \in \mathbb{P}_{n},z\in
V^{\otimes n}\left( \zeta \right) \right) }=U_{P}\left( V,c,\left[ -\right]
\right) .
\end{eqnarray*}
\end{proof}

\begin{remark}
\label{rem: Pareigig prim}Note that, by \cite[Lemma 5.3.3]{Montgomery}, the
projection $U_{P}\left( V,c,\left[ -\right] \right) \rightarrow U\left(
V,c,b\right) $ in Theorem \ref{teo: BrLie => PLie} becomes the identity
whenever $P\left( U_{P}\left( V,c,\left[ -\right] \right) \right) $
identifies with $V$ through the canonical map $i_{U_{P}}:V\rightarrow
U_{P}\left( V,c,\left[ -\right] \right) $ (in view of Theorem \ref{teo:
magnum}, the converse is true if $\mathrm{sdeg}\left( V,c\right) \leq 1$).
\end{remark}

Theorem \ref{teo: BrLie => PLie} leads to the following problem.

\begin{problem}
To determine under which conditions (\ref{form: Pi su}) is true for a given
braided vector space $(V,c)$. Note that Im$\left( \Pi _{\zeta }^{n}\right) $
is a vector subspace of $V^{\otimes n}\left( \zeta \right) $ so that a
necessary condition is $E_{n}\left( V,c\right) \subseteq \sum_{\zeta \in 
\mathbb{P}_{n}}V^{\otimes n}\left( \zeta \right) $.
\end{problem}

\begin{remark}
Assume $\mathrm{char}\left( K\right) \neq 2$ and let $(V,c)$ be a braided
vector space. Let us show (\ref{form: Pi su}) holds for $n=2$ i.e. that Im$%
\left( \Pi _{-1}^{2}\right) =E_{2}\left( V,c\right) .$ By definition, $%
E_{2}\left( V,c\right) =\left\{ x\in V^{\otimes 2}\mid c\left( x\right)
=-x\right\} $ and 
\begin{eqnarray*}
V^{\otimes 2}\left( -1\right) &=&\left\{ x\in V^{\otimes 2}\mid \left(
\varphi ^{-1}\tau _{1}^{2}\varphi \right) \vartriangleright x=x,\text{ for
every }\varphi \in \mathcal{B}_{2}\right\} \\
&=&\left\{ x\in V^{\otimes 2}\mid \tau _{1}^{2}\vartriangleright x=x\right\}
=\left\{ x\in V^{\otimes 2}\mid c^{2}\left( x\right) =x\right\}.
\end{eqnarray*}%
For each $x\in V^{\otimes 2}\left( -1\right) ,$ we have $\Pi _{-1}^{2}\left(
x\right) =\sum_{\sigma \in \mathcal{S}_{2}}\sigma \blacktriangleright
x=x-\tau _{1}\vartriangleright x=x-c\left( x\right) . $ Let $\gamma
:E_{2}\left( V,c\right) \rightarrow V^{\otimes 2}\left( -1\right) $ be
defined by $\gamma \left( x\right) =x/2.$ Then, for every $x\in E_{2}\left(
V,c\right) ,$ we have $\Pi _{-1}^{2}\gamma \left( x\right) =\frac{1}{2}%
\left( x-c\left( x\right) \right) =x$, where the last equality holds as $%
x\in E_{2}\left( V,c\right) $. Hence Im$\left( \Pi _{-1}^{2}\right)
=E_{2}\left( V,c\right) .$

Note that $\Pi _{-1}^{2}$ is not injective in general. For, if $x\in
V^{\otimes 2}\left( -1\right) $ we have that $\Pi _{-1}^{2}\left( x\right)
\in V^{\otimes 2}\left( -1\right) $ so that it makes sense to compute $\Pi
_{-1}^{2}\Pi _{-1}^{2}\left( x\right) =x-c\left( x\right) -c\left( x\right)
+c^{2}\left( x\right) =2\Pi ^{2}\left( x\right) .$ Thus injectivity of $\Pi
_{-1}^{2}$ implies $\Pi _{-1}^{2}\left( x\right) =2x$ i.e. $c\left( x\right)
=-x$ for every $x\in V^{\otimes 2}\left( -1\right) $ which is not true in
general (for example choose $V$ to be a two dimensional vector space and $c$
to be the canonical flip map defined by $c\left( v\otimes w\right) =w\otimes
v$ for every $v,w\in V$).
\end{remark}

In view of Theorem \ref{teo: BrLie => PLie}, we can now recover two
meaningful examples of Pareigis-Lie algebra.

\begin{corollary}
(cf. \cite[Corollaries 4.2 and 5.4]{Pareigis-On Lie}) Let $\left(
A,m_{A},u_{A},c_{A}\right) $ be a braided algebra. Then $\left( A,c_{A},%
\left[ -\right] \right) $ is a Pareigis-Lie algebra where $\left[ x\right]
_{\zeta }^{n}:=m_{A}^{n-1}\Pi _{\zeta }^{n}\left( x\right) $, for every $%
x\in V^{\otimes n}\left( \zeta \right) .$ Moreover, if $A$ is a connected
braided bialgebra, then the space $P\left( A\right) $ of primitive elements
of $A$ forms a Pareigis-Lie algebra too.
\end{corollary}

\begin{proof}
By Proposition \ref{pro: alg is Lie} $\left( A,c_{A},b_{A}\right) $ is a
braided Lie algebra, where $b_{A}\left( z\right) :=m_{A}^{t-1}\left(
z\right) ,$ for every $z\in E_{t}\left( A,c_{A}\right) .$ By Theorem \ref%
{teo: univ U} also $P\left( A\right) $ carries a Pareigis-Lie algebra
structure which is induced by that of $A.$ We conclude by applying Theorem %
\ref{teo: BrLie => PLie}.
\end{proof}

\noindent \textbf{Acknowledgements.} The author is grateful to V. K.
Kharchenko and A. Masuoka for valuable comments. He also thanks M. Gra\~{n}a
for meaningful remarks.

\appendix

\section{Further results on strongness degree}

In this section we collect some further results on strongness degree of
graded braided bialgebras. In particular we will see that under suitable
assumptions the strongness degree of a graded braided bialgebra has an upper
bound.\medskip\newline
The following technical results are needed to prove Theorem \ref{teo: S[inf]
is strong}.

\begin{lemma}
\label{lem: equality ker Deltas}Let $\left( B,m_{B},u_{B},\Delta
_{B},\varepsilon _{B},c_{B}\right) $ be a graded braided bialgebra. Assume
there exists $n\in 
\mathbb{N}
$ such that $\Delta _{B}^{t,n-t}:B^{n}\rightarrow B^{t}\otimes B^{n-t}$ is
injective for every $0\leq t\leq n.$ Then%
\begin{equation*}
\ker \left( \Delta _{B}^{n,1}\right) =\ker \left( \Delta _{B}^{n-1,2}\right)
=\cdots =\ker \left( \Delta _{B}^{1,n}\right) .
\end{equation*}%
Furthermore, for every $a,b\geq 1$ such that $a+b=n+1,$ we have $%
E_{n+1}\left( B\right) =\ker \left( \Delta _{B}^{a,b}\right) .$
\end{lemma}

\begin{proof}
For every $a,b\geq 1$ such that $a+b=n+1,$ we have%
\begin{equation*}
\ker \left( \Delta _{B}^{a,b}\right) =\ker \left[ \left( \Delta
_{B}^{a-1,1}\otimes B^{b}\right) \Delta _{B}^{a,b}\right] =\ker \left[
\left( B^{a-1}\otimes \Delta _{B}^{1,b}\right) \Delta _{B}^{a-1,b+1}\right]
=\ker \left( \Delta _{B}^{a-1,b+1}\right) .
\end{equation*}%
The last assertion follows by definition of $E_{n+1}\left( B\right) $.
\end{proof}

\begin{proposition}
\label{pro: Delta S[n]}Let $\left( B,m_{B},u_{B},\Delta _{B},\varepsilon
_{B},c_{B}\right) $ be a graded braided bialgebra. Let $n\in 
\mathbb{N}
$ and let $S_{n}:=S^{\left[ n\right] }\left( B\right) .$ For every $t\in 
\mathbb{N}
$, set $I_{t}^{S_{n}}:=\left( E\left( S_{n}\right) \right) \cap \left(
S_{n}\right) ^{t}.$ The following assertions hold.

\begin{enumerate}
\item[$i)$] $I_{t}^{S_{n}}=\left\{ 0\right\} $ for every $0\leq t\leq n+1$.

\item[$ii)$] $I_{n+2}^{S_{n}}=\ker \left( \Delta _{S_{n}}^{a,b}\right) ,$
for every $a,b\geq 1$ such that $a+b=n+2$.

\item[$iii)$] $\Delta _{S_{n}}^{a,b}$ is injective for every $a,b\geq 1$
such that $0\leq a+b\leq n+1$.
\end{enumerate}
\end{proposition}

\begin{proof}
Fix $n\in 
\mathbb{N}
$ and let us prove that $i)$ and $ii)$ follow by $iii)$. Set $A:=S_{n}=S^{%
\left[ n\right] }\left( B\right) .$

$iii)\Rightarrow i)$ We have $I_{0}^{A}:=\left( E\left( A\right) \right)
\cap A^{0}=\left\{ 0\right\} $ and $I_{1}^{A}:=\left( E\left( A\right)
\right) \cap A^{1}=\left\{ 0\right\} .$ From this and iii), since, $\Delta
_{A}^{t,1}\left( I_{t+1}^{A}\right) \subseteq A^{t}\otimes
I_{1}^{A}+I_{t}^{A}\otimes A^{1}$, for every $t\in 
\mathbb{N}
,$ one easily proves that i) is satisfied by induction on $t$.

$iii)\Rightarrow ii)$ Let $a,b\geq 1$ be such that $a+b=n+2.$ Using i), we
obtain $\Delta _{A}^{a,b}\left( I_{n+2}^{A}\right) \subseteq A^{a}\otimes
I_{b}^{A}+I_{a}^{A}\otimes A^{b}=\left\{ 0\right\} $ so that $%
I_{n+2}^{A}\subseteq \ker \left( \Delta _{A}^{a,b}\right) $. By iii) and
Lemma \ref{lem: equality ker Deltas}, we have $E_{n+2}\left( A\right) =\ker
\left( \Delta _{A}^{a,b}\right) $ so that $\ker \left( \Delta
_{A}^{a,b}\right) \subseteq \left( E_{n+2}\left( A\right) \right) \cap
A^{n+2}\subseteq \left( E\left( A\right) \right) \cap A^{n+2}=I_{n+2}^{A}.$
Hence ii) holds.

Let us check that iii) holds, by induction on $n\in 
\mathbb{N}
.$

For $n=0,$ by (\ref{c2}), there is nothing to prove. Let $n\geq 1$ and
assume that iii) is true for $n-1.$ Let us prove it for $n.$ Let $R:=S^{%
\left[ n-1\right] }\left( B\right) .$ Then, by inductive hypothesis and the
first part, we have

\begin{enumerate}
\item[$i\prime )$] $I_{t}^{R}=\left\{ 0\right\} $ for every $0\leq t\leq n$.

\item[$ii\prime )$] $I_{n+1}^{R}=\ker \left( \Delta _{A}^{a,b}\right) ,$ for
every $a,b\geq 1$ such that $a+b=n+1.$

\item[$iii\prime )$] $\Delta _{R}^{a,b}$ is injective for every $a,b\geq 1$
such that $0\leq a+b\leq n.$
\end{enumerate}

For every $i\in 
\mathbb{N}
$, let $\pi _{A}^{i}:R^{i}\rightarrow A^{i}:=R^{i}/I_{i}^{R}$ be the
canonical projection. Recall that $\Delta _{A}^{a,b}$ is uniquely defined by 
$\Delta _{A}^{a,b}\circ \pi _{A}^{a+b}=\left( \pi _{A}^{a}\otimes \pi
_{A}^{b}\right) \circ \Delta _{R}^{a,b}.$ By $i\prime )$, $\pi _{A}^{i}$ is
an isomorphism for every $0\leq i\leq n$. Using this fact and hypothesis $%
iii\prime )$, from the last displayed equality, we deduce that $\Delta
_{A}^{a,b}$ is injective for every $a,b\geq 1$ such that $0\leq a+b\leq n.$
To conclude it remains to check that $\Delta _{A}^{a,b}$ is injective for $%
a,b\geq 1$ such that $a+b=n+1$. By $ii\prime )$, we have that%
\begin{equation*}
A^{n+1}=\frac{R^{n+1}}{I_{n+1}^{R}}=\frac{R^{n+1}}{\ker \left( \Delta
_{R}^{a,b}\right) }
\end{equation*}%
so that $\Delta _{R}^{a,b}$ factors to a unique map $\overline{\Delta
_{R}^{a,b}}:A^{n+1}\rightarrow R^{a}\otimes R^{b}$ such that $\overline{%
\Delta _{R}^{a,b}}\circ \pi _{A}^{n+1}=\Delta _{R}^{a,b}$. Furthermore $%
\overline{\Delta _{R}^{a,b}}$ is injective. We have $\Delta _{A}^{a,b}\circ
\pi _{A}^{n+1}=\left( \pi _{A}^{a}\otimes \pi _{A}^{b}\right) \circ \Delta
_{R}^{a,b}=\left( \pi _{A}^{a}\otimes \pi _{A}^{b}\right) \circ \overline{%
\Delta _{R}^{a,b}}\circ \pi _{A}^{n+1}$ and hence $\Delta _{A}^{a,b}=\left(
\pi _{A}^{a}\otimes \pi _{A}^{b}\right) \circ \overline{\Delta _{R}^{a,b}}.$
Since $\pi _{A}^{a}$ and $\pi _{A}^{b}$ are isomorphisms and $\overline{%
\Delta _{R}^{a,b}}$ is injective, it is clear that also $\Delta _{A}^{a,b}$
is injective.
\end{proof}

\begin{lemma}
\label{lem: lim of gr bialg}Let $K$ be a field and let $\left( \left(
B_{i}\right) _{i\in 
\mathbb{N}
},\left( \xi _{i}^{j}\right) _{i,j\in 
\mathbb{N}
}\right) $ be a direct system of $K$-vector spaces, where, for $i\leq j,$ $%
\xi _{i}^{j}:B_{i}\rightarrow B_{j}.$ Assume that each $B_{i}$ is endowed
with a graded braided bialgebra structure $\left(
B_{i},m_{B_{i}},u_{B_{i}},\Delta _{B_{i}},\varepsilon
_{B_{i}},c_{B_{i}}\right) $ such that $\xi _{i}^{j}$ is a graded braided
bialgebra homomorphism$,$ for every $i,j\in 
\mathbb{N}
.$ Then $\underrightarrow{\lim }B_{i}$ carries a natural graded braided
bialgebra structure that makes it the direct limit of $\left( \left(
B_{i}\right) _{i\in 
\mathbb{N}
},\left( \xi _{i}^{j}\right) _{i,j\in 
\mathbb{N}
}\right) $ as a direct system of graded braided bialgebras. Furthermore the $%
n$-th graded component of $\underrightarrow{\lim }B_{i}$ is $%
\underrightarrow{\lim }B_{i}^{n},$ where $B_{i}^{n}$ denotes the $n$-th
graded component of $B_{i}$.
\end{lemma}

\begin{claim}
\label{claim: type one}Let $\left( B,m_{B},u_{B},\Delta _{B},\varepsilon
_{B},c_{B}\right) $ be a graded braided bialgebra. Then, by (\ref{c1}), (\ref%
{c2}), (\ref{gr1}), (\ref{gr2}) and (\ref{form: delta mult gen}) we have
that $\left( B^{0},m_{B}^{0,0},u_{B}^{0},\Delta _{B}^{0,0},\varepsilon
_{B}^{0},c_{B}^{0,0}\right) $ is a braided bialgebra. Denote by $%
_{B^{0}}^{B^{0}}\mathfrak{M}_{B^{0}}^{B^{0}}$ the category of Hopf bimodules
over $B^{0}.$ By (\ref{form: delta mult gen}), $\left(
B^{1},m_{B}^{0,1},m_{B}^{1,0},\Delta _{B}^{0,1},\Delta _{B}^{1,0}\right) $
is an object in $_{B^{0}}^{B^{0}}\mathfrak{M}_{B^{0}}^{B^{0}}.$ Thus it
makes sense to consider the tensor algebra $T_{B^{0}}\left( B^{1}\right) $
and the cotensor algebra $T_{B^{0}}^{c}\left( B^{1}\right) .$ Both of them
carry a graded braided bialgebra structure arising from their respective
universal properties. Moreover by the universal property of the tensor
algebra, there is a morphism 
\begin{equation*}
F:T_{B^{0}}\left( B^{1}\right) \rightarrow T_{B^{0}}^{c}\left( B^{1}\right)
\end{equation*}%
of graded algebras which is the identity on $B^{0}$ and $B^{1}$
respectively. $F$ is in fact a graded braided algebra homomorphism. Thus the
image of $F$ is a graded braided bialgebra which is called \textbf{the
braided bialgebra of Type one associated to }$B^{0}$\textbf{\ and }$B^{1}$
and denoted by $B^{0}\left[ B^{1}\right] .$ In case when $c_{B}$ is the
canonical flip map, this definition and notation goes back to \cite{Ni}. For
further details, the reader is referred also to \cite[Lemma 6.1 and Theorem
6.8]{AM- Type One}. The proofs there are performed inside braided monoidal
categories, but can be easily adapted to a non categorical framework. This
essentially works because of the definition of graded braided bialgebra
which includes the compatibility of $c_{B}$ with graded components of $B.$
\end{claim}

\begin{theorem}
\label{teo: S[inf] is strong}Let $\left( B,m_{B},u_{B},\Delta
_{B},\varepsilon _{B},c_{B}\right) $ be a graded braided bialgebra. Then

\begin{enumerate}
\item[i)] $S^{\left[ \infty \right] }\left( B\right) $ is a graded braided
bialgebra which is strongly $%
\mathbb{N}
$-graded as a coalgebra (i.e. $\mathrm{sdeg}\left( S^{\left[ \infty \right]
}\left( B\right) \right) =0$).

\item[ii)] The $n$-th graded component of $S^{\left[ \infty \right] }\left(
B\right) $ is given by 
\begin{equation*}
S^{\left[ \infty \right] }\left( B\right) ^{n}=\left\{ 
\begin{tabular}{ll}
$B^{0}$ & $\text{for }n=0,$ \\ 
$S^{\left[ n-1\right] }\left( B\right) ^{n},$ & $\text{for }n\geq 1,$%
\end{tabular}%
\right.
\end{equation*}%
so that $S^{\left[ \infty \right] }\left( B\right) =B^{0}\oplus S^{\left[ 0%
\right] }\left( B\right) ^{1}\oplus S^{\left[ 1\right] }\left( B\right)
^{2}\oplus \cdots .$

\item[iii)] There is a unique coalgebra homomorphism $\psi _{B}:B\rightarrow
T_{B^{0}}^{c}\left( B^{1}\right) =T^{c}$ such that $p_{0}^{T^{c}}\circ \psi
_{B}=p_{0}^{B}$ and $p_{1}^{T^{c}}\circ \psi _{B}=p_{1}^{B},$ where $%
p_{i}^{T^{c}}:T^{c}\rightarrow \left( T^{c}\right) ^{i}$ and $%
p_{i}^{B}:B\rightarrow B^{i}$ are the canonical projections. The map $\psi
_{B}$ is indeed a graded braided bialgebra homomorphism and as graded
braided bialgebras we have $\left( S^{\left[ \infty \right] }\left( B\right)
,\pi _{0}^{\infty }\right) \simeq $Im$\left( \psi _{B}\right) .$
\end{enumerate}
\end{theorem}

\begin{proof}
Set $S:=S^{\left[ \infty \right] }\left( B\right) $ and set $S_{i}:=S^{\left[
i\right] }\left( B\right) $ for each $i\in 
\mathbb{N}
$.

By Lemma \ref{lem: lim of gr bialg}, $S$ carries a natural graded braided
bialgebra structure that makes it the direct limit of $\left( \left(
S_{i}\right) _{i\in 
\mathbb{N}
},\left( \pi _{i}^{j}\right) _{i,j\in 
\mathbb{N}
}\right) $ as a direct system of graded braided bialgebras. Furthermore the $%
n$-th graded component of $S^{\left[ \infty \right] }\left( B\right) $ is $%
S^{n}:=\underrightarrow{\lim }\left( S_{i}\right) ^{n},$ where $\left(
S_{i}\right) ^{n}$ denotes the $n$-th graded component of $S_{i}$.

i) By construction, the family $\left( \Delta _{S}^{a,b}\right) _{a,b\in 
\mathbb{N}
}$ of morphisms $\Delta _{S}^{a,b}:S^{a+b}\rightarrow S^{a}\otimes S^{b}$ is
uniquely determined by $\Delta _{S}^{a,b}\circ \left( \pi _{i}^{\infty
}\right) ^{a+b}=\left( \left( \pi _{i}^{\infty }\right) ^{a}\otimes \left(
\pi _{i}^{\infty }\right) ^{b}\right) \circ \Delta _{S_{i}}^{a,b},$ for
every $i\in 
\mathbb{N}
,$ where $\left( \pi _{i}^{\infty }\right) ^{n}:\left( S_{i}\right)
^{n}\rightarrow S^{n}$ is the structural morphism of $S^{n}=\underrightarrow{%
\lim }\left( S_{i}\right) ^{n}.$

By Proposition \ref{pro: Delta S[n]}, $\Delta _{S_{i}}^{a,b}$ is injective
for every $a,b\geq 1$ such that $0\leq a+b\leq i+1.$ In particular for $%
i=a+b-1$ we get the equality $\Delta _{S}^{a,b}\circ \left( \pi
_{a+b-1}^{\infty }\right) ^{a+b}=\left( \left( \pi _{a+b-1}^{\infty }\right)
^{a}\otimes \left( \pi _{a+b-1}^{\infty }\right) ^{b}\right) \circ \Delta
_{S_{a+b-1}}^{a,b}$ where $\Delta _{S_{a+b-1}}^{a,b}$ is injective. Note
also that the system 
\begin{equation*}
\left( S_{0}\right) ^{t}\overset{\left( \pi _{0}^{1}\right) ^{t}}{%
\rightarrow }\left( S_{1}\right) ^{t}\overset{\left( \pi _{1}^{2}\right) ^{t}%
}{\rightarrow }\left( S_{2}\right) ^{t}\rightarrow \cdots \rightarrow \left(
S_{t-1}\right) ^{t}\overset{\left( \pi _{t-1}^{t}\right) ^{t}}{\rightarrow }%
\left( S_{t}\right) ^{t}\overset{\left( \pi _{t}^{t+1}\right) ^{t}}{\underset%
{\sim }{\rightarrow }}\left( S_{t+1}\right) ^{t}\overset{\left( \pi
_{t+1}^{t+2}\right) ^{t}}{\underset{\sim }{\rightarrow }}\cdots
\end{equation*}%
is stationary as, in view of Proposition \ref{pro: Delta S[n]}, the
projection%
\begin{equation*}
\left( S_{n}\right) ^{t}\overset{\left( \pi _{n}^{n+1}\right) ^{t}}{%
\rightarrow }\left( S_{n+1}\right) ^{t}=\frac{\left( S_{n}\right) ^{t}}{%
I_{t}^{S_{n}}}
\end{equation*}%
is an isomorphism for every $0\leq t\leq n+1.$ This entails that $\left( \pi
_{a+b-1}^{\infty }\right) ^{t}:\left( S_{a+b-1}\right) ^{t}\rightarrow S^{t}$
is an isomorphism for $t\in \left\{ a,b,a+b-1\right\} .$ Therefore $\Delta
_{S}^{a,b}=\left( \left( \pi _{a+b}^{\infty }\right) ^{a}\otimes \left( \pi
_{a+b}^{\infty }\right) ^{b}\right) \circ \Delta _{S_{a+b}}^{a,b}\circ \left[
\left( \pi _{a+b}^{\infty }\right) ^{a+b}\right] ^{-1}$ is injective for
every $a,b\geq 1\ $and hence for every $a,b\in 
\mathbb{N}
$. Since $S$ is a graded braided bialgebra, by Theorem \ref{teo: strongness}%
, it is strongly $%
\mathbb{N}
$-graded as a coalgebra.

ii) For every $t\in 
\mathbb{N}
,$ the morphism $\left( \pi _{t}^{t+1}\right) ^{0}:\left( S_{t}\right)
^{0}\rightarrow \left( S_{t+1}\right) ^{0}$ is an isomorphism so that we can
identify $S^{0}$ with $\left( S_{0}\right) ^{0}=B^{0}.$

The morphism $\left( \pi _{t}^{\infty }\right) ^{t+1}:\left( S_{t}\right)
^{t+1}\rightarrow S^{t+1}$ is an isomorphism so that we can identify $%
S^{t+1} $ with $\left( S_{t}\right) ^{t+1}$ for every $t\in 
\mathbb{N}
.$ In particular we have $S^{1}=\left( S_{0}\right) ^{1}=B^{1}.$

iii) The first part concerning $\psi _{B}$ follows by the universal property
of the cotensor coalgebra. By Theorem \ref{teo: strongness}, the canonical
map $\psi _{S}:S\rightarrow T_{S^{0}}^{c}\left( S^{1}\right)
=T_{B^{0}}^{c}\left( B^{1}\right) $ is injective. For $t=1,2$ we have%
\begin{eqnarray*}
p_{t}^{T_{S^{0}}^{c}\left( S^{1}\right) }\circ \psi _{S}\circ \pi
_{0}^{\infty } &=&p_{t}^{S}\circ \pi _{0}^{\infty }=\left( \pi _{0}^{\infty
}\right) ^{t}\circ p_{t}^{B} \\
&=&\left( \pi _{0}^{\infty }\right) ^{t}\circ p_{t}^{T_{B^{0}}^{c}\left(
B^{1}\right) }\circ \psi _{B}=p_{t}^{T_{S^{0}}^{c}\left( S^{1}\right) }\circ
\psi _{B}.
\end{eqnarray*}%
Since $\pi _{0}^{\infty }$ is a graded braided bialgebra homomorphism, the
universal property of cotensor coalgebra entails that $\psi _{S}\circ \pi
_{0}^{\infty }=\psi _{B}.$ Since $\psi _{S}$ is injective and $\pi
_{0}^{\infty }$ surjective we get that $S^{\left[ \infty \right] }\left(
B\right) =\mathrm{Im}\left( \psi _{B}\right) .$
\end{proof}

\begin{corollary}
\label{coro: type1}Let $\left( B,m_{B},u_{B},\Delta _{B},\varepsilon
_{B},c_{B}\right) $ be a graded braided bialgebra which is strongly $\mathbb{%
N}$-graded as an algebra. Then as graded braided bialgebras we have $S^{%
\left[ \infty \right] }\left( B\right) \simeq B^{0}\left[ B^{1}\right],$ the
braided bialgebra of Type one associated to $B^{0}$ and $B^{1}.$
\end{corollary}

\begin{proof}
Set $S:=S^{\left[ \infty \right] }\left( B\right) .$ By Theorem \ref{teo:
S[inf] is strong}, $S^{\left[ \infty \right] }\left( B\right) $ is a graded
braided bialgebra which is strongly $%
\mathbb{N}
$-graded as a coalgebra so that, by Theorem \ref{teo: strongness}, the
canonical map $\psi _{S}:S\rightarrow T_{S^{0}}^{c}\left( S^{1}\right)
=T_{B^{0}}^{c}\left( B^{1}\right) $ is injective. As a quotient of $B$,
which is strongly $%
\mathbb{N}
$-graded as an algebra, by \cite[Proposition 3.6]{AM- Type One}, $S$ is
strongly $%
\mathbb{N}
$-graded as an algebra too. By \cite[Theorem 3.11]{AM- Type One}, this means
that the canonical map $\varphi _{S}:T_{S^{0}}\left( S^{1}\right)
\rightarrow S$ is surjective. Now, the composition $\psi _{S}\circ \varphi
_{S}:T_{S^{0}}\left( S^{1}\right) \rightarrow T_{S^{0}}^{c}\left(
S^{1}\right) $ is the unique graded braided bialgebra homomorphism that
restricted to $S^{0}$ and $S^{1}$ gives the respective inclusion (compare
with \cite[Theorem 6.8]{AM- Type One}). Hence, its image is, by definition, $%
S^{0}\left[ S^{1}\right] $ that is the braided bialgebra of Type one
associated to $S^{0}$ and $S^{1}.$ We conclude by observing that $%
S^{0}=B^{0} $ and $S^{1}=B^{1}.$
\end{proof}

\begin{remark}
When $B$ is the tensor algebra $T(V,c)$ of a braided vector space $(V,c)$,
the previous result should be compared with \cite[Proposition 3.2]{Masuoka-
pre Nichols}.
\end{remark}

Next aim is to show that under suitable assumptions the strongness degree of
a graded braided bialgebra has an upper bound.

\begin{theorem}
\label{teo: finite strongness}Let $\left( B,m_{B},u_{B},\Delta
_{B},\varepsilon _{B},c_{B}\right) $ be a graded braided bialgebra. The
following assertions are equivalent for $n\in 
\mathbb{N}
$.

\begin{enumerate}
\item $\mathrm{sdeg}\left( B\right) \leq n.$

\item $\pi _{n}^{\infty }:S^{\left[ n\right] }\left( B\right) \rightarrow S^{%
\left[ \infty \right] }\left( B\right) $ is an isomorphism.

\item $S^{\left[ n\right] }\left( B\right) $ is strongly $%
\mathbb{N}
$-graded as a coalgebra.
\end{enumerate}
\end{theorem}

\begin{proof}
$\left( 1\right) \Rightarrow \left( 2\right) $ Let $N:=\mathrm{sdeg}\left(
B\right) .$ By definition $\pi _{i}^{i+1}$ is an isomorphism for every $%
i\geq N.$ In particular $\pi _{i}^{i+1}$ is an isomorphism for every $i\geq
n $ whence $\pi _{n}^{\infty }$ is an isomorphism too.

$\left( 2\right) \Rightarrow \left( 3\right) $ It follows by Theorem \ref%
{teo: S[inf] is strong}.

$\left( 3\right) \Rightarrow \left( 1\right) $ By Theorem \ref{teo:
strongness} $\pi _{n}^{n+1}:S^{\left[ n\right] }\left( B\right) \rightarrow
S^{\left[ n+1\right] }\left( B\right) $ is an isomorphism. Now, the
uniqueness of $S\left( \pi _{i}^{i+1}\right) $ entails $\pi
_{i+1}^{i+2}=S\left( \pi _{i}^{i+1}\right) .$ Thus, for every $i\in 
\mathbb{N}
,$ inductively one proves that $\pi _{i}^{i+1}$ is an isomorphism for every $%
i\geq n.$ Hence the direct system 
\begin{equation*}
B=S^{\left[ 0\right] }\left( B\right) \overset{\pi _{0}^{1}}{\rightarrow }S^{%
\left[ 1\right] }\left( B\right) \overset{\pi _{1}^{2}}{\rightarrow }S^{%
\left[ 3\right] }\left( B\right) \rightarrow \cdots
\end{equation*}%
is stationary and $\mathrm{sdeg}\left( B\right) \leq n.$
\end{proof}

\begin{theorem}
\label{teo: trunked}Let $\left( B,m_{B},u_{B},\Delta _{B},\varepsilon
_{B},c_{B}\right) $ be a graded braided bialgebra which is strongly $%
\mathbb{N}
$-graded as an algebra.

The following assertions are equivalent for $N\geq 1$.

\begin{enumerate}
\item $S^{\left[ \infty \right] }\left( B\right) ^{N}=0.$

\item $S^{\left[ \infty \right] }\left( B\right) ^{n}=0,$ for every $n\geq
N. $

\item $S^{\left[ N-1\right] }\left( B\right) ^{N}=0.$

\item $S^{\left[ N-1\right] }\left( B\right) ^{n}=0,$ for every $n\geq N.$
\end{enumerate}

If one of these conditions if fulfilled, then $\mathrm{sdeg}\left( B\right)
\leq N-1$.
\end{theorem}

\begin{proof}
Since $B$ is strongly $%
\mathbb{N}
$-graded as an algebra, by \cite[Proposition 3.6]{AM- Type One}, each graded
algebra quotient of $B$ is strongly $%
\mathbb{N}
$-graded too. This entails that both $S^{\left[ \infty \right] }\left(
B\right) $ and $S^{\left[ N\right] }\left( B\right) $ are strongly $%
\mathbb{N}
$-graded as an algebras. Clearly, if $A$ is a graded algebra which is
strongly $%
\mathbb{N}
$-graded then $A^{n}=A^{n-N}\cdot A^{N},$ for every $n\geq N.$ Hence we
deduce that $\left( 1\right) \Leftrightarrow \left( 2\right) $ and $\left(
3\right) \Leftrightarrow \left( 4\right) .$

By Theorem \ref{teo: S[inf] is strong}, $S^{\left[ \infty \right] }\left(
B\right) ^{N}=S^{\left[ N-1\right] }\left( B\right) ^{N}$ so that $\left(
1\right) \Leftrightarrow \left( 3\right) .$

Now assume $\left( 4\right) $ and let us prove that $\mathrm{sdeg}\left(
B\right) \leq N$. By Theorem \ref{teo: finite strongness}, it is enough to
prove that $S_{N-1}:=S^{\left[ N-1\right] }\left( B\right) $ is strongly $%
\mathbb{N}
$-graded as a coalgebra. By Proposition \ref{pro: Delta S[n]}, $\Delta
_{S_{N-1}}^{a,b}$ is injective for every $a,b\geq 1$ such that $0\leq
a+b\leq N.$ On the other hand $\left( S_{N-1}\right) ^{n}=0,$ for every $%
n\geq N,$ so that $\Delta _{S_{N-1}}^{a,b}:\left\{ 0\right\} =\left(
S_{N-1}\right) ^{a+b}\rightarrow \left( S_{N-1}\right) ^{a}\otimes \left(
S_{N-1}\right) ^{b}$ is injective for every for every $a,b\geq 1$ such that $%
a+b\geq N.$ In particular we have proved that $\Delta _{S_{N-1}}^{n,1}$ is
injective for every $n\in 
\mathbb{N}
$. By Theorem \ref{teo: strongness}, $S_{N-1}$ is strongly $%
\mathbb{N}
$-graded as a coalgebra.
\end{proof}

\begin{definition}
Recall that the \textbf{Hilbert-Poincar\'{e} sery} of a graded algebra $%
\left( B,m_{B},u_{B}\right) $ is by definition 
\begin{equation*}
\mathfrak{h}\left( B\right) :=\sum_{n\in 
\mathbb{N}
}\dim _{K}\left( B^{n}\right) X^{n}\in K\left[ \left[ X\right] \right] .
\end{equation*}
\end{definition}

\begin{corollary}
\label{coro: f.g. Type1}Let $\left( B,m_{B},u_{B},\Delta _{B},\varepsilon
_{B},c_{B}\right) $ be a graded braided bialgebra which is strongly $%
\mathbb{N}
$-graded as an algebra.

i) Assume that $B^{0}\left[ B^{1}\right] $ as a graded braided bialgebra
divides out $T_{B^{0}}\left( B^{1}\right) $ by relations in degree not
greater then $N.$ Then $\mathrm{sdeg}\left( B\right) \leq N-1.$

ii) If $B^{0}\left[ B^{1}\right] $ is finite dimensional, then $\mathrm{sdeg}%
\left( B\right) \leq \deg \left( \mathfrak{h}\left( B^{0}\left[ B^{1}\right]
\right) \right) -1\leq \dim _{K}B^{0}\left[ B^{1}\right] -1.$
\end{corollary}

\begin{proof}
Set $S:=S^{\left[ \infty \right] }\left( B\right) $ and set $S_{i}:=S^{\left[
i\right] }\left( B\right) $ for each $i\in 
\mathbb{N}
$. Since $B$ is strongly $%
\mathbb{N}
$-graded as an algebra, by Corollary \ref{coro: type1}, we have $B^{0}\left[
B^{1}\right] =S.$

i) By Proposition \ref{pro: Delta S[n]}, $\Delta _{S_{n}}^{a,b}$ is
injective for every $a,b\geq 1$ such that $0\leq a+b\leq n+1$. In particular 
$\Delta _{S_{N-1}}^{a,b}$ is injective for every $a,b\geq 1$ such that $%
0\leq a+b\leq N$. Now, by definition we have 
\begin{equation*}
S_{N}=S\left( S_{N-1}\right) =\frac{S_{N-1}}{\left( E\left( S_{N-1}\right)
\right) }=\frac{S_{N-1}}{\left( \bigoplus\limits_{n\in 
\mathbb{N}
}E_{n}\left( S_{N-1}\right) \right) }=\frac{S_{N-1}}{\left(
\bigoplus\limits_{n\geq N+1}E_{n}\left( S_{N-1}\right) \right) }.
\end{equation*}%
Since $S^{n}=\underrightarrow{\lim }S_{i}^{n}$ and $S$ is defined by
relations in degree not greater then $N$, from $E_{n}\left( S_{N-1}\right)
\subseteq \left( S_{N-1}\right) ^{n}$ we deduce that $\bigoplus\limits_{n%
\geq N+1}E_{n}\left( S_{N-1}\right) =0.$ Hence $S_{N-1}\simeq S_{N},$ so
that $S_{N-1}$ is strongly $%
\mathbb{N}
$-graded as a coalgebra. By Theorem \ref{teo: finite strongness}, $\mathrm{%
sdeg}\left( B\right) \leq N.$

ii) Since $B^{0}\left[ B^{1}\right] $ is finite dimensional, then $\deg
\left( \mathfrak{h}\left( B^{0}\left[ B^{1}\right] \right) \right) $ is
finite and $S^{\left[ \infty \right] }\left( B\right) ^{n}=0,$ for every $%
n\geq \deg \left( \mathfrak{h}\left( B^{0}\left[ B^{1}\right] \right)
\right) .$ Either by the first part or by Theorem \ref{teo: trunked}, we
conclude.
\end{proof}

\end{document}